\documentclass{amsproc}
\usepackage[utf8]{inputenc}
\usepackage{amssymb}
\usepackage{amsxtra}
\usepackage{stmaryrd }
\usepackage{amssymb,amsmath,amsthm, amsfonts, stmaryrd}
\usepackage{color}
\usepackage[all]{xy}
\usepackage{hyperref}

\usepackage{mathtools,leftindex,tensor}
\usepackage[version=4]{mhchem}
\usepackage{appendix}
\usepackage{cancel}
\usepackage[tracking,spacing]{microtype}
\theoremstyle{plain}

\microtypecontext{spacing=nonfrench}

\usepackage{tikz-cd}
\usetikzlibrary{arrows}

\usepackage{tikz}

\newtheorem{theorem}{Theorem}[section]
\newtheorem{lemma}[theorem]{Lemma}
\newtheorem{proposition}[theorem]{Proposition}

\newtheorem{definition}[theorem]{Definition} \theoremstyle{definition}
\newtheorem{example}[theorem]{Example}
\newtheorem{construction}[theorem]{Construction}
\newtheorem{remark}[theorem]{Remark}

\newcommand{\source}{\mathsf{s}}
\newcommand{\target}{{\mathsf{t}}}

\newcommand{\pt}{*}
  \DeclareMathOperator{\GL}{GL}

\DeclareMathOperator{\End}{End}

\DeclareMathOperator{\id}{id}

\DeclareMathOperator{\holim}{holim}

\DeclareMathOperator{\obj}{obj}
\DeclareMathOperator{\ho}{Ho}

\newcommand*{\emptycomment}[1]{}

\DeclareMathOperator{\Aut}{Aut}
\DeclareMathOperator{\pFun}{PsFun} 
\DeclareMathOperator{\picard}{Pic}

\newcommand*{\horiprod}{\mathbin{\cdot_\textup h}}%horizontal multiplication
\newcommand*{\vertprod}{\mathbin{\cdot_\textup v}}%vertical multiplication
\newcommand*{\horicirc}{\mathbin{\circ_\textup h}}%horizontal composition
\newcommand{\huaB}{\mathcal{B}}%{{\mathcal{E}}}%{\mathcal{B}}
\newcommand{\huaS}{\mathcal{S}}
\newcommand{\huaA}{\mathcal{A}}%{{\mathcal{F}}}%{\huaA}

\newcommand{\huaR}{\mathcal{R}}
\newcommand{\huaE}{\mathcal{E}}
\newcommand{\huaF}{\mathcal{F}}
\newcommand{\huaG}{\mathcal{G}}

\newcommand{\huaU}{\mathcal{U}}
\newcommand{\huaV}{\mathcal{V}}
\newcommand{\huaW}{\mathcal{W}}
\newcommand{\huaX}{\mathcal{X}}
\newcommand{\huaY}{\mathcal{Y}}
\newcommand{\huaQ}{\mathcal{Q}}
\newcommand{\huaP}{\mathcal{P}}
\newcommand{\huaC}{{\mathcal{C}}}%{\mathcal{C}}
\newcommand{\huaD}{\mathcal{D}}

\newcommand{\huaH}{\mathcal{H}}
\newcommand{\huaK}{\mathcal{K}}

\newcommand{\huaT}{\mathcal{T}}
\newcommand{\huaZ}{\mathcal{Z}}
\newcommand{\huaM}{\mathcal{M}}
\newcommand{\huaN}{\mathcal{N}}

\newcommand{\cat}{\textbf{\textsf{Cat}}}% categories
\newcommand{\Rep}{\textsf{Rep}}%representation cat
\newcommand{\Vect}{\textsf{Vect}}%cat of v.s.
\newcommand{\tVect}{\textsf{s2}$\huaV ect_k$}%cat of v.s.
\newcommand{\Auto}{\textsf{Auto}} %cat of autofunctors
\newcommand{\gpd}{\textsf{Grpd}} %cat of groupoids

\newcommand{\bgpd}{\textsf{Bibun}} %bicategory of groupoids, bibundles and bibundle maps
\newcommand{\lfun}{\textsf{LaxFun}} % category of lax functors, lax natrual transformation
\newcommand{\pfun}{\textsf{psFun}} % category of  pseudofunctors, strong natrual transformation
\newcommand{\pcolim}{\textsf{pseudocolim}}

\DeclareMathOperator{\Fun}{Fun} %all functors

\newcommand{\sel}{\mathrm{sel}}

\newcommand{\Z}{\ensuremath{\mathbb Z}}

\newcommand{\git}{/\!\!/}
\newcommand{\Mfd}{\mathsf{Mfd}}% category of smooth manifolds
\newcommand{\op}{\mathrm{op}}% opposite
\newcommand{\String}{\textup{String}}% string two-group
\newcommand{\Spin}{\textup{Spin}}% spin group
\newcommand{\Stringor}{\textup{stringor}}% string two-group

\newcommand{\ukm}{\underline{1}}
\newcommand{\sset}{\mathrm{sSet}}

\newcommand{\Mor}{\mathrm{Mor}}

\newcommand{\grp}{\textup{grp}}

\title{2-Equivariant 2-Vector bundles and 2K-theories}
\author{Zhen Huan}
\address{Zhen Huan,  School of Mathematics and Statistics,   Huazhong University of Science and Technology, Hubei 430074, China} \curraddr{}
\email{2019010151@hust.edu.cn}
\date{\today}
\begin{document}

\begin{abstract}

We define 2-vector bundles over a Lie groupoid as pseudofunctors into the bicategory of finite-dimensional super algebras, bimodules, and intertwiners. The resulting 2K-theory, obtained as the Grothendieck completion of the homotopy category, is a category in which ordinary K-theory appears as the endomorphism ring of the trivial object and twisted K-theories as morphisms from the trivial object to twistings. We establish a biequivalence between this pseudofunctor model and the stacky model of Kristel–Ludewig–Waldorf. The theory is extended equivariantly for coherent Lie 2-groups, with explicit computations for abelian and discrete 2-groups that recover the representation rings predicted by Lurie's 2-equivariant elliptic cohomology. For the stringor bundle of Stolz–Teichner, we show that it is the filtered pseudocolimit of finite-dimensional approximants, embedding it into our 2K-theoretic framework. Finally, weak groupoid objects internal to a bicategory yield 2-orbifold 2-vector bundles and their 2K-theory, unifying the ordinary and equivariant settings.

\end{abstract}

\maketitle

\tableofcontents

\section{Introduction}

The development of chromatic homotopy theory over recent decades has fundamentally reshaped our understanding of stable homotopy theory. Initiated by the work of Morava, the pursuit of higher chromatic phenomena has led to a rich interplay between homotopy theory, algebraic geometry, and mathematical physics. Among the most striking predictions to emerge from this program is Lurie's vision of a hierarchy of equivariant cohomology theories, indexed by higher groups, that interpolate between ordinary K-theory and elliptic cohomology \cite{Lurie_elliptic_survey}. A key missing piece in this picture has been a geometric, 2-categorical framework that simultaneously accommodates the full structure of 2-vector bundles, twisted K-theories, and the higher categorical data predicted by Lurie's program.

Existing approaches to 2-vector bundles exhibit a diversity of perspectives that, while each valuable in its own right, has left the underlying unity of the subject obscured.
The homotopy-theoretic model of Baas-Dundas-Rognes \cite{bdr} captures the stable aspects but lacks direct geometric transparency; the stacky model of Kristel-Ludewig-Waldorf \cite{KLW2Vect22} is geometrically explicit, built from cover data and descent, but does not naturally encode twisted K-theories as morphisms or reveal the categorical structure that Lurie's program requires; the Stolz-Teichner stringor bundle \cite{st}\cite{StolzTeichnerSpinor}\cite{KLW22_stringor} has resisted a finite-dimensional geometric interpretation; and Lurie's 2-equivariant elliptic cohomology program has remained a programmatic outline without a concrete geometric realization against which its predictions can be verified.

In this paper, we provide such a framework. We define 2-vector bundles over Lie groupoids as pseudofunctors into the bicategory \tVect \hspace{0.1mm} of finite-dimensional super algebras, bimodules, and intertwiners. This definition is not merely a reformulation: by packaging all local trivializations and transition bimodules into the coherence data of the pseudofunctor, it naturally incorporates twisted K-theory as morphisms from the trivial object to twistings, and ordinary K-theory as endomorphisms of the trivial object. The resulting $2K$-theory, obtained as the Grothendieck completion of the homotopy category, is a category, not merely a group: its objects form a Grothendieck group of 2-vector bundles, while its morphisms are Grothendieck groups of $1$-morphisms, and ring structures appear at the level of endomorphisms. This categorical enrichment is precisely the structure required by Lurie's higher-equivariant vision, and it is the first time such a structure has been realized geometrically.

Our framework delivers four principal contributions.

\textbf{First}, we establish a biequivalence between our pseudofunctor model and the stacky model of Kristel-Ludewig-Waldorf \cite{KLW2Vect22} (Theorem \ref{2models_equiv}). The latter constructs 2-vector bundles via hypercover data and descent; our theorem proves that the two models describe the same geometric objects up to Morita equivalence, bridging the global, functorial perspective with the geometrically explicit, cover-based approach. This comparison is the 2-categorical analogue of the familiar equivalence between vector bundles as locally free sheaves and vector bundles as descent data, and it ensures that our framework carries the full geometric content of the KLW model while offering a more flexible language.

\textbf{Second}, we extend the theory to the equivariant setting. For a Lie groupoid equipped with an action by a coherent Lie 2-group, we define 2-equivariant 2-vector bundles as pseudofunctors from the delooping of the 2-group to the bicategory of 2-vector bundles (Definition \ref{def:2eq2vectgrpd}), and develop the corresponding 2-equivariant $2K$-theory (Definition \ref{def:2equiv_2k}). This provides a 2-categorical framework for equivariant higher geometry that is both flexible and computationally tractable. Crucially, we show that classical equivariant K-theory is recovered as the endomorphism ring of the trivial object in the 2K-theory of the action groupoid (Example \ref{ex:act:eq}), and classical equivariant twisted K-theory is recovered as morphisms from the trivial object to a twisting (Example \ref{ex:act:eq:twisted}). These recoveries demonstrate that our framework is not an abstraction disconnected from established theory, but rather a genuine generalization that contains the classical theory as a special case.

\textbf{Third}, we establish a direct link to the stringor bundle of Stolz-Teichner, as rigorously constructed by Kristel-Ludewig-Waldorf \cite{KLW22_stringor}. We show that the stringor bundle determines a pseudofunctor on its associated \v{C}ech groupoid (Proposition \ref{str_bundle_psfunctor}), and, more substantially, that it is the filtered pseudocolimit of finite-dimensional approximants (Theorem \ref{Thm:embed:stringor}). This result provides a rigorous bicategorical formulation of the heuristic picture—anticipated by Stolz and Teichner \cite{st}\cite{StolzTeichnerSpinor}—that the stringor bundle should arise as a limit of finite-dimensional geometric data. The hyperfinite type III$_1$ structure of the fibre is not an obstacle to approximation; rather, it is precisely the type of data that arises from an AFD approximating system. Consequently, the full stringor bundle is encoded in the $2K$-theory of its \v{C}ech groupoid, embedding this key geometric object into our $2K$-theoretic framework.

\textbf{Fourth}, we take a step toward higher orbifold theory. Using weak groupoid objects internal to a bicategory, following the Segalic pseudofunctor model of Paoli \cite{Paoli2019}, we define 2-orbifold 2-vector bundles and their $2K$-theory (Definitions \ref{def:2orb2vectbl}, \ref{def:2orb_2k}). This framework simultaneously generalizes the ordinary $2K$-theory of Lie groupoids (Section \ref{sect_2vect_2K}) and the equivariant $2K$-theory of Lie 2-group actions (Section  \ref{sect_2eq_2vb_2k}), providing a unified higher-categorical setting for orbifold geometry. We show that when the weak groupoid object comes from a Lie groupoid, the theory reduces to the ordinary theory of Section \ref{sect_2vect_2K}; when it comes from a Lie groupoid equipped with an action by a coherent Lie 2-group, it reduces to the equivariant theory of Section \ref{sect_2eq_2vb_2k}. This unification demonstrates that orbifold $2K$-theory is not an additional structure but rather the native language of the framework.

Beyond these four contributions, our framework provides explicit computational verification of Lurie's predictions. For the 2-group $BA$, where $A$ is an abelian Lie group, we classify the internal equivalence classes in $2\Rep(BA)$ (Lemma \ref{obj_2rep_rep}, Proposition \ref{km_clf}). For $A = U(1)$, the classification recovers the representation ring $\Z[t, t^{-1}]$; for $A = \Z/n$, it recovers $\Z[t]/\langle t^n - 1 \rangle$  (Example \ref{ex_ba_k_obj}). These algebraic structures match the predictions of Lurie's 2-equivariant elliptic cohomology program \cite[Section 5]{Lurie_elliptic_survey} exactly. For a discrete 2-group $G$, we classify strong transformations in $\ho(2\Rep(G))$ (Examples \ref{ps_equiv_k}, \ref{ps_equiv_k_mor}), revealing that morphisms between objects correspond to ordinary or projective super representations of $G$, depending on the cohomology class of the compositor. In the trivial case, the endomorphism ring of the trivial object is exactly the ordinary equivariant K-theory $K_G(\pt)$, i.e., the representation ring $R(G)$. These computations demonstrate that the categorical structure of $2K$-theory encodes not only objects but also $1$-morphisms, consistent with the higher categorical structure anticipated by Lurie.

Our computations are carried out in the non-stable, categorical framework of 2-vector bundles and their Grothendieck completion. The fact that they reproduce the representation rings predicted by Lurie’s stable 2-equivariant elliptic cohomology program provides evidence that the non-stable $2K$-theory is a natural geometric precursor to the stable theory, and suggests a direct link between the categorical representation theory of 2-groups and the elliptic cohomology of higher classifying spaces.

We also construct the K-theory spectrum associated to the representation bicategory $2\Rep(\huaG_\bullet)$ (Appendix \ref{subsect:2k_spectrum}), adapting Osorno's construction \cite{Osorno2012SpectraAGT} to the equivariant setting. This spectrum provides the stable homotopy-theoretic underpinning of the $2K$-theory defined in the body of the paper, and its zeroth space recovers the Grothendieck completion of the homotopy category.

The framework developed here naturally suggests two future directions. The first concerns the transgression of group cohomology classes. Given a 2-group $\Gamma^\alpha$ arising from a 4-class $\alpha$ on a compact Lie group $G$, the higher inertia groupoid $\lfun(B\Z, B\Gamma^\alpha)$ encodes the $S^1$-transgression of $\alpha$. We show that the horizontal composition of $1$-morphisms in this higher inertia groupoid recovers the transgression map explicitly (Theorem \ref{transgression_inertia}), providing a categorical perspective on transgression that we expect to extend to a precise relation between the $\String (G)$-equivariant $2K$-theory of the loop groupoid and the twisted $G$-equivariant Tate K-theory of the free loop space. This connection, which would realize a long-standing heuristic expectation, is the subject of ongoing work. The second direction concerns the analytic and spectral development of 2-orbifold 2-vector bundles, including the construction of a proper equivariant spectrum lifting the Grothendieck completion defined in Section \ref{sect_orb_vb_2K}, which would provide a richer invariant and a natural home for the stringor bundle over suitable loop groupoids.

The paper is organized as follows. Section  \ref{sect_prelim} reviews the necessary preliminaries on super algebras, bimodules, and bicategorical conventions. Section \ref{sect_2vect_2K}  constructs 2-vector bundles over Lie groupoids as pseudofunctors, develops their bicategorical structure, defines the $2K$-theory via Grothendieck completion, and establishes the embedding of ordinary and twisted K-theories as endomorphisms and morphisms of the trivial object. Section \ref{sect_2eq_2vb_2k} extends the theory to the equivariant setting: we define 2-equivariant 2-vector bundles over Lie groupoids for coherent Lie 2-groups, develop the corresponding 2-equivariant $2K$-theory, and recover classical equivariant and twisted equivariant K-theories via action groupoids. Section \ref{bieq_2model} establishes the biequivalence between the pseudofunctor model and the stacky model of Kristel–Ludewig–Waldorf, providing a comparison of the two approaches. Section \ref{sect:stringor} relates the stringor bundle to the $2K$-theoretic framework: we show that it defines a pseudofunctor on a \v{C}ech groupoid and arises as a filtered pseudocolimit of finite-dimensional approximants. Section \ref{sect:comp_lurie} contains explicit computations of $2K$-theory for abelian 2-groups $BA$ and discrete 2-groups $G$, recovering the representation rings predicted by Lurie's program and establishing the correspondence between morphisms and (projective) super representations. Section \ref{future_direct} discusses related work and future directions, including higher transgression and 2-orbifold $2K$-theory. Appendix \ref{subsect:2k_spectrum} constructs the K-theory spectrum associated to the representation bicategory $2\Rep(\huaG_\bullet)$.

\section{Preliminaries} \label{sect_prelim}
\subsection{2-vector spaces} \label{sect_2vect}

In this section we give the necessary background of 2-vector spaces that we need in this paper.

%There are several models of 2-vector spaces, including Baez-Crans 2-vector space \cite{baezcrans_VI}, Kapranov-Voevodsky 2-vector space \cite{KV94_2},  $k$-linear categories, finite-dimensional super algebras over $k$ \cite{KLW2Vect22},  each of which forms a bicategory. A full list of models of 2-vector spaces constructed before 2015 is given in \cite{BDSV2015}. 

We use the bicategory of 2-vector spaces applied in \cite{KLW2Vect22}. It is the symmetric monoidal bicategory \tVect \hspace{0.1mm} of finite-dimensional super algebras over a field $k$, which are exactly finite-dimensional, $\Z/2$-graded, unital, associative algebras over $k$. We will use the notation $A= A_0\oplus A_1$ to denote the graded components. The $1$-morphisms $A\rightarrow B$ in \tVect    \hspace{0.1mm}  
are $\Z/2$-graded, finitely generated $B-A$-bimodules. The $2$-morphisms are parity-preserving bimodule intertwiners. 

The composition of $1$-morphisms is defined by relative tensor product. More explicitly, the composition of the $1$-morphisms $M: A\rightarrow B$ and $N: B\rightarrow C$ is the $C-A$-bimodule $N\otimes_B M$. The identity bimodule of an object $A$ is itself, considered as an $A-A$-bimodule in the obvious way. The associator, left and right unitors are defined from the isomorphisms of bimodules below. 
\begin{equation} \label{bimodule_2mor} % alr order correct
    (K\otimes_C N)\otimes_B M  \xrightarrow{a_{K, N, M}}   K\otimes_C (N\otimes_B M) ; \quad M\otimes_A A \xrightarrow{r_{M}} M; \quad B\otimes_B M
   \xrightarrow{l_{M}} M.
\end{equation} where $M$ is any $B-A$-bimodule, $N$ is any $C-B$-bimodule, and $K$ is any $D-C$-bimodule.

\begin{definition} \label{def_alg_iso}
    An $A-B$-bimodule $\huaM$ is  invertible if it is invertible in the bicategory \tVect, i.e. there exists a $B-A$-bimodule $N$ such that there is an invertible even $A-A$-intertwiner  $M\otimes_{B}N\cong A$ and an invertible even $B-B$-intertwiner $N\otimes_AM \cong B$.

 In the bicategory \tVect, two super algebras $A$, $B$ are said to be internally equivalent (or Morita equivalent) if there exists an invertible $A-B$-bimodule $M$.
\end{definition}

Invertible bimodules play a central role in this paper: they serve as the 1-morphism fibres of 2-vector bundles over Lie groupoids (Definition \ref{def:2vect:obj}). %and they correspond to twistings of K-theory (Section \ref{rel_orb_twisted_k}).

Below we give an explicit construction of   a pseudofunctor from the category $\textsf{s}\huaA lg_k$ of super algebras and  super algebra homomorphisms to the bicategory \tVect. 
\begin{example} \label{ex_2vect_pseudo}
    %This example is preliminary for Example \ref{pseudo_IX}.

    Let $\huaA$ and $\huaB$ denote two super algebras  and $f: \huaA\longrightarrow \huaB$ a super algebra homomorphism. As shown in \cite[(1.27)]{Freed_vienna}, we can associate a $\huaB-\huaA$-bimodule to $f$ whose underlying vector space is $\huaB$ and, for any $b, b'\in \huaB$, any $a\in \huaA$,
    \begin{equation}
        b\cdot b'\cdot a := bb'f(a). 
    \end{equation} We use the symbol $\huaB_f$ to denote this bimodule.
    
    In this way we can define  a pseudofunctor $\huaF$ from the category $\textsf{s}\huaA lg_k$  to the bicategory \tVect.  It sends a super algebra $\huaA$ to itself and sends an algebra homomorphism $f:\huaA\rightarrow \huaB$ to the $\huaB-\huaA$-bimodule $\huaB_{f}$. 

The identity $\id_{\huaA}: \huaA\rightarrow \huaA$ is sent to the $\huaA-\huaA$-bimodule $\huaA$, which is the identity map in $\textsf{s}\huaA lg_k$. Thus, $\huaF^0_{\huaA}$ is the identity $2$-morphism.

For any composition of algebra homomorphisms $ \huaA\xrightarrow{f} \huaB \xrightarrow{g} \huaC$, $$\huaF^2_{g, f}: \huaC_{g}\otimes_{\huaB}\huaB_{f}\longrightarrow \huaC_{g\circ f} $$ is defined by \[ c\otimes_{\huaB} b \mapsto cg(b). \]%given by the left unit constraint $l_{\huaC}$.
It is an invertible $\huaC-\huaA$-intertwiner. %-bimodule structure on $\huaC$ is given by: for any $c'\in\huaC$ and $a'\in \huaA$, \[ c'\cdot c \cdot a' = c'c(g\circ f)(a'). \] 
It is straightforwards to check that $\huaF^2$ and $\huaF^0$ satisfy lax associativity and lax left and right unity. 
%checked. correct.

A bundlization of this example is given in Example \ref{pseudo_IX_pre}.
   
\end{example}

%\comment{add the double category structure later.}

%\begin{remark} \cite[Section 2.1, page 8]{KLW2Vect22}
%    One can show that every super algebra $A$ is dualizable with respect to the symmetric monoidal structure, and that dual objects are provided by opposite super algebras $A^{op}$. Furthermore, a super algebra is fully dualizable if and only if it is semisimple.

%Of particular interest are the invertible objects, i.e., super algebras $A$ for which there exists a super algebra $B$ such that $A\otimes   B = 1$ and $B \otimes A = 1$.  \cite[Theorem 2.1.5]{KLW2Vect22}: The invertible objects in the bicategory of 2-vector spaces  are precisely the central simple super algebras. \end{remark}

The bicategory \tVect \hspace{0.1mm} is symmetric monoidal with respect to the direct sum $\oplus$ and the tensor product $\otimes$ of super algebras and bimodules; the monoidal structure will be used in Section \ref{2K_spectra} to define the K-theory spectrum.

\subsection{Conventions}  \label{sec:conventions}

%In the paper, we use $-\horiprod -$ to denote the horizontal product of two $2$-morphisms in the given bicategory, and use $-\vertprod - $ to denote 
%the vertical product of $2$-morphisms. 
%In addition, we use $- \horicirc -$ to denote the multiplication of two elements in a coherent Lie 2-group in computation.

Throughout the paper, $k$ denotes a field of characteristic zero (typically $\mathbb{R}$ or $\mathbb{C}$). All categories and bicategories are assumed to be essentially small when needed.  Smooth manifolds are assumed to be paracompact and finite-dimensional. Lie 2-groups are assumed to be coherent   unless explicitly stated as strict.
Throughout the paper, we use the  symbols  as listed in Figure \ref{symbol_list}.

\begin{figure} \begin{center} 
\begin{tabular}{|c | c | } 
 \hline 
 \textbf{Symbol}  & \textbf{Meaning} 
\\ \hline
%$(B_0, B_1, B_2)$ & a bicategory with objects $B_0$,\\ & $1$-morphisms $B_1$, and $2$-morphisms $B_2$ \\ \hline 
 $-\horiprod -$  & the horizontal product of two $2$-morphisms  \\
  & in the given bicategory \\
 \hline
 $-\vertprod - $ & the vertical product of two $2$-morphisms \\
  & in the given bicategory  \\ \hline
  $- \horicirc -$ &  the multiplication of two elements  \\
   & in a coherent Lie 2-group \\ \hline
  $-\horicirc -$  & the horizontal composition of $1$-morphisms \\
    &in a bicategory \\ \hline
   $1_{f}$  & the identity $2$-morphism \\
    & at the $1$-morphism $f$ \\
   \hline $\id_x$ & the identity $1$-morphism \\
   & at the object $x$ \\ \hline
   $ \hom_{\huaC}(A, B)$ & the set of morphisms from object $A$ \\
    & to $B$ in a category $\huaC$\\ \hline
       $ \End_{\huaC}(A)$ &  $\hom_{\huaC}(A, A)$\\ \hline
   $\Fun(\huaA, \huaB)$   & the category of functors from $\huaA$ to $\huaB$\\
   \hline
  $ \pfun(-, -) $& the bicategory of pseudofunctors, \\ & strong transformations, and modifications \\ \hline
  % $\smap(-, -)$ & simplicial mapping space \\ \hline
   $\Auto(-)$ & the category of automorphisms of an object\\ & in a bicategory \\ \hline
   $\Aut(-)$ & the group of automorphisms of an object \\ & in a category \\ \hline
\end{tabular} \caption{  The Symbol List }  \label{symbol_list}
\end{center} \end{figure}

\section{2-Vector bundles and 2K-theory} \label{sect_2vect_2K}

\subsection{The bicategory of 2-vector bundles} \label{sss_2stk_2vb}

In this section we introduce the central geometric objects of this paper: 2-vector bundles over a Lie groupoid. We adopt the pseudofunctorial definition: a 2-vector bundle is a pseudofunctor from the Lie groupoid to the bicategory \tVect \hspace{0.1mm} of super algebras, bimodules, and intertwiners. This formulation packages all local trivialization data and transition bimodules into the coherence data of the pseudofunctor, thereby avoiding explicit cover-dependent gluing conditions. We show that these objects, together with their 1- and 2-morphisms, form a bicategory. Their relation to the covering-based models of Kristel–Ludewig–Waldorf \cite{KLW2Vect22} will be established in Section \ref{2modeq:thm}.

\begin{definition}\label{def:2vect:obj}
    
 A super 2-vector bundle over $\huaX_\bullet$ is a pseudofunctor
\[ \huaV: \huaX_\bullet \longrightarrow \text{\tVect} .\]
Explicitly, it consists of the data below:
\begin{itemize}
    \item  For each object $x \in \huaX_0$, a super algebra $\huaV(x) \in \obj(\text{\tVect})$;
    \item  For each arrow $ ( \gamma: x \to y )\in \huaX_1$, a $\huaV(y)-\huaV(x)$-bimodule $\huaV(\gamma)$;
\item  For each composable pair $(\gamma_1, \gamma_2) \in \huaX_2$, an invertible intertwiner
  \[ \mu_{\gamma_1,\gamma_2}: \huaV(\gamma_1) \otimes_{\huaV(\target(\gamma_2))} \huaV(\gamma_2) \longrightarrow \huaV(\gamma_1 \circ \gamma_2); \] 
\item  For each object $x \in \huaX_0$, an invertible intertwiner
  \[ u_x: \huaV(x) \longrightarrow \huaV(\id_x). \]

\item 
These data are required to satisfy the standard pseudofunctor coherence conditions: the lax associativity pentagon and the lax left/right unit triangles.

    \begin{itemize}
        \item \textbf{Lax Associativity}

        For any composable $(\gamma_3, \gamma_2, \gamma_1)\in \huaX_3$, 
    \begin{equation} \label{2vect_obj_pentagon} \xymatrix{(\huaV(\gamma_{3})\otimes_{\huaV(\source(\gamma_3))} \huaV(\gamma_{2}) )\otimes_{\huaV(\source(\gamma_2))} \huaV(\gamma_{1}) \ar@/_1pc/[dd]| {a_{\huaV(\gamma_{3}), \huaV(\gamma_{2}), \huaV(\gamma_{1}) }} \ar@/^1pc/[r]^{\mu_{\gamma_3, \gamma_2} \horiprod \id_{\huaV(\gamma_{1})}} & \huaV(\gamma_{3}\circ \gamma_{2}) \otimes_{\huaV(\source(\gamma_2))} \huaV(\gamma_{1})  
    \ar@/^1pc/[rd]^>>>>>{\mu_{\gamma_3 \circ \gamma_2, \gamma_1}  }
    &\\
            & &  \huaV(\gamma_3\circ \gamma_2\circ \gamma_1) \\
            \huaV(\gamma_{3})\otimes_{\huaA(\source(\gamma_3))} (\huaV(\gamma_{2}) \otimes_{\huaV(\source(\gamma_2))} \huaV(\gamma_{1}))  \ar@/_1pc/[r]_-{ \id_{\huaV(\gamma_{3})} \horiprod \mu_{\gamma_2, \gamma_1} } & 
            \huaV(\gamma_{3})\otimes_{\huaV(\source(\gamma_3))} \huaV(\gamma_{2}\circ \gamma_{1}) \ar@/_1pc/[ru]_-{ \mu_{\gamma_3, \gamma_2\circ \gamma_1} }
              } 
            \end{equation}

\item \textbf{Lax Left and Right Unity}

For any $\gamma \in \huaX_1$, 
            \begin{equation} \label{2vect_obj_unity}
                \xymatrix{  \huaV(\target(\gamma))\otimes_{\huaV(\target(\gamma))} \huaV(\gamma) \ar[r]^-{l_{\huaV(\gamma)}}
                \ar[d]_-{u_{\target(\gamma)} \horiprod 1_{\huaV(\gamma)}  }   & \huaV(\gamma) \\
                \huaV(\id_{\target(\gamma)}) \otimes_{\huaV(\target(\gamma))} \huaV(\gamma) \ar[ru]_-{\mu_{\id_{\target(\gamma)}, \gamma} }  }
                \quad 
                 \xymatrix{ \huaV(\gamma)\otimes_{\huaV(\source(\gamma))} \huaV(\source(\gamma))  \ar[r]^-{r_{\huaV(\gamma)}}
                \ar[d]_-{ 1_{\huaV(\gamma)} \horiprod u_{\source(\gamma)}  }   & \huaV(\gamma) \\
              \huaV(\gamma)   \otimes_{\huaV(\source(\gamma))}   \huaV(\id_{\source(\gamma)}) \ar[ru]_-{\mu_{ \gamma,  \id_{\source(\gamma)}} }  }
            \end{equation}
            \end{itemize}
\end{itemize}
\end{definition}

In this paper, all pseudofunctors are assumed to be smooth in the sense of Remark \ref{rmk:smooth}; the smoothness conditions will be made explicit there. The algebraic core of the construction does not depend on this choice. 

\begin{remark}
    The assignment $\gamma \mapsto \huaV(\gamma)$ in Definition \ref{def:2vect:obj} is part of a functor on the underlying groupoid; the 2-morphism level is trivial because the base groupoid has only identity 2-morphisms. Equivalently, the coherence data are already included in the definition of a 2-vector bundle.
\end{remark}

%\cite[Definition 2.6]{Nikolaus-Schweigert}

\begin{definition}\label{def:2vect:1mor}
Let $ \huaV $ and $ \huaV'  $ denote two super 2-vector bundles over $\huaX_{\bullet}$.
A 1-morphism \[\delta: \huaV\rightarrow \huaV'\] is a strong transformation which consists of: 
\begin{itemize}

    \item For each $x\in \huaX_{0}$, a $\huaV'(x)-\huaV(x)$-bimodule $\delta_x$;
    % functor on the level of objects

    \item For any $x\xrightarrow{\gamma} y$ in $\huaX_1$, an invertible $\huaV'(y)-\huaV(x)$-intertwiner $$\delta_{\gamma}:   
    \huaV'(\gamma)\otimes_{\huaV' (x)}    
    \delta_x   \longrightarrow      \delta_y \otimes_{\huaV(y) } \huaV(\gamma)$$
 \begin{equation} \label{def:1mor_2vb_phi}
         \xymatrix{ \huaV(x)  \ar[d]_{\huaV(\gamma)} \ar[r]^{\delta_x} &  \huaV'(x)   \ar[d]^{\huaV'(\gamma)} \\
    \huaV(y)  \ar@{}@<+0.5ex>[ru]|{\buildrel{\delta_\gamma}\over\Leftarrow} \ar[r]_{\delta_y} &  \huaV'(y) . } 
    \end{equation}

    \item The data above make diagrams below commute.

    \begin{itemize}
        \item \textbf{Lax Naturality}

            For any $(\gamma', \gamma)\in \huaX_2$, 
    \begin{equation} \label{equiv-obj-1M3}
         \xymatrix{ \huaV(\target(\gamma'))      \ar[r]^{\delta_{\target(\gamma')}} &  \huaV'(\target(\gamma'))   & \\
          &\huaV(\source(\gamma')   )    \ar[lu]|{\huaV(\gamma')}    \ar@{}@<+0.5ex>[l]|>>>>>>>>{\buildrel{\mu_{\gamma', \gamma}}\over\Leftarrow}   \ar@{}@<+0.5ex>[u]|{\buildrel{\delta_{\gamma'}}\over\Leftarrow}   \ar[r]^{\delta_{\source(\gamma')}} &  \huaV'(\source(\gamma') ) =    \ar[lu]_{\huaV'(\gamma')}   \\
    \huaV(\source(\gamma))  \ar[ru]|{\huaV(\gamma)}    \ar[uu]^{\huaV(\gamma'\circ \gamma)}   \ar@{}@<+0.5ex>[rru]|{\buildrel{\delta_{\gamma}}\over\Leftarrow} \ar[r]_{\delta_{\source(\gamma)}}  & \huaV'( \source(\gamma))  \ar[ru]_{\huaV'(\gamma)}    &} 
             \xymatrix{ \huaV( \target(\gamma'))   \ar[r]^{\delta_{\target(\gamma')}} &  \huaV'(\target(\gamma'))    & \\
          &  &  \huaV'(\source(\gamma'))    \ar@{}@<+0.5ex>[l]|>>>>>>>>>>{\buildrel{\mu'_{\gamma', \gamma}}\over\Leftarrow}  
          \ar[lu]_{\huaV'(\gamma')}   \\
    \huaV( \source(\gamma))   \ar[uu]|{\huaV(\gamma'\circ \gamma)}  \ar@{}@<+0.5ex>[ruu]|{\buildrel{\delta_{\gamma'\circ \gamma}}\over\Leftarrow} \ar[r]_{\delta_{\source(\gamma)}} 
    & \huaV'( \source(\gamma) )  \ar[ru]_{\huaV'(\gamma)}  \ar[uu]|{\huaV'(\gamma'\circ \gamma)}    & } 
    \end{equation}

\bigskip

    \item \textbf{Lax Unity}

    For any $x\in \huaX_0$,

    \begin{equation} \label{def:1mor_2vectbd_unity}
       \xymatrix{ \huaV(x) \ar[dd]_-{\huaV(\id_x)} \ar[r]^{\delta_x} &   \huaV'(x)  \ar[dd]|{\huaV'(\id_x)}  \ar@/^4pc/[dd]^{\huaV'(x)}  & \\ & & \\
    \huaV(x)   \ar@/_1pc/[r]_-{\delta_x} 
      \ar@{}@<+0.5ex>[ruu]|{\buildrel{\delta_{\id_x}}\over \Leftarrow} &  \huaV'(x)  \ar@{}@<+0.5ex>[ruu]|{\buildrel{u'_{x}}\over\Leftarrow} &\\ } =
    \xymatrix{ &\huaV(x)\ar[rdd]|{\delta_x} \ar@/_4pc/[dd]_{\huaV(\id_x)}  \ar[r]^-{\delta_x} \ar[dd]|{\huaV(x)} &   \huaV'(x) \ar[dd]|{\huaV'(x)}    \\ &  \ar@{}@<+0.5ex>[ru]|{\buildrel{l_{\delta_x}}\over\Leftarrow} & \\
    &\huaV(x)   \ar@{}@<+0.5ex>[luu]|{\buildrel{u_{x}}\over\Leftarrow}  
    \ar@/_1pc/[r]_-{\delta_x}
      \ar@{}@<+0.5ex>[ru]|{\buildrel{r^{-1}_{\delta_x}}\over\Leftarrow} &  \huaV'(x)   \\ } .
\end{equation}  

\end{itemize}
   % \begin{equation}
    % \xymatrix{ (\huaM')_{g'_1(z'')} \otimes_{ \huaA'_{g'(t(z'))}}(\huaM')_{g'_1(z')} \otimes_{ \huaA'_{g'(t(z))}}  \huaP_{z} \ar[d]_-{\id\otimes \phi_{z', z}} \ar[rr]^-{(\mu_{2 })_{g'(z''), g'(z') }\otimes \id}   &&(\huaM')_{g'_1(z''\circ z')} \otimes_{ \huaA'_{g'(t(z))}}  \huaP_{z} \ar[dd]^-{\phi_{z'', z}} \\     (\huaM')_{g'_1(z'')} \otimes_{ \huaA'_{g'(t(z'))}} P_{z'} \otimes_{\huaA_{g(s(z'))}} (\huaM)_{g_1(z)}  \ar[d]_-{\phi_{z'', z'}\otimes \id}  &&\\   P_{z''} \otimes_{\huaA_{g(s(z''))}} (\huaM)_{g_1(z')} \otimes_{\huaA_{g(s(z'))}} (\huaM)_{g_1(z)} \ar[rr]_-{\id\otimes (\mu_1)_{g(z'), g(z) }} & &  P_{z''} \otimes_{\huaA_{g(s(z''))}} (\huaM)_{g_1(z'\circ z)}} \end{equation} for all $(z'', z', z)\in \huaZ_3$.
\end{itemize}

\end{definition}

\begin{remark}
    The composition of $1$-morphisms between 2-vector bundles is exactly the horizondal composition of strong transformation. An explicit formula of it is given \cite[Definition 4.2.15]{JY:Bicat}. In addition, the identity $1$-morphisms $1_{\huaV}$ at a 2-vector bundle $\huaV$ is exactly the identity transformation defined in \cite[Definition 4.2.14]{JY:Bicat}.
\end{remark}

%\begin{remark} Here is a little explanation of the bimodule bundles $ \partial^*_1 \huaP$ and $ \partial^*_0 \huaP$ over $\huaZ_1$. Each $ \partial^*_i \huaP$ are the pull-back bundle along  \[ \huaZ_1\xrightarrow{d_i} \huaZ_0 \hookrightarrow \huaZ_1, \text{ for } i=1, 2, \] where the second functor sends each object $z_0$ to the identity map $\id_{z_0}$.   For each $z\in \huaZ_1$, the fiber $ (\partial^*_1 \huaP)_{z}$ is $\huaP_{\id_{s(z)}}$, and the fiber $( \partial^*_0 \huaP)_{z}$ is $\huaP_{\id_{t(z)}}$.  \end{remark}

%%\cite[Definition 2.6]{Nikolaus-Schweigert}
\begin{definition}\label{def:2mor_2vectbd}

Let \[\delta, \delta': \huaV\longrightarrow \huaV' \] denote two 1-morphisms  between the  super 2-vector bundles $ \huaV $ and $ \huaV'  $.

A $2$-morphism $\kappa: \delta\Rightarrow \delta'$ consists of 
an intertwiner \[\kappa_x: \delta_x\longrightarrow \delta'_x,\] for each $x\in \huaX_0$, such that for each $\gamma\in \huaX_1$
\begin{equation} \label{def:2mor_2vectbd_eq}
       \xymatrix{ &\huaV(\source(\gamma))         \ar[r]^-{\huaV_{\gamma}} \ar[dd]|{\delta_{\source(\gamma)}} &   \huaV(\target(\gamma))  \ar[dd]|{\delta_{\target(\gamma)}}      \ar@/^4pc/[dd]^{\delta'_{\target(\gamma)}}   \ar@{}@<+6ex>[dd]|{\buildrel{\kappa_{\target(\gamma)}}\over\Rightarrow} \\ & & \\
   & \huaV'(\source(\gamma))  
   \ar@/_1pc/[r]_-{ \huaV'_{\gamma}} 
      \ar@{}@<+0.5ex>[ruu]|{\buildrel{\delta_{\gamma}}\over\Rightarrow}      &  \huaV'(\target(\gamma)) 
    \\ } =
    \xymatrix{ \huaV(\source(\gamma))     \ar@/_5pc/[dd]_{\delta_{\source(\gamma)}}       \ar[r]^-{\huaV_{\gamma}} \ar[dd]|{\delta'_{\source(\gamma)}} &  \huaV(\target(\gamma))    \ar[dd]|{\delta'_{\target(\gamma)}}  & \\ & & \\
    \huaV'(\source(\gamma))          \ar@{}@<+7ex>[uu]|{\buildrel{\kappa_{\source(\gamma)}}\over\Rightarrow}     \ar@/_1pc/[r]_-{\huaV'_{\gamma}} 
      \ar@{}@<+0.5ex>[ruu]|{\buildrel{\delta'_{\gamma}}\over\Rightarrow}   & \huaV'(\target(\gamma))   &\\ } .
\end{equation}

\end{definition}

\begin{remark} \label{2mor_vert_comp}

The vertical composition and the horizontal composition of two $2$-morphisms between $1$-morphisms between $2$-vector bundles are exactly the vertical composition and the horizontal composition of modifications, whose explicit formulas are given in \cite[Definition 4.4.5]{JY:Bicat}. And the identity $2$-morphism at a $1$-morphism is the identity modification defined in \cite[Definition 4.4.5]{JY:Bicat}.

\end{remark}

\begin{definition}
    The super 2-vector bundles over $\huaX_{\bullet}$, the $1$-morphisms between super 2-vector bundles over $\huaX_{\bullet}$, and the $2$-morphisms between $1$-morphisms, form the bicategory 
    $\textsf{s}2\huaV \huaB dl_k (\huaX_{\bullet}) $   of super 2-vector bundles over $\huaX_{\bullet}$.
\end{definition}

\begin{example} \label{pseudo_IX_pre}
Let $\huaX_{\bullet}$ be a Lie groupoid. Let $\textsf{s}\huaA lg\huaB dl_k(\huaX_\bullet)$ denote the category of super algebra bundles over $\huaX_{\bullet}$ and even algebra bundle homomorphisms. We can define a pseudofunctor $\tilde{\huaF}_{\huaX_{\bullet}}: \textsf{s}\huaA lg\huaB dl_k(\huaX_\bullet) \longrightarrow \textsf{s}2\huaV \huaB dl_k(\huaX_{\bullet}) $ by composing the pseudofunctor $\huaF:   \textsf{s}\huaA lg_k\rightarrow  \text{\tVect}$ defined in Example \ref{ex_2vect_pseudo}, i.e. sending a  functor $\huaX_{\bullet}\xrightarrow{g}  \textsf{s}\huaA lg_k$ to the composition
\[ \huaX_{\bullet}\xrightarrow{g}  \textsf{s}\huaA lg_k \xrightarrow{\huaF} \text{\tVect}.\] This is the bundlization of the algebraic embedding.

\end{example}

\begin{theorem} \label{biequiv_2vectbdl}
\leavevmode
\begin{enumerate}
    \item Each Lie homomorphism (i.e. strict functor) $F: \huaY_{\bullet} \rightarrow \huaX_{\bullet}$ induces a restriction pseudofunctor 
    \[ F^{\ast}: \textsf{s}2\huaV \huaB dl_k(\huaX_{\bullet}) \longrightarrow \textsf{s}2\huaV \huaB dl_k(\huaY_{\bullet})\] by precomposing each pseudofunctor $\huaX_{\bullet}\to \text{\tVect}$ by $F$.
    \item 
    If $\huaX_{\bullet}$ and $\huaY_{\bullet}$ are biequivalent in the bicategory $\bgpd$ of Lie groupoids, bibundles and equivariant bibundle maps, then the bicategories   $\textsf{s}2\huaV \huaB dl_k (\huaX_{\bullet}) $ and  $\textsf{s}2\huaV \huaB dl_k (\huaY_{\bullet}) $ are biequivalent.
    \end{enumerate}
\end{theorem}

The theorem can be proved by definitions routinely. %The proof of Theorem \ref{biequiv_2vectbdl} (ii) is a generalization of that of \cite[Theorem 3.1]{Huan:2Rep_2Vect}. 
We sketch the proof below.

\begin{proof}
\leavevmode

\begin{enumerate}
    \item 
  Consider the precomposition operation \[ F^*: \pfun(\huaX_\bullet, \text{\tVect}) \longrightarrow \pfun(\huaY_\bullet, \text{\tVect}), \qquad
\huaV \longmapsto \huaV \circ F.\]

We verify the pseudofunctoriality level by level:
\begin{enumerate}
    \item 
On objects: If $\huaV: \huaX_\bullet \to \text{\tVect}$ is a pseudofunctor, then the composite $\huaV \circ F: \huaY_\bullet \to \text{\tVect}$ is again a pseudofunctor, since pseudofunctors are closed under composition.
\item  On 1-morphisms: If $\alpha: \huaV \Rightarrow \huaW$ is a strong transformation, then $F^*(\alpha) := \alpha \horicirc \id_F$ (i.e., the pre-whiskering of $\alpha$ with $F$) is a strong transformation $\huaV\circ F \Rightarrow \huaW\circ F$. 
\item  On 2-morphisms: If $\Gamma: \alpha \Rightarrow \beta$ is a modification, then $F^*(\Gamma) := \Gamma \horiprod 1_{\id_F}$ is a modification $\alpha \horicirc \id_F \Rightarrow \beta  \horicirc \id_F.$
\item  Preservation of composition and units: Precomposition strictly preserves  vertical composition of strong transformations  as well as identity transformations, by the middle four exchange law. Hence $F^*$ is strict on 1- and 2-morphism composition (i.e., $F^*(\alpha \circ \beta) = F^*(\alpha) \circ F^*(\beta)$, and $F^*(\id_\huaV) = \id_{F^*(\huaV)}).$
\end{enumerate} 
Therefore, $F^*$ is a pseudofunctor from $\pfun(\huaX_\bullet, \text{\tVect})$ to $\pfun(\huaY_\bullet, \text{\tVect})$. Since $\pfun(\huaX_\bullet, \text{\tVect}) = \textsf{s}2\huaV \huaB dl_k(\huaX_\bullet)$, the claim follows. 
    
    \item 
    Consider the pseudofunctor \[ \pfun(-, \text{\tVect}): \bgpd^{\op} \longrightarrow \textbf{Bi}\cat\] which sends each Lie groupoid $\Gamma_{\bullet}$ to the bicategory $\pfun(\Gamma_{\bullet}, \text{\tVect})$ of pseudofunctors, and sends each bibundle to the pullback pseudofunctor given by precomposition. The 2-functoriality at the level of $1$-morphisms follows from the standard properties of pseudofunctor composition, and at the level of $2$-morphisms form whiskering. The coherence conditions are satisfied by \cite[Definition 4.1.26, Lemma 4.1.29]{JY:Bicat}.
    By the assumed biequivalence, there exist $1$-morphisms \[ F: \huaX_{\bullet}\rightarrow \huaY_{\bullet}, \quad G: \huaY_{\bullet}\rightarrow \huaX_{\bullet}\] in $\bgpd$, together with invertible $2$-morphisms 
    \[ \eta: \id_{\huaX_{\bullet}} \Rightarrow G\circ F, \quad \epsilon: F\circ G\Rightarrow \id_{\huaY_{\bullet}}.\]
    Applying the pseudofunctor $ \pfun(-, \text{\tVect})$ to these data yields:
    \begin{itemize}
        \item $1$-morphisms:
        \begin{align*}
            F^*:  \pfun(\huaY_{\bullet}, \text{\tVect}) &\longrightarrow  \pfun( \huaX_{\bullet}, \text{\tVect}), \\
            G^*:  \pfun(\huaX_{\bullet}, \text{\tVect}) &\longrightarrow  \pfun( \huaY_{\bullet}, \text{\tVect}),
        \end{align*}
        \item invertible $2$-morphisms
            \begin{align*}
                &\tilde{\eta}: \quad \id_{\pfun(\huaX_{\bullet}, \text{\tVect})}  \Rightarrow (G\circ F)^*, \\
                &\tilde{\epsilon}: \quad (F\circ G)^*  \Rightarrow \id_{\pfun(\huaY_{\bullet}, \text{\tVect})} 
            \end{align*} obtained by post-whiskering $\eta$ and $\epsilon$, respectively.
  
    \end{itemize}

    By the $2$-functoriality of pseudofunctors, the $2$-morphisms
    $\tilde{\eta}$ and $\tilde{\epsilon}$ inherit the triangle identities from $\eta$ and $\epsilon$. Thus, 
\[ (F^*, G^*, \tilde{\eta}, \tilde{\epsilon})\] forms a biequivalence between $\pfun(\huaY_{\bullet}, \text{\tVect}) $  and $  \pfun( \huaX_{\bullet}, \text{\tVect})$. 
    Thus, the result follows. 
    \end{enumerate}
    
\end{proof}

\subsection{2K-theory and its spectra} \label{2K_spectra}

In this section we  define the K-theory spectrum of \tVect .\mbox{ } Later in Section \ref{sect:def_2K}, we define the corresponding $2K$-theory. The principle of the construction is based on \cite{Osorno2012SpectraAGT}, which extends the idea of the definition of the $2K$-theory in \cite{bdr}
 to a more general setting. %In \cite{bdr} the 2-vector space model is taken to be that by Kapranov and Voevodsky \cite{KV94_2}. 

Recall we have the properties for the operations $\oplus$ and $\otimes$ in \tVect \mbox{ }below.

For any $B-A$-bimodule $\huaM$, $B'-A'$-bimodule $\huaM'$, $C-B$-bimodule $\huaN$, and $C'-B'$-bimodule $\huaN'$, 
the direct sum $\huaM\oplus \huaM'$ is a $(B\oplus B')- (A\oplus A')$-bimodule, and the tensor product of the direct sums of bimodules satisfies 
\begin{equation} \label{bimod_distr}
    (\huaN\oplus \huaN')\otimes_{B\oplus B'} (\huaM\oplus \huaM') = (\huaN\otimes_{B}\huaM)\oplus (\huaN'\otimes_{B'}\huaM').
\end{equation} %\comment{Note, this is the unique well-defined way. Thus, the $F^2$ of $\oplus$ can be identity.}

\begin{lemma} The bicategory \tVect  \hspace{0.2mm}  is a strict symmetric monoidal bicategory in the sense of \cite[Definition 1.3]{Osorno2012SpectraAGT}. \label{2vect_strict}\end{lemma}

\begin{proof}
The direct sum \[ \oplus: \text{\tVect} \times \text{\tVect} \longrightarrow \text{\tVect}\] is a strict functor as defined below.
\begin{itemize}
    \item on the level of objects: $(A, B)$ is sent to the direct sum $A\oplus B$;
    \item on the level of $1$-morphisms, $( A\xrightarrow{\huaM_A} A', \quad B\xrightarrow{\huaM_B} B')$ is sent to $\huaM_A\oplus \huaM_B$, which is a $(A'\oplus B')$-$(A\oplus B)$-bimodule. 
    \item The direct sum of a $A'-A$-intertwiner $ \phi_A: \huaM_A\rightarrow \huaM'_A$ and a $B'-B$-intertwiner $\phi_B:\huaM_B\rightarrow \huaM'_B$  is a $(A'\oplus B')$-$(A\oplus B)$-intertwiner
    \[ \phi_A\oplus \phi_B: \huaM_A\oplus \huaM_B\longrightarrow \huaM'_A\oplus \huaM'_B \]
 \end{itemize}
The functor $\oplus$ is strict: on objects it is the usual direct sum of super algebras, on 1-morphisms it is the direct sum of bimodules and on 2-morphisms it is the direct sum of intertwiners. Strict associativity and strict unitality follow componentwise from the corresponding properties of direct sums of vector spaces.  The unit on \tVect \mbox{ }is the zero vector space.  The symmetry constraint $\beta$ is induced by the usual flip isomorphism $A \oplus B \cong B \oplus A$ in each component, which satisfies the hexagon identities. Hence $(
\text{\tVect}, \oplus, 0, \beta)$ is a strict symmetric monoidal bicategory in the sense of \cite[Definition 1.3]{Osorno2012SpectraAGT}.
\end{proof}

%since the direct sum of vector spaces is strictly associative and strictly unital and we have the equation \eqref{bimod_distr}. 
%This bifunctor is defined by the direct sum of vector spaces, thus, is strict. \comment{checked, correct.}

The spectrum construction of Osorno requires a symmetric monoidal bicategory. To apply it to \tVect, we need to replace \tVect \hspace{0.1mm} by its maximal 2-subgroupoid $M(\text{\tVect})$ (i.e., keeping only invertible 1- and 2-morphisms). This is necessary because the Duskin nerve of a bicategory only captures the invertible higher morphisms. A naive choice of Morita representatives would introduce an artificial dependence on the choice of representatives; however, different choices yield biequivalent bicategories, and hence weakly equivalent spectra. To avoid this artificial choice, we define $M(\text{\tVect})$ directly on all objects of \tVect.

\begin{definition} \label{def:m2vect}
Let $M(\text{\tVect})$ denote the bicategory defined as below.

\begin{itemize}
    \item Its objects are the objects in \tVect.

    \item The $1$-morphisms in $M(\text{\tVect})$ are all the invertible $1$-morphisms in \tVect.

    \item The $2$-morphisms are taken to be  all the $2$-isomorphisms in \tVect. %the invertible intertwiners between endomorphisms. 
\end{itemize} The structure 2-isomorphisms $a$, $l$ and $r$ are inherited from those in \tVect.
\end{definition}

%\begin{remark}Note, if $\huaM$ is an invertible $A-B$-bimodule, $\huaM'$ is an invertible $A'-B'$-bimodule, then $\huaM\oplus \huaM'$ is an invertible $(A\oplus A')-(B\oplus B')$-bimodule.  And the direct sum of invertible intertwiners is still invertible intertwiner.   The bicategory $M(\text{\tVect})$  inherits the strict symmetric monoidal structure from \tVect. \end{remark}

%\begin{remark}
%    Note that the 2-nerve of $M(\text{\tVect})$ is isomorphic to \[ \coprod_{A} B(\End(A))\] where $A$ goes over all the objects in \tVect and $B(\End(A))$ is the delooping of the category $\End(A)$.  \end{remark}

\begin{example} \label{bdr_2vect_embed}
Let $\textsf{s}\huaV ect_k$ denote the  symmetric bimonoidal category of super vector spaces and super linear maps  over the field $k$.
%Let $\textsf{s}\huaV ect_k$ denote the symmetric bimonoidal category of super vector bundles over the field $k$. 
Let $GL_n(\textsf{s}\huaV ect_k)$   denote the category of weakly invertible square matrices defined in \cite[Definition 3.6]{bdr}.
Let $BGL(\textsf{s}\huaV ect_k)$ denote the bicategory whose objects are $n\in N$ and the category of $1$-morphisms from $n$ to itself is  $GL_n(\textsf{s}\huaV ect_k)$.
    There is a normalized  pseudofunctor \[ I: BGL(\textsf{s}\huaV ect_k)\longrightarrow M(\text{\tVect})\] 
    \begin{itemize}
        \item sending each $n$ to $k^{\oplus n}$;
        \item sending each weakly invertible matrix $\bigl( V_{ij} \bigr)_{n\times n}$ to the $k^{\oplus n}$-$k^{\oplus n}$-bimodule $\bigl( V_{ij} \bigr)_{n\times n}$;
        \item sending a  linear isomorphism  $\bigl( \phi_{ij} \bigr)_{n\times n}$ from  $\bigl( V_{ij} \bigr)_{n\times n}$ to itself to the corresponding invertible intertwiner.
    \end{itemize} $I$ is strong monoidal. For any weakly invertible $n\times n$-matrices $V$, $V'$,  $I$ defines a bijection between the set of linear isomorphisms of an weakly invertible matrix $V$ and the set of the invertible intertwiners from $I(V)$ to itself.  In addition, 
    as we discussed invertible $k^{\oplus n}$-$k^{\oplus n}$-bimodules in detail further in Section \ref{2k_wr_general}, each  invertible $k^{\oplus n}$-$k^{\oplus n}$-bimodules can be interpreted as a weakly invertible matrix.  Thus, $I$ is  fully faithful.

\end{example}

To apply  \cite[Theorem 2.1]{Osorno2012SpectraAGT}, we need $M(\text{\tVect})$ to be a strict symmetric monoidal bicategory, which follows from Lemma \ref{2vect_strict}; restricting to invertible $1$- and $2$-morphisms preserves the strictness of the monoidal structure. Hence $M(\text{\tVect})$ satisfies the hypotheses of \cite[Theorem 2.1]{Osorno2012SpectraAGT}.

By \cite[Theorem 2.1]{Osorno2012SpectraAGT}, the classifying space of $M(\text{\tVect})$ is an infinite loop space upon group completion.
Then we apply the construction of classifying space of bicategory  \cite{CarrascoCegarraGarzon2010} and define the zeroth space $\huaK(\text{\tVect}) $ of the K-theory spectrum   using Duskin nerve (i.e. geometric nerve) $N^D_{\bullet}$ and the 2-vector space model in \cite{KLW2Vect22} instead of that in \cite{KV94_2} (see Definition \ref{def:infloop}). 
\begin{definition} \label{def:infloop}
    \[\huaK(\text{\tVect}) : = \Omega B (|N^D_{\bullet} (M(\text{\tVect}))|)   \]
\end{definition}

\begin{remark} 
The Duskin nerve is one of several nerve constructions for bicategories. By Carrasco–Cegarra–Garzón \cite{CarrascoCegarraGarzon2010}, all such nerves have homotopy equivalent geometric realizations. Hence the choice of nerve does not affect the classifying space or the resulting spectrum. We use the Duskin nerve throughout because of its full faithfulness and 3-coskeletal properties, which underpin the classification theorem in Section \ref{classify_thm}.
\end{remark}

\begin{remark} 

Applying \cite[Theorem 2.1]{Osorno2012SpectraAGT} to $M(\text{\tVect})$   yields a spectrum. We define the K-theory spectrum of  \tVect  \hspace{0.1mm} to be this spectrum:
 \[ \mathbb{K}(\text{\tVect}) : = \mathbb{A}(M(\text{\tVect})).\] The zeroth space of this spectrum is the infinite loop space $\huaK(\text{\tVect})$ in Definition \ref{def:infloop}. We apply the same symbol $\mathbb{A}$ as in \cite{Osorno2012SpectraAGT} to denote the functor from the category of strict symmetric monoidal bicategories to the category of spectra.
\end{remark}

\subsection{A nerve-theoretic classification}\label{classify_thm}

We now give a nerve-theoretic classification of 2-vector bundles over a Lie groupoid. The main result (Theorem \ref{thm:classify:2vectbdl}) identifies isomorphism classes of such bundles with homotopy classes of simplicial maps between the Duskin nerves of the Lie groupoid and the maximal 2-subgroupoid of the target 2-vector space.

Fix a bicategory $2\Vect$ of 2-vector spaces, for instance the bicategory of finite-dimensional super algebras, bimodules, and intertwiners.

In our (general) definition of 2-vector bundles over a Lie groupoid $\huaX_{\bullet}$, i.e. the bicategory of pseudofunctors \[ \pfun( \huaX_{\bullet} , 2\Vect),  \]
 the domain $\huaX_{\bullet}$ is a Lie groupoid, hence has only  $1$-isomorphisms and no non-identity 2-morphisms, only the invertible $1$-morphisms and  the 2-isomorphisms in $2\Vect$ play a role in the bicategory  $\pfun( \huaX_{\bullet} , 2\Vect)$  of pseudofunctors. Consequently, replacing $2\Vect$ by its maximal 2-subgroupoid $M(2\Vect)$ (i.e., keeping only invertible $1$- and  $2$-morphisms, as that in Definition \ref{def:m2vect}) does not change the bicategory \[ \pfun( \huaX_{\bullet} , 2\Vect) = \pfun( \huaX_{\bullet} , M(2\Vect)).  \] 
 In addition, it is standard that the bicategory of pseudofunctors is biequivalent to the bicategory of normalized pseudofunctors \cite[Remark 4.2]{benabou:bicat} \cite{Campbell2020}.

Below we provide a nerve-theoretic classification result in term of Duskin nerve \cite{Street1987}. 
%This observation is independent of the spectrum construction in Section \ref{2K_spectra} and does not extend to the full bicategory. It is included for completeness, as it provides a  categorical description of the invertible part of the 2K-theory.
We recall the following fundamental result due to Bullejos, Faro and Blanco \cite{BFB_2nerve}.

\begin{theorem}[Full faithfulness of the Duskin nerve] % \cite[Proposition 3]{BFB_2nerve}] 
Let $\textbf{Bi}\cat^{\mathrm{sup}}$ be the category of small bicategories and normal lax functors. There exists a nerve functor
\[ 
N^D_{\bullet}: \textbf{Bi}\cat^{\mathrm{sup}} \longrightarrow \sset \]
which is fully faithful. Explicitly, for any bicategory $\huaC$ and $\huaA$, there is a bijection between 
the normal lax functors $\huaC\longrightarrow  \huaA$ and the simplicial maps $N^D_{\bullet}(\huaC)   \longrightarrow  N^D_{\bullet}(\huaA)$.
\label{bfbprop3}
\end{theorem}

\begin{theorem} %[\cite[Proposition 5]{BFB_2nerve}]
If $\huaA$ is a 2-groupoid, then two pseudofunctors $F, G: \huaC \to \huaA$ have homotopic Duskin nerves if and only if there exists a lax  transformation between them. 

\label{bfbprop5}
\end{theorem}
Note that when the target bicategory $\huaA$ is a 2-groupoid,  all the lax transformations $F\Rightarrow G$ are strong transformations.

By \cite{duskin2}, the Duskin nerve $N^D_{\bullet}$ completely encodes all higher categorical structure of a bicategory and is 3-coskeletal. Thus, modifications correspond to higher homotopies. Thus, 
 when the target bicategory $\huaA$ is a 2-groupoid, the nerve functor $N^D_{\bullet}$ induces a bijection
 between the isomorphism classes of strong transformations $F\Rightarrow G$  and the equivalence classes of simplicial homotopies between the simplicial maps $N^D_{\bullet}(F)$ and $N^D_{\bullet}(G)$.
 
% of Hom-categories:\[ \textbf{Bi}\cat^{\mathrm{sup}} (F, G) \cong \sset(N(F), N(G)).\] 

%For fixed pseudofunctors $F, G: \huaX_{\bullet}\longrightarrow M (2\Vect)$, the Duskin nerve $N$ induces a bijection between 
%\begin{itemize}
%    \item strong transformations $\alpha: F\Rightarrow G$ (up to modifications), and 
%    \item homotopies between the simplicial maps $N(F)$ and $N(G)$. \end{itemize}

Applying Theorem \ref{bfbprop3}, \ref{bfbprop5} and the idea in \cite{duskin2}, we obtain:

\begin{theorem}[Classification of  2-vector bundles over a Lie groupoid] Let $\huaX_\bullet$ be a Lie groupoid. The Duskin nerve functor $N^D_{\bullet}$ induces an equivalence between the homotopy categories.
\[ 
\ho\bigl(\pFun(\huaX_{\bullet },        M (2\Vect)) \bigr)
\cong
\ho\bigl( \sset( N^D_{\bullet}(\huaX_{\bullet }),  N^D_{\bullet}(M (2\Vect))) \bigr).  \] 
\label{thm:classify:2vectbdl}\end{theorem}

In particular,   two  pseudofunctors $\huaX_{\bullet}\longrightarrow  M (2\Vect))$ are equivalent (as objects in the homotopy category) if and only if their nerves are homotopy equivalent.

\begin{remark}
Under geometric realization, the simplicial maps above correspond to continuous maps between the corresponding topological spaces. Thus the isomorphism of homotopy categories can be interpreted as a classification of 2-vector bundles by homotopy classes of maps
\[ |N^D_{\bullet} (\huaX_\bullet)| \longrightarrow |N^D_{\bullet}(M (2\Vect))|.\]
This interpretation is not needed for the proof of the classification theorem, but it provides a topological perspective on the result.
The simplicial-level isomorphism already gives the complete classification.
\end{remark}

\subsection{2K-theory of Lie groupoids} \label{sect:def_2K}

We now define the $2K$-theory of a Lie groupoid by applying the Grothendieck completion to the homotopy category of 2-vector bundles. The result is a category whose objects are formal differences of 2-vector bundles and whose morphisms are formal differences of 1-morphisms; ordinary and twisted K-theories embed as endomorphisms and morphisms of the trivial object.

We recall the standard  Grothendieck completion $K_0$ of a symmetric monoidal category (see e.g.  \cite[Chapter II, \S 5]{Weibel2013}  or \cite{Quillen1973}),
and apply it to the homotopy category $\ho(\textsf{s}2\huaV \huaB dl_k(\huaX_{\bullet}))$. This yields the following construction of $2K$-theory.

\begin{definition}[Grothendieck completion of a symmetric monoidal bicategory] \label{def:gr_comp}
Let $M$ be a symmetric monoidal bicategory and let $\ho(M)$ denote its homotopy category.
The \emph{Grothendieck completion} $K_0(M)$ is defined as follows:
\begin{itemize}
    \item \textbf{Objects:} the set of objects of $K_0(M)$ is the quotient of $\obj(M) \times \obj(M)$ 
    by the relation $(m,n) \simeq (m',n')$ if and only if there exists an object $c \in M$ 
    such that $m \oplus n' \oplus c \cong m' \oplus n \oplus c$ in $\ho(M)$.
    
    \item \textbf{Morphisms:} given two objects $(m,n)$ and $(p,q)$ in $K_0(M)$, 
    a morphism from $(m,n)$ to $(p,q)$ is an equivalence class of pairs 
    $(f,g) \in \Mor(\ho(M)) \times \Mor(\ho(M))$ with $f: m \to p$, $g: n \to q$.
    Two such pairs $(f,g)$ and $(f',g')$ are equivalent if there exist objects $c, d \in M$ 
    and a morphism $\alpha: c \to d$ in $\ho(M)$ such that the diagram
    \[
\xymatrix{
    m \oplus n' \oplus c  \ar[d]_{ f \oplus g' \oplus \alpha} \ar[r]^{\cong} & m' \oplus n \oplus c \ar[d]^{f' \oplus g \oplus \alpha}\\
    p \oplus q' \oplus d  \ar[r]_{\cong} & p' \oplus q \oplus d
    }
    \]
    commutes  in $\ho(M)$, where the horizontal arrows are 
    the isomorphisms witnessing the equivalences $(m,n) \simeq (m',n')$ and 
    $(p,q) \simeq (p',q')$.
\end{itemize}
Intuitively, for two $1$-morphisms $m \xrightarrow{f} p$ and $n \xrightarrow{g} q$ in $\ho(M)$, 
the pair $(f,g)$ represents a formal difference ``$f - g$'' from ``$m - n$'' to ``$p - q$''.
\end{definition}

The homotopy category $\ho(\textsf{s}2\huaV \huaB dl_k(\huaX_{\bullet}))$ inherits the 
symmetric monoidal structure from $\textsf{s}2\huaV \huaB dl_k(\huaX_{\bullet})$ 
with respect to the direct sum $\oplus$.

\begin{definition}[$2K$-theory] \label{def:2k}
Let $\huaX_\bullet$ be a Lie groupoid. 
The \emph{$2K$-theory} $2K(\huaX_{\bullet})$ of $\huaX_{\bullet}$ is the 
Grothendieck completion of a skeleton of the homotopy category 
$\ho(\textsf{s}2\huaV \huaB dl_k(\huaX_{\bullet}))$.
\end{definition}

Here the Grothendieck completion is taken in the sense of Definition \ref{def:gr_comp}; the objects and morphisms of $2K(\huaX_{\bullet})$
are thus formal differences of $2$-vector bundles and of the equivalence classes of $1$-morphisms, respectively. In the rest of the paper we use the same symbol to denote a category and any skeleton of it when there is no confusion.

\begin{remark}
    Unlike ordinary K-theory, our $2K$-theory is a category; its object set is a group, while the ring structure is recovered at the level of endomorphisms.
\end{remark}

%As shown in Theorem \ref{Thm_2k_2vect}, the $2K$-theory $2K(\huaX_{\bullet})$  gives a classification of the internal equivalence classes of 2-vector bundles over $\huaX_{\bullet}$. A direct corollary is that the $2K$-theory of a smooth manifold classifies the internal equivalence classes of 2-vector bundles over it.

\begin{example}
Let $X$ denote a smooth manifold, viewed as a Lie groupoid $X\git X$ with only identity arrows. 
Let  $\textsf{s}\huaA lg\huaB dl_k(X)$ denote the category of super algebra bundles over $X$ and  even  algebra bundle homomorphisms. It is   a   symmetric monoidal category with respect to the direct sum $\oplus$.
Denote by  $K^{\text{alg}}(X)$  the Grothendieck group of the isomorphism classes of super algebra bundles over  $X$. 
The embedding  $\tilde{\huaF}_{X\git X}:   \textsf{s}\huaA lg\huaB dl_k(X) \longrightarrow \textsf{s}2\huaV \huaB dl_k(X\git X)$ from Example \ref{pseudo_IX_pre} induces a group homomorphism  \[K^{\text{alg}}(X) \longrightarrow \obj(2K(X\git X)),  \] i.e. each Grothendieck class of super algebra bundles gives an object of the $2K$-theory. This homomorphism is compatible with the symmetric monoidal structures on both sides.
\end{example}

\subsection{Relation with orbifold K-theory and twisted K-theories} \label{rel_orb_twisted_k}

In this section we show the relation between orbifold K-theory and $2K$-theory and that between twisted K-theory and $2K$-theory.

A relation between orbifold K-theory and the $2K$-theory of a groupoid is given in Example \ref{ex_rel_orb_2K}.

\begin{example} \label{ex_rel_orb_2K}
Let $\huaX_{\bullet}$ be a proper \'etale  Lie groupoid. Let $\textsf{s}\huaV ect_k$ denote the category of super vector spaces and super linear maps. Let $\Fun(\huaX_{\bullet}, \textsf{s}\huaV ect_k)$  denote the symmetric monoidal category of orbifold super vector bundles over $\huaX_{\bullet}$, with the direct sum $\oplus$   and the tensor product $\otimes$   defined pointwise.

Let  
$B \Fun(\huaX_{\bullet}, \textsf{s}\huaV ect_k)$ denote the bicategory with a single object $\ast$, whose morphisms are the objects  of 
$\Fun(\huaX_{\bullet}, \textsf{s}\huaV ect_k)$ (i.e. orbifold super vector bundles), and whose $2$-morphisms are natural transformations. The composition of $1$-morphisms is given by the tensor product $\otimes$.

We define a strict 2-functor \[ I_K : B\Fun(\huaX_{\bullet}, \textsf{s}\huaV ect_k)  \longrightarrow \pFun(\huaX_{\bullet}, \text{\tVect})\]  as follows.
\begin{itemize}
    \item The single object $\ast$ is sent to the trivial $2$-vector bundle $\underline{k} = (k, k, l_k, \id_k)$. 

    \item Each functor $\huaX_{\bullet}\xrightarrow{f}   \textsf{s}\huaV ect_k$ is sent to a strong transformation $I_K(f): \underline{k}\rightarrow \underline{k} $ with 
    \begin{itemize}
        \item for each $x\in \huaX_0$, $I_K(f)_x:= f(x)$, which is a $k-k$-bimodule;
        \item for each $x\xrightarrow{\gamma} y$ in $\huaX_1$, $I_K(f)_{\gamma}: k\otimes_k f(x) \rightarrow f(y)\otimes_k k$ is defined to be the composition \[ r_{f(y)^{-1}} \circ f(\gamma) \circ l_{f(x)}.\]
    \end{itemize}

    \item Each natural transformation $\alpha: f\Rightarrow f'$ is sent to a modification $I_K(\alpha): I_K(f) \Rightarrow I_K(f')$ defined by $I_K(\alpha)_x = \alpha_x$ for each $x\in \huaX_0$.
\end{itemize}

In addition, the compositor and unitor are both defined to be the identity. Note that $I_K$ sends the associator, left unitor and right unitor of the tensor products of vector spaces to those of $k-k$-bimodules respectively. 

Moreover, from any strong transformation $\delta$ from $\underline{k}$ to itself, we can define a functor $J_K(\delta): \huaX_{\bullet}\rightarrow \textsf{s}\huaV ect_k$
with $J_K(\delta)(x)= \delta_x$, for any $x\in \huaX_0$, and, for any $x\xrightarrow{\gamma} y$ in $\huaX_1$, 
$J_K(\delta)(\gamma) = r_{\delta_y}\circ \delta_{\gamma}\circ l_{\delta_x}^{-1}$. For any $(\gamma, \gamma')\in \huaX_2$, $J_K(\delta)(\gamma\circ \gamma') = J_K(\delta)(\gamma)\circ J_K(\delta)(\gamma')$ by the lax naturality of $\delta$; and $J_K(\delta)(\id_x) = \id_x $ by the lax unity of $\delta$. 

And, from each modification $\alpha: \delta\Rightarrow \delta$ in $\pFun(\huaX_{\bullet}, \text{\tVect})$, we can define a natural transformation 
$J_K(\alpha): J_K(\delta)\Rightarrow J_K(\delta)$ by $J_K(\delta)_x := \alpha_x$, for each $x\in \huaX_0$.  It is indeed a natural transformation by the modification axiom of $\alpha$. 

\bigskip
Thus, we obtain an isomorphism of categories
\begin{equation}
    \Fun(\huaX_{\bullet}, \textsf{s}\huaV ect_k) \cong \End_{\textsf{s}2\huaV \huaB dl_k (\huaX_{\bullet})} (\underline{k}).
\end{equation}
After taking the homotopy categories of both the delooping $B \Fun(\huaX_{\bullet}, \textsf{s}\huaV ect_k)$ and the bicategory $$\pFun(\huaX_{\bullet}, \text{\tVect}),$$ the isomorphism classes of vector bundles over $\huaX_{\bullet}$ are precisely $\pi_0( \End_{\textsf{s}2\huaV \huaB dl_k  (\huaX_{\bullet})} (\underline{k})) $.

Since both sides are symmetric monoidal under direct sum, the isomorphism of categories induces an isomorphism of Grothendieck groups: $$
K_0(\ho(I_K)): K_{orb}(\huaX_{\bullet}) \longrightarrow  K_0 ( \pi_0(\End_{\textsf{s}2\huaV \huaB dl_k  (\huaX_{\bullet})} (\underline{k}))).$$

To see that this map preserves the direct sum, note that for any orbifold super vector bundles $V$ and $W$, their direct sum $V \oplus W$ has a natural $k-k$-bimodule structure (with $k$ acting diagonally on both factors). Under the embedding $I_K$, this corresponds to the direct sum of the corresponding endomorphisms of $\underline{k}$:
\[I_K(V \oplus W) \cong I_K(V) \oplus I_K(W), \]
where the direct sum on the right is the usual direct sum of $k-k$-bimodules. Hence the induced map on Grothendieck groups is compatible with the symmetric monoidal structures. In particular, the endomorphism category $\End_{\textsf{s}2\huaV \huaB dl_k (\huaX_\bullet)}(\underline{k})$ is closed under direct sums, so its Grothendieck completion is well-defined and embeds into the full $2K$-theory of $\huaX_\bullet$.

\end{example}

We now turn to a second instance of the same phenomenon: twisted K-theory also embeds into $2K$-theory, not as endomorphisms of the trivial object, but as morphisms from the trivial object to a twisting. 

By \cite[Definition 1.78]{Freed_vienna} and \cite[Section 3.2]{KLW2Vect22}, 
a twisting of complex K-theory over a Lie groupoid $\huaX_\bullet$ 
corresponds precisely to a super 2-vector bundle 
$\huaV = ( \huaA, \huaM, \mu, u)$ over 
$\huaX_\bullet$ (Definition \ref{def:2vect:obj}) whose 
underlying super algebra bundle $\huaA$ is Azumaya (i.e. an invertible algebra bundle over $\huaX_0$).

%In addition, we reinterpret twisted vector bundle in the setting of this paper. 

Recall from  \cite[Definition 1.86]{Freed_vienna} that, for a twisting $\tau$ for K-theory represented by an invertible 2-vector bundle $$\huaA_{\tau}= (\huaA, \huaB, \lambda, \id)$$ over $\huaX_{\bullet}$,  a $\tau$-twisted vector bundle over $\huaX_{\bullet}$ consists of an $\huaA $-module $E_0 \to \huaX_0$  and an isomorphism
\[ \psi: \huaB\otimes_{s^*\huaA} s^*E_0 \longrightarrow t^* E_0\] of $t^*\huaA$-modules over $\huaX_1$, which satisfy a cocycle condition on $\huaX_2$.

\begin{example} \label{ex:twisted:2K}
    Let $\huaE:= (E_0, \psi)$ be a $\tau$-twisted vector bundle over a Lie groupoid $\huaX_{\bullet}$. Note that each fiber $(E_0)_x$ is naturally a left $\huaA_x$-module and a right $k$-module, hence an $\huaA_x-k$-bimodule. It gives precisely a strong transformation $T(\huaE): \underline{k}\Rightarrow \huaA_{\tau}$ with $T(\huaE)_{x} = (E_0)_x$, for each $x\in \huaX_0$, and, for any $x\xrightarrow{\gamma} y$ in $\huaX_1$,  $$T(\huaE)_{\gamma}:= r^{-1}_{(E_0)_y} \circ \psi_{\gamma}: B_{\gamma}\otimes_{\huaA_x} (E_0)_x \rightarrow (E_0)_y
    \otimes_k k \cong (E_0)_y.$$ The cocycle condition of the twisted vector bundle is precisely the lax naturality for strong transformation, while the unity axiom holds automatically. Thus, a $\tau$-twisted vector bundle over $\huaX_{\bullet}$ is a strong transformation from $\underline{k}$ to $ \huaA_{\tau}$. 

    In the other direction, given any strong transformation $\alpha: \underline{k}\Rightarrow \huaA_{\tau}$, we can define a $\tau$-twisted vector bundle $Q(\alpha) = (E_{\alpha}, \psi_{\alpha})$ over $\huaX_{\bullet}$ immediately with $(E_{\alpha})_x:= \alpha_x$, for any $x\in \huaX_0$, and the isomorphism $$(\psi_{\alpha})_\gamma:=  r_{\alpha_y}  \circ   \alpha_\gamma$$ for any $x\xrightarrow{\gamma} y$ in $\huaX_1$. 

    Thus, $\tau$-twisted vector bundles over $\huaX_{\bullet}$   are in bijection with $1$-morphisms $$\underline{k}\Rightarrow \huaA_{\tau}$$ in the bicategory $\textsf{s}2\huaV \huaB dl_k  (\huaX_{\bullet})$ of super 2-vector bundles over $\huaX_{\bullet}$. 
    
    In addition, a $\huaA_{\tau}$-module map between two $\tau$-twisted vector bundles $$\huaE = (E_0, \psi)\longrightarrow  \huaE' = (E'_0, \psi')$$ consists of a bundle map $f: E_0 \rightarrow E'_0$ which makes the diagram below commutes
    \begin{equation}\label{twist:mor:mod}
      \xymatrix{  \huaB\otimes_{s^*\huaA} s^*E_0 \ar[d]_{\id\otimes f} \ar[r]^{\psi} &t^* E_0 \ar[d]^{f}
\\       \huaB\otimes_{s^*\huaA} s^*E_0' \ar[r]_{\psi'} \ar@{}[ru]|{\circlearrowright} & t^* E_0'}
    \end{equation} $f$ corresponds a modification $T(f): T(\huaE)\longrightarrow T(\huaE')$ defined by $T(f)_{x} = f_x$ for each $x\in \huaX_0$; the diagram \eqref{twist:mor:mod} corresponds  exactly to the  modification axiom.
    For the other direction, given a modification $T(\huaE)\longrightarrow T(\huaE')$, we can recover the $\huaA$-module map $\huaE\rightarrow \huaE'$ in the canonical way. Thus, there is a one-to-one correspondence between the modifications  $T(\huaE)\longrightarrow T(\huaE')$ and  the $\huaA_{\tau}$-module maps $\huaE\rightarrow \huaE'$.

    Thus, the isomorphism classes of $\tau$-twisted vector bundles are in bijection with $$\pi_0(\hom_{\textsf{s}2\huaV \huaB dl_k  (\huaX_{\bullet})} (\underline{k}, \huaA_{\tau})).$$ After taking Grothendieck group of both, we have \[ K^{\tau}(\huaX_{\bullet}) \cong K_0(\pi_0(\hom_{\textsf{s}2\huaV \huaB dl_k  (\huaX_{\bullet})} (\underline{k}, \huaA_{\tau}))).\] Here the Grothendieck group on the right is taken with respect to the direct sum of $1$-morphisms, which corresponds to the direct sum of twisted vector bundles.
\end{example}

\section{2-equivariant 2-Vector bundles and 2-equivariant $2K$-theory} \label{sect_2eq_2vb_2k}

\subsection{2-equivariant 2-Vector bundles over a Lie groupoid} \label{2eq_2vb}

%It may be a little weird to construct the 2-vector bundles over the category of Lie groupoids whereas we construct equivariant 2-vector bundles in the bicategory of Lie groupoids acted by a coherent Lie 2-group $\huaG$, which is a sub-bicategory of the bicategory of weak groupoid objects in $\bgpd$. But that is the case we need.

In this section we construct a bicategory of $\huaG_{\bullet}$-equivariant 2-vector bundles over a Lie groupoid where $\huaG_{\bullet}$ is a coherent Lie 2-group.  Our definition of equivariant 2-vector bundles separates the action on the base from the action on the fibres. This generality is essential for several applications discussed in this paper: it covers the classical equivariant case (where the base is an action groupoid), the representation-theoretic case (where the base is a point), and the gauge-theoretic case (where the 2-group acts only on the fibres).

We first  define the main object of interest: a $\huaG_{\bullet}$-equivariant super 2-vector bundle. 
\begin{definition} \label{def:2eq2vectgrpd}

Let $\huaG_\bullet$ be a coherent 2-group and $\huaX_{\bullet}$  a Lie groupoid. 
A $\huaG_\bullet$-equivariant super 2-vector bundle over $\huaX_\bullet$ is a pseudofunctor
\[
\huaT: B\huaG_{\bullet}  \longrightarrow \textsf{s}2\huaV \huaB dl_k(\huaX_\bullet),
\]
where $B\huaG_{\bullet}$ is the delooping bicategory of $\huaG_\bullet$, and $\textsf{s}2\huaV \huaB dl_k(\huaX_\bullet)$ is the bicategory of super 2-vector bundles over $\huaX_\bullet$.

Explicitly,  $\huaT$ consist of:
\begin{itemize}
    \item  a super 2-vector bundle $\huaV: = \huaT(\ast) $ over $\huaX_\bullet$;
    \item for each $g_0 \in \huaG_0$,  a 1-morphism
    \[
    \huaT(g_0) : \huaV \longrightarrow \huaV
    \]
    in $\textsf{s}2\huaV \huaB dl_k(\huaX_\bullet)$; %, such that for each $x \in \huaX_0$,\[    (\huaT(g_0))_{x}: \huaV(x) \longrightarrow \huaV(g_0\cdot x)    \]    is a $ \huaV(g_0\cdot x)-\huaV(x)$-bimodule;
    \item for each $1$-morphism $ (h_1: g_0 \to g'_0 ) \in \huaG_1$, a 2-morphism
    \[
    \huaT(h_1): \huaT(g_0) \Longrightarrow \huaT(g'_0)
    \]
    in $\textsf{s}2\huaV \huaB dl_k(\huaX_\bullet)$;
    \item for each   pair $(g_2, g_1)$ of objects in $\huaG_{\bullet}$, a 2-isomorphism  in $\textsf{s}2\huaV \huaB dl_k(\huaX_\bullet)$
    \[
    \huaT^2_{g_2, g_1}: \huaT(g_2) \circ \huaT(g_1) \Longrightarrow \huaT(g_2 \circ g_1);
    \]
    \item  a 2-isomorphism  in $\textsf{s}2\huaV \huaB dl_k(\huaX_\bullet)$
    \[
    \huaT^0_{\ast} : \id_{\huaV} \Longrightarrow \huaT(e);
    \] where $e$ is the unit object in $\huaG_{\bullet}$.
   
\end{itemize}
These data satisfy the pseudofunctor axioms: the lax associativity pentagon and the lax left/right unity triangles.
\end{definition}

\begin{remark} \label{rmk:motiv:2eqectbdl}
By the standard correspondence between 2-group actions and pseudofunctors from the delooping (see \cite[Section 8.3, 10.7]{JY:Bicat}), a pseudofunctor
\[ \huaT: B\huaG_\bullet \longrightarrow \textsf{s}2\huaV \huaB dl_k(\huaX_\bullet) \] 
is equivalently a $\huaG_\bullet$-action on the 2-vector bundle $ \huaV:= \huaT(\ast)$. Thus Definition \ref{def:2eq2vectgrpd} indeed captures the notion of a $\huaG_\bullet$-equivariant 2-vector bundle.

Note that this action is fiberwise: $\huaG_\bullet$ acts on the 2-vector bundle $\huaV$ itself, not on the base Lie groupoid $\huaX_\bullet$. This is the categorical analogue of a group action on the fibres of a vector bundle, without requiring an action on the base.

Classical equivariance over a base with a group action is recovered as a special case by taking $\huaX_\bullet $ to be the action groupoid of a smooth $G$-action on a manifold $X$ and requiring that the pseudofunctor $\huaT: B G  \longrightarrow \textsf{s}2\huaV \huaB dl_k( \huaX_{\bullet})$ is compatible with the groupoid structure of the action groupoid; see Remark \ref{ex_equiv_2vectbdl} for an explicit calculation.
\end{remark}

\begin{remark}[The classical equivariant case] 

\label{ex_equiv_2vectbdl}

To see how Definition \ref{def:2eq2vectgrpd} recovers the familiar notion of equivariant bundles, suppose   $\huaG_{\bullet} =G$  is an ordinary  Lie group $G$ (viewed as a trivial $2$-group) and let   $\huaX_{\bullet} = (G\ltimes X    \rightrightarrows X) $  be the action groupoid of a smooth left $G$-action on a manifold $X$. The recovery is not automatic: it requires that the pseudofunctor data $\{\huaT(g)\}_{g\in G}$ are compatible with the action groupoid structure, as we now spell out.

In this case, a pseudofunctor \[ \huaT: B G  \longrightarrow \textsf{s}2\huaV \huaB dl_k( \huaX_{\bullet})\] consists of

\begin{itemize}
   \item  a fixed super 2-vector bundle $\huaV: = \huaT(\ast) $ over $ \huaX_{\bullet}$;
    \item for each $g \in G $,  a 1-morphism $    \huaT(g) : \huaV \longrightarrow \huaV    $;
    \item for each   pair $(g_2, g_1)$,  the compositor  $\huaT^2_{g_2, g_1}: \huaT(g_2) \circ \huaT(g_1) \Longrightarrow \huaT(g_2 \circ g_1);$
    \item  the unitor  $
    \huaT^0_{\ast} : \id_{\huaV} \Longrightarrow \huaT(e);$
\end{itemize}
Since the base is the action groupoid, $\huaV$ itself contains the equivariant data: for each arrow $(g,x): x \to g\cdot x$ in $G \ltimes X$, the invertible bimodule $\huaT(g)_x: \huaV(x) \to \huaV(g\cdot x)$ encodes the ``pullback" isomorphism (see \cite[Section 3.8]{KLW2Vect22}).

The lax associativity of $\huaT$ reduces to the cocycle condition. 
Over each $(g_1, g_2, x)\in G\ltimes G\ltimes X$, $\huaT^2_{g_1, g_2}$ is an invertible modification in $\textsf{s}2\huaV \huaB dl_k( \huaX_{\bullet})$, and $(\huaT^2_{g_1, g_2})_x$ is a 2-isomorphism from 
$\huaT(g_1)_{(g_2\cdot x)} \circ \huaT(g_2)_x$ to $\huaT(g_1g_2)_x$.  In addition, the lax associativity is exactly
\begin{equation}
    (\huaT^2_{g_1g_2, g_3})_x \vertprod ( (\huaT^2_{g_1, g_2})_{(g_3\cdot x)}  \horiprod 1_{\huaT(g_3)_x}) = (\huaT^2_{g_1, g_2g_3})_x \vertprod (1_{\huaT(g_1)_{((g_2g_3)\cdot x)}} \horiprod (\huaT^2_{g_2, g_3})_x). 
\end{equation}
And, the lax left and right unity over each  $(g, x)\in G\ltimes X$ are given by   \begin{align*}
(\huaT^2_{e, g})_x \vertprod  ( (\huaT^0_{\ast})_{(g\cdot x)} \horiprod (1_{\huaT(g)})_x   )   & =  l_{\huaT(g)_x}\\
(\huaT^2_{g, e})_x \vertprod  (  (1_{\huaT(g)})_x  \horiprod (\huaT^0_{\ast})_{ x }    )   & =  r_{\huaT(g)_x}
\end{align*} Geometrically, these two identities ensure that the parallel transport encoded by the compositor is compatible with the identity arrows of the action groupoid. 

This example illustrates the mechanism by which Definition \ref{def:2eq2vectgrpd} subsumes the equivariant 2-vector bundle in the sense of \cite{KLW2Vect22}, where the equivariance is encoded by bundle isomorphisms over the action groupoid.
\end{remark}

\bigskip

Next we define $1$-morphisms between $\huaG_{\bullet}$-equivariant 2-vector bundles over $\huaX_{\bullet}$.    
\begin{definition} \label{def:2eq2vb_1mor}

Let $\huaT, \huaS: B\huaG_\bullet \to \textsf{s}2\huaV \huaB   dl(\huaX_\bullet)$ be two $\huaG_{\bullet}$-equivariant super 2-vector bundles over $\huaX_{\bullet}$. 
A 1-morphism $\alpha: \huaT \Rightarrow \huaS$ is a strong transformation between them, consisting of:
\begin{itemize}
\item for the only object $\ast$ in $B\huaG_\bullet$, a $1$-morphism
\[ \alpha_{\ast}: \huaT(\ast)\rightarrow \huaS(\ast )\]   in $\textsf{s}2\huaV \huaB   dl(\huaX_\bullet)$;
    \item for each $g_0 \in \huaG_0$, a 2-isomorphism
    \[
    \alpha_{g_0}: \huaS(g_0) \circ \alpha_{\ast} \Longrightarrow \alpha_{\ast } \circ \huaT(g_0)
    \]
    in $\textsf{s} 2\huaV \huaB   dl(\huaX_\bullet)$,
\end{itemize}
such that the lax associativity and lax unity hold.
\end{definition}

\begin{definition} \label{def:2eq2vb_2mor}

Let $\alpha, \beta: \huaT \Rightarrow \huaS$ be two 1-morphisms of $\huaG_{\bullet}$-equivariant super 2-vector bundles over $\huaX_{\bullet}$. 
A 2-morphism $\gamma: \alpha \Rightarrow \beta$ is a modification, consisting of:
\begin{itemize}
    \item for   the only object $\ast$ in $B\huaG_{\bullet}$, a 1-morphism
    \[
    \gamma_{\ast}: \alpha_{\ast} \Longrightarrow \beta_{\ast}
    \]
    in $\textsf{s}2\huaV \huaB dl_k(\huaX_\bullet)$,
\end{itemize}
such that for each $g_0 \in \huaG_0$, the corresponding modification axiom is satisfied.

\begin{equation} \label{def:2mor_2eq2vectbd}
       \xymatrix{ &\huaT(\ast)         \ar[r]^-{\huaT(g_0)} \ar[dd]|{\alpha_\ast} &   \huaT(\ast)  \ar[dd]|{\alpha_{\ast}}      \ar@/^4pc/[dd]^{\beta_{\ast}}   \ar@{}@<+6ex>[dd]|{\buildrel{\gamma_{\ast}}\over\Rightarrow} \\ & & \\
   & \huaS(\ast)  
   \ar@/_1pc/[r]_-{ \huaS(g_0)} 
      \ar@{}@<+0.5ex>[ruu]|{\buildrel{\alpha_{g_0}}\over\Rightarrow}      &  \huaS(\ast) 
    \\ } =
    \xymatrix{ \huaT(\ast)     \ar@/_5pc/[dd]_{\alpha_{\ast}}       \ar[r]^-{\huaT(g_0)} \ar[dd]|{\beta_{\ast}} &  \huaT(\ast)    \ar[dd]|{\beta_{\ast}}  & \\ & & \\
    \huaS(\ast)          \ar@{}@<+7ex>[uu]|{\buildrel{\gamma_{\ast}}\over\Rightarrow}     \ar@/_1pc/[r]_-{\huaS(g_0)} 
      \ar@{}@<+0.5ex>[ruu]|{\buildrel{\beta_{g_0}}\over\Rightarrow}   & \huaS(\ast)   &\\ } 
\end{equation}

\end{definition}

%The bicategory of $\huaG_{\bullet}$-equivariant 2-vector bundles over $\huaX_{\bullet}$ is precisely the Grothendieck construction $ \int F$ of the pseudofunctor $F: B\huaG_{\bullet} \longrightarrow 2\huaV \huaB dl_k(\huaX_{\bullet})  $ that classifies the $\huaG_{\bullet}$-action.

The $\huaG_{\bullet}$-equivariant 2-vector bundles over $\huaX_{\bullet}$ (Definition \ref{2eq_2vb}), the $1$-morphisms (Definition \ref{def:2eq2vb_1mor}) and the $2$-morphisms (Definition \ref{def:2eq2vb_2mor}) together define the bicategory $$(\textsf{s}2\huaV \huaB  dl_k)_{\huaG_{\bullet}}(\huaX_{\bullet}) $$ of $\huaG_{\bullet}$-equivariant 2-vector bundles over $\huaX_{\bullet}$. 
The vertical and horizontal compositions are inherited from those of strong transformations and modifications in $\textsf{s}2\huaV \huaB  dl_k(\huaX_\bullet)$; the associativity and unit constraints are defined from those of the underlying bicategory.

\begin{theorem}

\leavevmode
    \begin{enumerate}
        \item If the coherent Lie 2-groups $\huaG_{\bullet}$ and $\huaG'_{\bullet}$ are biequivalent, then the bicategory $(\textsf{s}2\huaV \huaB  dl_k)_{\huaG_{\bullet}}(\huaX_{\bullet}) $ is biequivalent to $(\textsf{s}2\huaV \huaB  dl_k)_{\huaG'_{\bullet}}(\huaX_{\bullet}) $ for each Lie groupoid $\huaX_{\bullet}$. 

        \item  If $\huaX_{\bullet}$ and $\huaY_{\bullet}$ are biequivalent in the bicategory $\bgpd$ of Lie groupoids, bibundles and equivariant bibundle maps, then the bicategories  $(\textsf{s}2\huaV \huaB  dl_k)_{\huaG_{\bullet}}(\huaX_{\bullet}) $ and  $(\textsf{s}2\huaV \huaB  dl_k)_{\huaG_{\bullet}}(\huaY_{\bullet}) $ are biequivalent.
            \end{enumerate}
\end{theorem}

\begin{proof}
\leavevmode

   \begin{enumerate}
       \item %The proof is routine and analogous step by step to that of \cite[Theorem 3.1]{Huan:2Rep_2Vect}.

 Suppose $B\huaG_\bullet$ and $B\huaG'_\bullet$ are biequivalent, with pseudofunctors
\[ \Phi: B\huaG_\bullet \to B\huaG_\bullet', \qquad \Psi: B\huaG_\bullet' \to B\huaG_\bullet \]
and invertible strong transformations
\[ \eta: \id_{B\huaG_\bullet} \Rightarrow \Psi \circ \Phi, \qquad \epsilon: \Phi \circ \Psi \Rightarrow \id_{B\huaG'_\bullet}\]
satisfying the triangle identities.

Consider the precomposition pseudofunctors \begin{align*}
& \Phi^*: \pfun(B\huaG'_\bullet, \textsf{s}2\huaV \huaB  dl_k(\huaX_\bullet)) \longrightarrow \pfun(B\huaG_\bullet, \textsf{s}2\huaV \huaB  dl_k(\huaX_\bullet)), \\
& \Psi^*: \pfun(B\huaG_\bullet, \textsf{s}2\huaV \huaB  dl_k(\huaX_\bullet)) \longrightarrow \pfun(B\huaG'_\bullet, \textsf{s}2\huaV \huaB  dl_k(\huaX_\bullet)), \end{align*}
given by $F \mapsto F \circ \Phi$ and $G \mapsto G \circ \Psi$, respectively. By the 2-functoriality of $\pfun(-, \huaC)$ (see the proof of Theorem \ref{biequiv_2vectbdl}), $\eta$ and $\epsilon$ induce invertible strong transformations
\[ \tilde{\eta}: \id \Rightarrow (\Psi \circ \Phi)^*, \qquad \tilde{\epsilon}: (\Phi \circ \Psi)^* \Rightarrow \id, \]
which inherit the triangle identities from $\eta$ and $\epsilon$. Hence $(\Phi^*, \Psi^*, \tilde{\eta}, \tilde{\epsilon})$ forms a biequivalence. Since
\[ (\textsf{s}2\huaV \huaB  dl_k)_{\huaG_\bullet}(\huaX_\bullet) = \pfun(B\huaG_\bullet, \textsf{s}2\huaV \huaB  dl_k(\huaX_\bullet)),\]
the claim follows.

\item 
   %Via Theorem \ref{biequiv_2vectbdl}.   
   
   Suppose $\huaX_\bullet$ and $\huaY_\bullet$ are biequivalent in $\bgpd$. By Theorem \ref{biequiv_2vectbdl}  (2), we have a biequivalence
\[ \textsf{s}2\huaV \huaB  dl_k(\huaX_\bullet) \simeq \textsf{s}2\huaV \huaB  dl_k(\huaY_\bullet). \]
Let \[ 
G^*: \textsf{s}2\huaV \huaB  dl_k(\huaX_\bullet) \to \textsf{s}2\huaV \huaB  dl_k(\huaY_\bullet), \qquad
F^*: \textsf{s}2\huaV \huaB  dl_k(\huaY_\bullet) \to \textsf{s}2\huaV \huaB  dl_k(\huaX_\bullet) \]
be the pseudofunctors inducing this biequivalence, with unit and counit strong transformations inherited from the biequivalence in $\bgpd$.

Define pseudofunctors
\begin{align*}
&(G^*)_*: (\textsf{s}2\huaV \huaB  dl_k)_{\huaG_\bullet}(\huaX_\bullet) \longrightarrow (\textsf{s}2\huaV \huaB  dl_k)_{\huaG_\bullet}(\huaY_\bullet), \\
&(F^*)_*: (\textsf{s}2\huaV \huaB  dl_k)_{\huaG_\bullet}(\huaY_\bullet) \longrightarrow (\textsf{s}2\huaV \huaB  dl_k)_{\huaG_\bullet}(\huaX_\bullet)     
\end{align*}
by post-composition with $G^*$ and $F^*$, respectively:
\[ (G^*)_*(\huaT) := G^* \circ \huaT, \qquad (F^*)_*(\huaT') := F^* \circ \huaT'. \]
The unit and counit of the biequivalence $G^* \dashv F^* $ induce, by whiskering with objects of the representation bicategories, the required unit and counit strong transformations between $(G^*)_*$ and $(F^*)_*$. These satisfy the triangle identities because they are inherited from the biequivalence between $\textsf{s}2\huaV \huaB  dl_k(\huaX_\bullet)$  and $\textsf{s}2\huaV \huaB  dl_k(\huaY_\bullet)$. Therefore $(G^*)_* $ and $(F^*)_*$ form a biequivalence.  
     \end{enumerate} 
\end{proof}

\subsection{2-equivariant 2K-theory and examples} \label{subsect:def:2eq2k}

We now extend the $2K$-theory construction of Section \ref{sect:def_2K} to the equivariant setting. Given a coherent Lie 2-group $\huaG_\bullet$ and a Lie groupoid $\huaX_\bullet$, we replace the bicategory $\textsf{s}2\huaV \huaB dl_k(\huaX_\bullet)$ by its $\huaG_\bullet$-equivariant version $(\textsf{s}2\huaV \huaB dl_k)_{\huaG_\bullet}(\huaX_\bullet)$ from Definition \ref{def:2eq2vectgrpd}. Applying the same Grothendieck completion procedure to its homotopy category yields the 2-equivariant $2K$-theory.

\begin{definition}[2-Equivariant $2K$-theory] \label{def:2equiv_2k} 
    The $\huaG_{\bullet}$-equivariant  $2K$-theory $2K_{\huaG_{\bullet}}(\huaX_{\bullet})$ of $\huaX_{\bullet}$ is defined to be  the 
Grothendieck completion of the homotopy category 
$\ho( (\textsf{s}2\huaV \huaB dl_k)_{\huaG_{\bullet}}  (\huaX_{\bullet}) )$.
\end{definition}

\begin{remark}
    The Grothendieck completion defined in Definition \ref{def:2equiv_2k}  gives the 0-dimensional part of a putative 2-equivariant $2K$-theory. A genuine equivariant spectrum lifting this definition would require an equivariant version of Osorno's $\mathbb{A}$-functor, which is beyond the scope of this paper. In Appendix \ref{subsect:2k_spectrum} we construct a spectrum associated to the $\huaG_{\bullet}$-representation bicategory $2\Rep(\huaG_\bullet)$, not to the $\huaG_{\bullet}$-equivariant $2$-vector bundles over a general Lie groupoid $\huaX_{\bullet}$.
\end{remark}

\begin{remark}
    The present $2K$-theory is the additive $0$-dimensional part of a higher categorical K-theory; the full multiplicative structure, which would make it a rig category in the sense of \cite{bdr} and a commutative ring in the stable limit, is deferred to future work.
\end{remark}

Below we give some  examples.

%\begin{example} Let $$K_G(X)$$ denote the Grothendieck group of the isomorphism classes of $G$-equivariant super algebra bundles over a smooth manifold $X$ acted by a compact Lie group $G$.  There is a group homomorphism $K_G(X) \longrightarrow 2K_{G\git G}(X\git X)$ which sends the isomorphism classes of $G$-equivariant super algebra bundles over $X$ to the internal equivalence  classes of $G\git G$-equivariant 2-vector bundles over $X\git X$. \end{example}

\begin{example}\label{2rep_concrete}  
    Let $\ast$ be the terminal groupoid and $\huaG_\bullet$ a coherent 
Lie $2$-group. 
Recall the bicategory of  $\huaG_\bullet$-$2$-representation is defined to be the bicategory $\pFun( B\huaG_{\bullet },  \text{\tVect})$ of pseudofunctors $B\huaG_\bullet \to \text{\tVect}$, strong transformations and modifications. %, which is of the same principle as \cite[Definition 2.4]{Huan:2Rep_2Vect}. 
 A $\huaG_\bullet$-equivariant $2$-vector bundle over 
$\ast$ is precisely a $\huaG_\bullet$-$2$-representation.

The homotopy category of $\huaG_\bullet$-equivariant $2$-vector bundles 
over $\ast$ is equivalent to $\ho(2\Rep(\huaG_\bullet))$. 
Applying the Grothendieck completion yields the identification 
$2K_{\huaG_\bullet}(\ast) \cong K_0(\ho(2\Rep(\huaG_\bullet)))$.
\end{example}

\begin{example}[Classical equivariant K-theory via the action groupoid] \label{ex:act:eq}

This is precisely the equivariant version of Example \ref{ex_rel_orb_2K}: we consider the specific case that  the base Lie groupoid $\huaX_\bullet$ is given by the action groupoid $G \ltimes X \rightrightarrows X$.

Let $G$ be a compact Lie group acting smoothly on a manifold $X$. We just use  $$G \ltimes X$$ to denote the associated action groupoid. Consider the non-equivariant $2K$-theory $2K(G \ltimes X)$. Then, by the same argument as in Example \ref{ex_rel_orb_2K}, we have an isomorphism of categories
\begin{equation} (\textsf{s}\huaV \huaB dl_k)_G(X) \cong \End_{\textsf{s}2\huaV \huaB dl_k(G \ltimes X)}(\underline{k}),
\end{equation}
where $(\textsf{s}\huaV \huaB dl_k)_G(X)$ denotes the category of $G$-equivariant super vector bundles over $X$. Passing to Grothendieck groups gives
\begin{equation}  
K_{G}(X) \cong  K_0 ( \pi_0(\End_{\textsf{s}2\huaV \huaB dl_k  (G \ltimes X)} (\underline{k}))) . \end{equation}

Thus classical equivariant K-theory is recovered as the endomorphism ring of the trivial object in the non-equivariant 2K-theory of the action groupoid. When $X = \pt$, this gives $K_G(\pt) \cong R(G)$, the representation ring of $G$.
\end{example}

\begin{example}[Classical equivariant twisted K-theory via the action groupoid] \label{ex:act:eq:twisted}

Similarly, let $\tau$ be a $G$-equivariant twisting over $X$, represented by an invertible 2-vector bundle $\huaA_\tau$ over the action groupoid $G \ltimes X$. By the same argument as in Example \ref{ex:twisted:2K}, applied to the action groupoid $G \ltimes X$, we have a bijection between
$\{G\text{-equivariant }\tau\text{-twisted vector bundles over }X\}$ and  $\hom_{\textsf{s}2\huaV \huaB dl_k(G \ltimes X)}(\underline{k}, \huaA_\tau).$ 
Taking Grothendieck groups yields
\[ 
K_G^\tau(X) \cong K_0\Big( \pi_0\big( \hom_{\textsf{s}2\huaV \huaB dl_k(G \ltimes X)}(\underline{k}, \huaA_\tau) \big) \Big). \]
Thus classical equivariant twisted K-theory is recovered as the morphisms from the trivial object to a twisting in the non-equivariant 2K-theory of the action groupoid.
\end{example}

\begin{remark}
We emphasize that these recoveries, as in Example \ref{ex:act:eq} and \ref{ex:act:eq:twisted}, take place in the non-equivariant $2K$-theory of the action groupoid, which should not be confused with the $\huaG_\bullet$-equivariant $2K$-theory of Definition \ref{def:2equiv_2k}, which captures a higher notion of equivariance.

\end{remark}

\section{ Relation with the Stacky Definition of 2-vector bundles} \label{bieq_2model}

\subsection{A Stacky model for 2-vector bundles over Lie groupoids}\label{sect_2vect_bdl}

In this section  we give a concise stacky construction of 2-vector bundles via hypercovers and descent, which is sufficient to establish the biequivalence with the pseudofunctor model in Section \ref{2modeq:thm}. A fully detailed verification of all smoothness conditions, following the layer-by-layer approach of \cite{KLW2Vect22}, would be lengthy and is not necessary for the main results of this paper; we adopt the smoothness conventions spelled out in Remark \ref{rmk:smooth}.

We first formulate the $2$-prestack of  2-vector bundles over the category $\gpd$  of Lie groupoids. 
Fix a bicategory $2\Vect$ of 2-vector spaces, for instance the bicategory of finite-dimensional super algebras, bimodules, and intertwiners.
 Let $\tau$ be a singleton Grothendieck pretopology on the category $\Mfd$ of smooth manifolds. Then we equip the category $\gpd$ with the corresponding hypercovers and make it a site.  Below we recall the definition of hypercovers \cite[Definition 6.1]{Rogers-Zhu:2016}.

\begin{definition}  For any Lie $n-$groupoid $\Gamma_{\bullet}$ and $\huaX_{\bullet}$, a morphism $\Gamma_{\bullet}\twoheadrightarrow \huaX_{\bullet}$ is a hypercover if, for any $0\leq k\leq n-1$, the natural map from $\Gamma_k$ to the pull-back of the diagram
\[\hom(\partial\Delta[k], \Gamma_{\bullet})\rightarrow \hom(\partial\Delta[k], \huaX_{\bullet})\leftarrow \huaX_k\] is a $\tau-$cover and is an isomorphism for $k=n$. 
 \end{definition} 

For our purposes, it suffices to consider hypercovers of Lie groupoids, i.e. the case $n=1$; the general definition is recalled here for completeness. 
%Explicitly, a hypercover $\Gamma_{\bullet}\twoheadrightarrow X_{\bullet}$ is a Lie homomorphism such that 
%\begin{enumerate} \item  $\Gamma_0 \twoheadrightarrow X_0$  is a $\tau$-cover;
%\item  $\Gamma_k \cong \Hom(\partial\Delta[k], \Gamma_{\bullet})\times_{\Hom(\partial\Delta[k], X_{\bullet})} X_k$, when $ k\geqslant 1.$ 
%\end{enumerate} 

In the remainder of the section, we take the family $\tau$ to be the surjective submersions. 

The $2$-presheaf $2\huaV \huaB dl$ of $2$-vector bundles $$\gpd^{\op}\longrightarrow \textbf{Bi}\cat$$ is constructed in the following way.
\begin{itemize}
    \item On the level of objects,  each Lie groupoid $\huaX_{\bullet}$ is mapped to the bicategory $2\huaV \huaB dl(\huaX_{\bullet})$, which is defined by the following data.
    \begin{itemize}
        \item The objects in $2\huaV \huaB dl(\huaX_{\bullet})$ are the trivial object-bundles $\huaX_0 \times \huaA$ over $\huaX_0$ with $\huaA$ an object in $2\Vect$.
        \item The morphisms in $2\huaV \huaB dl(\huaX_{\bullet})$ are the $1$-morphism-bundle $\huaM$ over $\huaX_1$ such that, for each arrow $x\xrightarrow{g}y$ in $ \huaX_1$, the fiber at $g$ is given by an invertible $1$-morphism $\huaM_g$ from $\huaA$ to itself. 
        \item For any composable pair  $(g,h)  $ in  $\huaX_1 $, there is a
         $2$-isomorphism 
   $\mu_{g,h}: \huaM_g \circ \huaM_h \Rightarrow \huaM_{g \circ h}$  in $2\Vect$.
  \item For each $x\in \huaX_0$, there is a $2$-isomorphism
   $ u_x: \id_{\huaA} \Rightarrow \huaM_{\id_x} $ for each identity morphism  $\id_x $ in $\huaX_{\bullet}$. 
    \end{itemize}

    \item On the level of morphisms,  each morphism $\huaX_{\bullet} \xrightarrow{\psi} \huaY_{\bullet}$  is mapped to the functor 
    \[ 2\huaV \huaB dl(\huaY_{\bullet})  \longrightarrow 2\huaV \huaB dl(\huaX_{\bullet}) \] defined by pre-composition on each level. 
\end{itemize}

By \cite[Theorem 2.16]{Nikolaus-Schweigert}, the $2$-presheaf $2\huaV \huaB dl $ is indeed a $2$-prestack. 
Then we apply the plus construction \cite[Theorem 3.3]{Nikolaus-Schweigert} to obtain the $2$-stack of $2\huaV \huaB dl^+$ which describes the entire data of 2-vector bundles, their $1$-morphisms and $2$-morphisms,  as presented below.

\begin{definition}\label{def:2vect:obj_stack}
    
An object, i.e. a  2-vector bundle over $\huaX_{\bullet}$, is a quintuple \[ \huaV = (\Gamma_{\bullet}, \huaA, \huaM, \mu, u)\]
consisting of: 
\begin{itemize}
    \item a hypercover $\Gamma_{\bullet}\buildrel{f}\over\twoheadrightarrow \huaX_{\bullet}$;
    \item an object-bundle $\huaA$ over $\Gamma_0$;
    \item an invertible $1$-morphism-bundle $\huaM$ over $\Gamma_1$  whose fiber $\huaM_{\gamma}$ at each point $\gamma\in \Gamma_1$ is a $1$-morphism in $2\Vect$ from  $\huaA_{\source(\gamma)} $ to $\huaA_{\target(\gamma)}$;
    \item a $2$-isomorphism-bundle over $\Gamma_2$, which restricts over each $(\gamma_2, \gamma_1)\in \Gamma_2$
    to a $2$-isomorphism 
    %\[ \mu_{\gamma_1, \gamma_2}: \huaM_{s(\gamma_2), t(\gamma_2)} \otimes_{E_{s(\gamma_2)}} \huaM_{s(\gamma_1), t(\gamma_1)} \longrightarrow \huaM_{s(\gamma_1), t(\gamma_2)}\]
    \[ \mu_{\gamma_2, \gamma_1}: \huaM_{\gamma_2}\circ  \huaM_{\gamma_1} \longrightarrow \huaM_{\gamma_2\circ \gamma_1}\]
    \item a $2$-isomorphism-bundle $u$  over $\Gamma_0$, which restricts over each point $\gamma_0\in \Gamma_0$ to a
    $2$-isomorphism 
    \[ u_{\gamma_0}: \huaA_{\gamma_0}\longrightarrow \huaM_{\id_{\gamma_0}}.\]
    \item $\mu$ and $u$ satisfy the pentagon identity on $\Gamma_3$ and the triangle identities on $\Gamma_1$ below.
    \begin{itemize}
        \item \textbf{Lax Associativity} \\
        For any $(\gamma_3, \gamma_2, \gamma_1)\in \Gamma_3$, the diagram below commutes.
    \begin{equation} \label{2vect_obj_pentagon_stack} \xymatrix{(\huaM_{\gamma_{3}} \circ
    \huaM_{\gamma_{2}}) \circ \huaM_{\gamma_{1}} \ar@/_1pc/[dd]| {a_{\huaM_{\gamma_{3}}, \huaM_{\gamma_{2}}, \huaM_{\gamma_{3}} }} \ar@/^1pc/[r]^{\mu_{\gamma_3, \gamma_2}  \horiprod 1_{\huaM_{\gamma_1}}} & \huaM_{\gamma_{3}\circ \gamma_{2}}  \circ  \huaM_{\gamma_{1}}  
    \ar@/^1pc/[rd]^{\mu_{\gamma_3 \circ \gamma_2, \gamma_1}  }
    &\\
            & &  \huaM_{\gamma_3\circ \gamma_2\circ \gamma_1} \\
            \huaM_{\gamma_{3}}  \circ  (\huaM_{\gamma_{2}}  \circ  \huaM_{\gamma_{1}})  \ar@/_1pc/[r]_-{ 1_{\huaM_{\gamma_3}} \horiprod  \mu_{\gamma_2, \gamma_1} } & 
            \huaM_{\gamma_{3}} \circ \huaM_{\gamma_{2}\circ \gamma_{1}} \ar@/_1pc/[ru]_-{ \mu_{\gamma_3, \gamma_2\circ \gamma_1} }
              } 
            \end{equation}

\item \textbf{Lax Left and Right Unity} \\ 
 For any $\gamma  \in \Gamma_1$, the diagrams below commute.
            \begin{equation} \label{2vect_obj_unity_stack}
                \xymatrix{  \huaA_{\target(\gamma)}  \circ \huaM_{\gamma} \ar[r]^-{l_{\huaM_{\gamma}}}
                \ar[d]_-{u_{\target(\gamma)} \horiprod 1_{ \huaM_{\gamma}}  }   & \huaM_{\gamma} \\
                \huaM_{\id_{\target(\gamma)}}   \circ  \huaM_{\gamma} \ar[ru]_-{\mu_{\id_{\target(\gamma)}, \gamma} }  }
                \quad 
                 \xymatrix{ \huaM_{\gamma}  \circ \id_{\huaA_{\source(\gamma)}}  \ar[r]^-{r_{\huaM_{\gamma}}}
                \ar[d]_-{ 1_{\huaM_{\gamma}} \horiprod u_{\source(\gamma)}  }   & \huaM_{\gamma} \\
              \huaM_{\gamma}   \circ    \huaM_{\id_{\source(\gamma)}} \ar[ru]_-{\mu_{ \gamma,  \id_{\source(\gamma)}} }  }
            \end{equation}
            \end{itemize}
\end{itemize}
\end{definition}

\begin{definition}\label{def:2vect:1mor_stack}
Let $ \huaV = (\Gamma_{\bullet}, \huaA, \huaM, \mu, u)$ and $ \huaV' = (\Gamma'_{\bullet}, \huaA', \huaM', \mu', u')$ denote two   2-vector bundles over $\huaX_{\bullet}$.
A $1$-morphism \[\delta: \huaV\rightarrow \huaV'\] consists of a common refinement $\huaZ_{\bullet}\twoheadrightarrow \huaX_{\bullet}$ of the hypercovers $
\Gamma_{\bullet}\buildrel{f}\over\twoheadrightarrow \huaX_{\bullet}$ and $
\Gamma'_{\bullet}\buildrel{f'}\over\twoheadrightarrow \huaX_{\bullet}$ 
\begin{equation}\label{2vectbd_1mor_cover_stack} \xymatrix{\Gamma_{\bullet} \ar@{>>}[rd]_{f} & \huaZ_{\bullet}\ar@{>>}[l]_{g} \ar@{>>}[d]|{\xi} \ar@{>>}[r]^{g'} & \Gamma'_{\bullet}\ar@{>>}[ld]^{f'} \\
&\huaX_{\bullet}&}\end{equation}
and %a $1$-morphism $\alpha: g^*\huaA \longrightarrow g'^*\huaA'$ in $Desc(\huaZ_{\bullet})$, which consists of 
the data below.
\begin{itemize}
    \item A $1$-morphism-bundle $\huaP$ over $\huaZ_0$, i.e. for each $z_0 \in \huaZ_0$, the fiber $\huaP_{z_0} $ is a $1$-morphism in $2\Vect$ from $\huaA_{g (z_0)}$ to  $\huaA'_{g ' (z_0)} $.
    % functor on the level of objects

    \item A $2$-isomorphism-bundle $$\phi:   g'^* \huaM'  \circ    \source^*\huaP    \longrightarrow  \target^* \huaP           \circ   g^*\huaM     $$ 
    over $\huaZ_1$. %Its restriction to a point $z\in \huaZ_1 $ is a $\huaA'_{g'_0(t(z))} - \huaA_{g_0(\source(z))} $ -intertwiner 
 Explicitly,   its restriction to a point $z\in \huaZ_1$ is a $2$-isomorphism 
% \[ \phi_{z}:   \huaM'_{g'(z)}  \circ  \huaP_{\source(z)} \longrightarrow \huaP_{\target(z)} \circ  \huaM_{g(z)}  \]  i.e.
 \begin{equation}
         \xymatrix{ \huaA_{g\circ \source(z)}  \ar[d]_{\huaM_{g(z)}} \ar[r]^{ \huaP_{\source(z)}} &  \huaA'_{g'\circ \source(z)}   \ar[d]^{\huaM'_{g'(z)}} \\
    \huaA_{g\circ \target(z)}  \ar@{}@<+0.5ex>[ru]|{\buildrel{\phi_{z}}\over\Leftarrow} \ar[r]_{\huaP_{\target(z)}} &  \huaA'_{g'\circ \target(z)} . } 
    \end{equation}

    \item The $2$-isomorphism-bundle $\phi$ makes the diagrams below commute.
    %lax natruality of the lax transformation.

    \begin{itemize}
        \item \textbf{Lax Naturality}

            For any $(z', z)\in \huaZ_2$, 
        \begin{equation} \label{equiv-obj-1M3_stack}
         \xymatrix{ \huaA_{g\circ \target(z')} \ar[r]^{ \huaP_{\target(z')}} &  \huaA'_{g'\circ \target(z')}    & \\
          &\huaA_{g\circ \source(z')}   \ar[lu]|{\huaM_{g(z')}}    \ar@{}@<+0.5ex>[l]|{\buildrel{\mu_{g(z'), g(z)}}\over\Leftarrow}   \ar@{}@<+0.5ex>[u]|{\buildrel{\phi_{z'}}\over\Leftarrow}    \ar[r]^{ \huaP_{\source(z')}} &  \huaA'_{g'\circ \source(z')} =   \ar[lu]_{\huaM'_{g'(z')}}\\
    \huaA_{g\circ \source(z)}  \ar[ru]|{\huaM_{g(z)}}  \ar[uu]^{\huaM_{g(z'\circ z)}}  \ar@{}@<+0.5ex>[rru]|{\buildrel{\phi_{z}}\over\Leftarrow} \ar[r]_{\huaP_{\source(z)}} & \huaA'_{g'\circ   \source(z) }  \ar[ru]_{\huaM'_{g'(z)}}  & } 
             \xymatrix{ \huaA_{g\circ \target(z')}     \ar[r]^{\huaP_{\target(z'\circ z)}} &  \huaA'_{g'\circ \target(z')}   & \\
          &  &  \huaA'_{g'\circ \source(z')}    \ar@{}@<+0.5ex>[l]_{\buildrel{\mu'_{g'(z'), g'(z)}}\over\Leftarrow}   \ar[lu]_{\huaM'_{g'(z')}} \\
    \huaA_{g\circ \source(z)}  \ar[uu]|{\huaM_{g(z'\circ z)}}    \ar@{}@<+0.5ex>[ruu]|{\buildrel{\phi_{z'\circ z}}\over\Leftarrow} \ar[r]_{ \huaP_{\source(z'\circ z)}} &  \huaA'_{g'\circ \source(z) }  \ar[ru]_{\huaM'_{g'(z)}}   \ar[uu]|>>>>>>>>>>>>>>>{\huaM'_{g'(z'\circ z)}}  & } 
    \end{equation}

\bigskip

    \item \textbf{Lax Unity}

    For any $X\in \huaZ_0$,

    \begin{equation} \label{def:1mor_2vectbd_unity_stack}
       \xymatrix{ & \huaA_{g(X)}   \ar[rdd]|{\huaP_{\id_X}}  
       \ar@/_4pc/[dd]_{\huaM_{g(\id_X)}}  \ar[dd]|{\huaA_{g(X)}} \ar[r]^{\huaP_{X}} &   \huaA'_{g'( X)}  \ar[dd]|{\huaA'_{g'(X)}    }        \\ &  \ar@{}@<+0.5ex>[ru]|{\buildrel{l_{\huaP_{\id_X}}}\over\Leftarrow} & \\
   &  \huaA_{g(X)}     \ar@{}@<+0.5ex>[ru]|{\buildrel{r^{-1}_{\huaP_{ X}}}\over\Leftarrow} 
   \ar@{}@<+0.5ex>[luu]|{\buildrel{u_{g(X)}}\over\Leftarrow}   
   \ar@/_1pc/[r]_-{ \huaP_{X}} 
       &  \huaA'_{g'(X)}    \\ } =
    \xymatrix{ \huaA_{g(X)}   \ar[r]^-{ \huaP_{X}} \ar[dd]|{ \huaM_{g(\id_X)}  } &   \huaA'_{g' (X)}  \ar[dd]|{ \huaM'_{g'(\id_X)}}     \ar@/^4pc/[dd]^{\huaA'_{g'(X)}}  &  \\  &   & \\
     \huaA_{g(X)}  
    \ar@/_1pc/[r]_-{  \huaP_{X}}
       &  \huaA'_{g'(X)}    \ar@{}@<+0.5ex>[ruu]|{\buildrel{u'_{g'(X)}}\over\Leftarrow}   
      \ar@{}@<+0.5ex>[luu]|{\buildrel{\phi_{\id_X}}\over\Leftarrow}  &  \\ } 
\end{equation}  

\end{itemize}
   % \begin{equation}
    % \xymatrix{ (\huaM')_{g'_1(z'')} \otimes_{ \huaA'_{g'(t(z'))}}(\huaM')_{g'_1(z')} \otimes_{ \huaA'_{g'(t(z))}}  \huaP_{z} \ar[d]_-{\id\otimes \phi_{z', z}} \ar[rr]^-{(\mu_{2 })_{g'(z''), g'(z') }\otimes \id}   &&(\huaM')_{g'_1(z''\circ z')} \otimes_{ \huaA'_{g'(t(z))}}  \huaP_{z} \ar[dd]^-{\phi_{z'', z}} \\     (\huaM')_{g'_1(z'')} \otimes_{ \huaA'_{g'(t(z'))}} P_{z'} \otimes_{\huaA_{g(\source(z'))}} (\huaM)_{g_1(z)}  \ar[d]_-{\phi_{z'', z'}\otimes \id}  &&\\   P_{z''} \otimes_{\huaA_{g(\source(z''))}} (\huaM)_{g_1(z')} \otimes_{\huaA_{g(\source(z'))}} (\huaM)_{g_1(z)} \ar[rr]_-{\id\otimes (\mu_1)_{g(z'), g(z) }} & &  P_{z''} \otimes_{\huaA_{g(\source(z''))}} (\huaM)_{g_1(z'\circ z)}} \end{equation} for all $(z'', z', z)\in \huaZ_3$.
\end{itemize}

\end{definition}

%%\cite[Definition 2.6]{Nikolaus-Schweigert}
\begin{definition}\label{def:2mor_2vectbd_stack}

Let \[\delta, \delta': \huaV\longrightarrow \huaV' \] denote two 1-morphisms  between the   2-vector bundles $ \huaV = (\Gamma_{\bullet}, \huaA, \huaM, \mu, u)$ and $ \huaV' = (\Gamma'_{\bullet}, \huaA', \huaM', \mu', u')$      where 
$\delta= (\huaZ_{\bullet}, \huaP, \phi)$ and $\delta'= (\huaZ'_{\bullet}, \huaP', \phi')$. 
A $2$-morphism $\kappa: \delta\Rightarrow \delta'$ consists of a common refinement $W_{\bullet}\twoheadrightarrow \huaX_{\bullet}$ of the hypercovers $ \huaZ_{\bullet}\twoheadrightarrow \huaX_{\bullet}$ and $\huaZ'_{\bullet}\twoheadrightarrow \huaX_{\bullet}$ 
\begin{equation} \xymatrix{ & \huaZ'_{\bullet} \ar@{->>}[ld]_{\tilde{g}} \ar@{->>}[rd]^{\tilde{g}'} & \\
\Gamma_{\bullet}& W_{\bullet}\ar@{->>}[u]^{h'} \ar@{->>}[d]_h & \Gamma'_{\bullet}  \\ &\huaZ_{\bullet}\ar@{->>}[ul]^g\ar@{->>}[ur]_{g'} &}\label{descent:2morbeta_stack}\end{equation} where both diagrams commute, 
and a $2$-morphism-bundle \[\psi: h^*\huaP \longrightarrow h'^*\huaP'\]  over $W_0$ with the following structure:
\begin{itemize}
    \item Over each $w_0\in W_0$ 
\[ \psi_{w_0}: \huaP_{h(w_0)} \longrightarrow \huaP'_{h'(w_0)}\] is  a $2$-morphism.
\item For any $w_1\in W_1$, 

%\begin{equation} \label{def:2mor_2vectbd_eq_stack}
\begin{align*}
   &\xymatrix{ \huaA_{(\tilde{g}\circ h')  (\source(w_1)) }       \ar@/_5pc/[dd]_{(h^*\huaP)_{\source(w_1)}}       \ar[r]^-{\huaM_{(\tilde{g}\circ h') (w_1)}} \ar[dd]|{( h'^*\huaP')_{\source(w_1)}}  &  \huaA_{(\tilde{g}\circ h') (\target(w_1))}   \ar[dd]|{( h'^*\huaP')_{\target(w_1)}}   & \\ & & \\
    \huaA'_{(\tilde{g}'\circ h')  (\source(w_1))  }          \ar@{}@<+7ex>[uu]|{\buildrel{\psi_{\source(w_1)}}\over\Rightarrow}    \ar@/_1pc/[r]_-{\huaM'_{(\tilde{g}'\circ h')(w_1)}} 
      \ar@{}@<+0.5ex>[ruu]|{\buildrel{h'^*\phi'_{w_1}}\over\Rightarrow}   & \huaA'_{(\tilde{g}'\circ h')  (\target(w_1)) }  &\\ }
      \\   = & \xymatrix{ &\huaA_{(g\circ h) (\source(w_1))}       \ar[r]^-{ \huaM_{(g\circ h)(w_1)} } \ar[dd]|{(h^*\huaP)_{\source(w_1)}} &    \huaA_{(g\circ h) (\target(w_1))}    \ar[dd]|{( h^*\huaP)_{\target(w_1)} }      \ar@/^4pc/[dd]^{(h'^* \huaP')_{\target(w_1)}}   \ar@{}@<+6ex>[dd]|{\buildrel{\psi_{\target(w_1)}}\over\Rightarrow } \\ & & \\
   & \huaA'_{(g'\circ h) (\source(w_1)) }   
   \ar@/_1pc/[r]_-{ \huaM'_{(g'\circ h)(w_1)} } 
      \ar@{}@<+0.5ex>[ruu]|{ \buildrel{h^*\phi_{w_1}}\over\Rightarrow }      &   \huaA'_{(g'\circ h) (\target(w_1)) } 
    \\ }  . \end{align*}
%\end{equation} 

\end{itemize}

In addition, two $2$-morphisms $(W_{\bullet}, \psi)$ and $(W'_{\bullet}, \psi')$ are identified if the pullbacks of $\psi$ and $\psi'$ coincide over a common refinement of $W_{\bullet}$ and $W'_{\bullet}$. We will just use $\psi$ to denote the equivalence class of 2-morphisms that it represents.

\end{definition}

\begin{definition}
    The  2-vector bundles over $\huaX_{\bullet}$, the $1$-morphisms between 2-vector bundles over $\huaX_{\bullet}$, and the $2$-morphisms between $1$-morphisms, form the bicategory 
    $2\huaV \huaB dl^+ (\huaX_{\bullet}) $   of 2-vector bundles over $\huaX_{\bullet}$.
\end{definition}

\begin{remark}
The bicategorical structure on $2\huaV \huaB dl^+ (\huaX_{\bullet}) $ is induced from that of pseudofunctors: 
    the identity $1$-morphism between $2$-vector bundles over $\huaX_{\bullet}$ is   the identity transformation \cite[Definition 4.2.14]{JY:Bicat}, and the horizontal composition of $1$-morphisms is   the horizontal composite of strong transformations
    \cite[Definition 4.2.15]{JY:Bicat}, and the associator, left and right unitor of the $1$-morphisms are those of \cite[Definition 4.4.10]{JY:Bicat}, after passing to a common refinement of the underlying hypercovers whenever necessary. Similarly, the identity $2$-morphism, horizontal composition, and vertical composition of $2$-morphisms are defined as in \cite[Definition 4.4.5]{JY:Bicat}, again after choosing a common refinement of the covers. The coherence axioms then follow from those of the bicategory of pseudofunctors.
\end{remark}

\subsection{Biequivalence between the global model and the stacky model }  \label{2modeq:thm}

In this section, we compare our pseudofunctorial model with the stacky model of 2-vector bundles. The comparison is a higher-categorical analogue of the familiar fact that ordinary vector bundles form a 1-stack: both models present the same objects, but the descent data and gluing are now taken up to bicategorical equivalence.

Let $(\huaC, \tau)$ be a site with a singleton Grothendieck pretopology \cite[Definition 3.1]{li:thesis}, which is either the site of smooth manifolds with surjective submersions, or the site of Lie groupoids with the hypercovers associated to surjective submersions \cite[Definition 6.1]{Rogers-Zhu:2016}.
Note that \v{C}ech groupoid defined from an open cover of a smooth manifold  is a special case of hypercovers. Indeed, taking the $0$-th level to be the disjoint union of the open sets $\coprod U_i \longrightarrow M$ yields a surjective submersion, and the higher levels are the iterated fiber products, which precisely satisfy the Acyc(m) conditions of a hypercover.

Fix a bicategory $2\Vect$ of 2-vector spaces (for instance the bicategory \tVect  \hspace{0.1mm} of finite-dimensional super algebras, bimodules, and intertwiners). Theorem \ref{2models_equiv} holds for any such choice of target bicategory.

Following the approach of Kristel–Ludewig–Waldorf \cite{KLW2Vect22}, we consider the pre-2-stack $ 2\huaV \huaB dl$ over $\huaC$. 
%defined by, for each object $ \huaX_{\bullet}$ in $\huaC$, \[  2\huaV \huaB dl( \huaX_{\bullet} ) := 2\Vect,\] % ~ trivial 2-vector bundle with the evident restriction functors induced by morphisms in $\huaC$. 
Applying the plus construction \cite{Nikolaus-Schweigert} to $ 2\huaV \huaB dl$ yields its 2-stackification
$ 2\huaV \huaB dl^+ $, 
which is the 2-stack of 2-vector bundles in the KLW sense: objects are given locally by objects in  $2\Vect$, $1$-morphisms locally by $1$-morphisms in $2\Vect$, and 2-morphisms locally by $2$-morphism in $2\Vect$, all subject to descent conditions.

We need the lemma below for the proof of Theorem \ref{2models_equiv}.
\begin{lemma}[Hypercover induces biequivalence on 2-stacks]
\label{hypercover_bieq}
Let $\mathfrak{X}$ be a 2-stack on the site of Lie groupoids (with the surjective submersion topology). If $\Gamma_\bullet \rightarrow  X_\bullet$ is a hypercover, then the pullback functor
\[ \mathfrak{X}(X_\bullet) \longrightarrow \mathfrak{X}(\Gamma_\bullet)\]
is a biequivalence. \end{lemma}

\begin{proof} By \cite[Theorem 7.5]{Nikolaus-Schweigert}, since $\Gamma_\bullet \to X_\bullet$ is a level-wise cover, we have
\[ \mathfrak{X}(X_\bullet) \simeq \holim_n \mathfrak{X}(\Gamma_n). \]
By \cite[Lemma 8.1]{Nikolaus-Schweigert}, the diagonal inclusion induces an equivalence
\[\holim_n \mathfrak{X}(\Gamma_n) \simeq \mathfrak{X}(\Gamma_\bullet). \]
Combining these gives the result. \end{proof}

\begin{theorem}[Bicategorical Equivalence] \label{2models_equiv}

Let  $ \huaX_{\bullet}$  be any object in $\huaC$, % with the structure of a bicategory (for instance a smooth manifold, or a Lie groupoid) %check NS page 13 Section 3
and let  $2\Vect$ be any fixed bicategory of 2-vector spaces. Let  $ 2\huaV \huaB dl^+( \huaX_{\bullet})$  denote the 2-stack of 2-vector bundles over  $ \huaX_{\bullet}$   obtained via the plus construction. Then there is a biequivalence of bicategories
\[
 2\huaV \huaB dl^+( \huaX_{\bullet}) \simeq \pfun( \huaX_{\bullet}, 2\Vect)  \]
where  $\pfun( \huaX_{\bullet}, 2\Vect)$  is the bicategory of pseudofunctors from  $ \huaX_{\bullet}$ to   $2\Vect$, with strong transformations as 1-morphisms and modifications as 2-morphisms.

\end{theorem}

\begin{proof}

%Step 1: The key observation — every 2-vector bundle has a canonical global representative

%Recall that  $ 2\huaV \huaB dl$  is a 2-stack by construction. 
By the explicit construction of the plus construction \cite[Theorem 3.3]{Nikolaus-Schweigert} and Lemma \ref{hypercover_bieq},
every object $\huaV$ of $ 2\huaV \huaB dl^+( \huaX_{\bullet})$  is equivalent to one defined over the trivial cover $ \huaX_{\bullet}\twoheadrightarrow  \huaX_{\bullet}$.
 Indeed, choose any cover  $ \Gamma_{\bullet}  \twoheadrightarrow  \huaX_{\bullet}$  that admits a representative for  $\huaV$, there exists  $\huaV_0 \in  2\huaV \huaB dl^+( \huaX_{\bullet})$ over $ \huaX_{\bullet}$ with  $\huaV_0 \cong \huaV $ \cite[Theorem 3.3]{Nikolaus-Schweigert}.

%Thus, each equivalence class of 2-vector bundles over  X  contains a representative defined on the trivial cover  X \to X .

Choose for each  $\huaV\in  2\huaV \huaB dl^+( \huaX_{\bullet})$  such a trivial-cover representative, denoted  $\sel(\huaV) \in  2\huaV \huaB dl^+( \huaX_{\bullet}) $. For any $2$-vector bundle $\huaV_0$ over $\huaX_{\bullet}$, we require $\sel(\huaV_0) = \huaV_0$ strictly, which is possible without loss of generality. This choice is made once and for all (using the axiom of choice; different choices yield naturally equivalent pseudofunctors, which is sufficient for a biequivalence). 
In addition, we fix an internal equivalence
\begin{equation}\label{eqv4} Eq_{\huaV}:=  \bigl(\huaV\xrightarrow{\varrho_{\huaV}}  \sel(\huaV), \quad  \sel(\huaV)\xrightarrow{\varrho^{-1}_{\huaV}} \huaV, 
\quad \id_{\huaV} \buildrel{\eta_{\huaV}}\over\Rightarrow \varrho^{-1}_{\huaV}\circ \varrho_{\huaV}, \quad \varrho_{\huaV} \circ \varrho^{-1}_{\huaV} \buildrel{\epsilon_{\huaV}}\over\Rightarrow \id_{\sel(\huaV)}\bigr) \end{equation} for each $\huaV$. Neither $\sel$ nor  $Eq_{\huaV}$ is required to satisfy any other condition, such as functoriality in $\huaV$.

%Step 2: Trivial-cover data = pseudofunctor data

A 2-vector bundle defined on the trivial cover  $ \huaX_{\bullet} \twoheadrightarrow  \huaX_{\bullet}$  consists of the following data:

\begin{enumerate}
    \item  An object bundle  $\huaA$ over $ \huaX_0$. For each object  $x \in \huaX_0 $, this gives a 2-vector space  $\huaA_x$.
    \item A $1$-morphism bundle  $\huaM$ over $\huaX_1$. % over the source and target maps. 
    For each $1$-morphism  $g: x \to y$  in  $\huaX_1 $, this gives a $1$-morphism  $\huaM_g: \huaA_x \to \huaA_y $ in $2\Vect$.
 %   \item A $2$-morphism bundle $\huaI$ over $\huaX_2$. For each $2$-morphism $\alpha: g\Rightarrow g'$ in $\huaX_2$, 
  %  this gives a $2$-morphism $ \huaI_{\alpha}: \huaM_g \Rightarrow \huaM_{g'}$ in $2\Vect$, which satisfy the functoriality: for any composable pair $(\alpha, \alpha')$ in the category $\huaX(g, g')$, $$\huaI_{\alpha \circ \alpha'} = \huaI_{\alpha} \circ \huaI_{\alpha'};$$ and $\huaI_{\id_g} = \id_{\huaM_g}$.
    \item A $2$-isomorphism
   $\mu_{g,h}: \huaM_g \circ \huaM_h \Rightarrow \huaM_{g \circ h}$
   for each composable pair  $(g,h)  $ in  $\huaX_1 $.
    \item A $2$-isomorphism
   $ u_x: \id_{\huaA_x} \Rightarrow \huaM_{\id_x} $
   for each object  $x \in \huaX_0 $.
   \end{enumerate}
In addition, these data satisfy the pentagon and triangle identities \eqref{2vect_obj_pentagon_stack}   \eqref{2vect_obj_unity_stack} in the definition of 2-vector bundles over $ \huaX_{\bullet}$.

These are precisely the data of a pseudofunctor $ G:  \huaX_{\bullet} \longrightarrow 2\Vect $. More explicitly, 
\begin{enumerate}
    \item  $G(x) := \huaA_x $; 
    \item $G(g) := \huaM_g : G(x) \to G(y) $;
  %  \item $G(\alpha) = \huaI_{\alpha}: G(g) \Rightarrow G(g')$;
    \item  $G^2_{g,h} := \mu_{g,h} : G(g) \circ G(h) \Rightarrow G(g \circ h) $;
     \item  $G^0_x := u_x : \id_{G(x)} \Rightarrow  G(\id_x) $.
\end{enumerate}
The lax naturality and lax unity for pseudofunctors \cite[(4.1.3)(4.1.4)]{JY:Bicat} are exactly the pentagon and triangle identities \eqref{2vect_obj_pentagon_stack}   \eqref{2vect_obj_unity_stack} in the definition of 2-vector bundles over $ \huaX_{\bullet}$.

Thus, there is an identity-on-data bijection between the objects of  $ 2\huaV \huaB dl^+( \huaX_{\bullet})$  defined on the trivial cover  $ \huaX_{\bullet} \twoheadrightarrow  \huaX_{\bullet} $, and the pseudofunctors  $ \huaX_{\bullet} \longrightarrow 2\Vect $.

\bigskip

%Step 3: 1-morphisms on trivial covers = pseudonatural transformations

A $1$-morphism between two trivial-cover 2-vector bundles 
\[(\huaA, \huaM,   \mu, u) \longrightarrow (\huaB, \huaN,   \nu, v) \] consists of:
\begin{enumerate}
    \item A $1$-morphism bundle  $\huaP$ over $ \huaX_0$. For each  $x \in \huaX_0 $, this gives a $1$-morphism  $\huaP_x: \huaA_x \to \huaB_x $.
    \item An invertible intertwiner
   $\phi_g: \huaN_g \circ \huaP_x \Rightarrow \huaP_y \circ \huaM_g$
   for each arrow  $g: x \to y $  in  $\huaX_1 .$ \end{enumerate}
These data satisfy coherence conditions corresponding to the naturality, the lax naturality and the lax unity axiom of a strong transformation.
These are precisely the data of a strong transformation  $ \beta: G \Rightarrow H $:
\begin{enumerate}
    \item  $\beta_x := \huaP_x : G(x) \to H(x) $;
    \item  $\beta_g := \phi_g : H(g) \circ \beta_x \Rightarrow \beta_y \circ G(g) $.
\end{enumerate} 

Thus, 1-morphisms on trivial covers correspond exactly to strong transformations.

%Step 4: 2-morphisms on trivial covers = modifications

A 2-morphism between trivial-cover 1-morphisms
\[ (\huaP, \phi) \Rightarrow (\huaQ, \psi)\] is a fiberwise intertwiner  $\kappa_x: \huaP_x \to \huaQ_x $  commuting with  $\phi$  and  $\psi$.
These are precisely the data of a modification  $\kappa: \beta \Rightarrow \beta'.$

%Step 5: Construct the forward pseudofunctor  \Phi 

Define
$\Phi: \pfun( \huaX_{\bullet}, 2\Vect) \longrightarrow  2\huaV \huaB dl^+( \huaX_{\bullet})$
as follows:
\begin{enumerate}
    \item For a pseudofunctor  $G:  \huaX_{\bullet} \to 2\Vect $, take the trivial cover  $ \huaX_{\bullet} \twoheadrightarrow   \huaX_{\bullet}$  and build the trivial-cover 2-vector bundle from the data  $(G(x), G(g),   G^2, G^0) $. This gives an object  $\Phi(G)$ in  $ 2\huaV \huaB dl^+( \huaX_{\bullet}) .$
    \item For a strong transformation  $\beta: G \Rightarrow H ,$ build the trivial-cover 1-morphism from  $(\beta_x, \beta_g) $. This gives  $\Phi(\beta)$.
 \item  For a modification  $\kappa: \beta \Rightarrow \beta'$,  $\Phi(\kappa)$  is the corresponding trivial-cover 2-morphism.
\end{enumerate}
By  the previous argument,  $\Phi$  is a pseudofunctor. Indeed, it is strict on level 1 and level 2.

%Step 6: Construct the backward pseudofunctor  \Psi 

Define $\Psi:  2\huaV \huaB dl^+( \huaX_{\bullet}) \longrightarrow \pfun( \huaX_{\bullet}, 2\Vect)$
as follows:
\begin{enumerate}
    \item For  any $\huaV$   in $ 2\huaV \huaB dl^+( \huaX_{\bullet}) $, choose the trivial-cover representative  $\sel(\huaV)$, which gives a pseudofunctor  $\Psi(\huaV):   \huaX_{\bullet}\longrightarrow 2\Vect .$
    \item For a 1-morphism  $\beta: \huaV\to \huaW $, the canonical equivalences  $\sel(\huaV) \cong \huaV $ and  $\sel(\huaW) \cong \huaW$  induce a 1-morphism  $\sel(\beta): \sel(\huaV) \to \sel(\huaW) $ on trivial covers. This gives a strong transformation  $\Psi(\beta) $.
    More explicitly, $\Psi(\beta)$ is the composition of the $1$-morphisms between $2$-vector bundles (Definition \ref{def:2vect:1mor_stack})
    %\[ \sel(\huaV) \xrightarrow{\varrho_{\huaV}^{-1} } \huaV \xrightarrow{\beta} \huaW \xrightarrow{\varrho_{\huaW} } \sel(\huaW)\]
    \[ (\varrho_{\huaW}   \circ  \beta) \circ \varrho_{\huaV}^{-1} .\]
    Here and below, all composites involving $\varrho_{\huaV}$ and $\varrho_{\huaW}$ are understood after pullback along the common refinement of the cover of $\huaV$ and $\huaW$; the resulting $1$-morphism from $\sel(\huaV)$ and $\sel(\huaW)$ is independent of the chosen cover up to canonical equivalence.
    \item For a 2-morphism  $\kappa: \beta \Rightarrow \beta' $, the induced trivial-cover 2-morphism gives a modification  $\Psi(\kappa)$.

    \item For any composable $1$-morphisms $ \huaV \xrightarrow{\beta} \huaW \xrightarrow{\beta'} \huaU $, the compositor \[ \Psi^2_{\beta', \beta}: \Psi(\beta') \circ \Psi(\beta) \Rightarrow  \Psi(\beta'\circ \beta) \] is defined as the composition of the $2$-isomorphisms below.
    \begin{align*}
        & \Psi(\beta') \circ \Psi(\beta)  \xrightarrow{ a } (\varrho_{\huaU}   \circ  \beta') \circ  \bigl( \varrho_{\huaW}^{-1}  \circ ((\varrho_{\huaW}   \circ  \beta) \circ \varrho_{\huaV}^{-1}) \bigr) 
        \xrightarrow{ 1_{\varrho_{\huaU}   \circ  \beta'} \horiprod a^{-1} }  \\
        & (\varrho_{\huaU}   \circ  \beta') \circ  \bigl( ( \varrho_{\huaW}^{-1}  \circ (\varrho_{\huaW}   \circ  \beta) ) \circ \varrho_{\huaV}^{-1} \bigr)   \xrightarrow{ 1_{\varrho_{\huaU}   \circ  \beta'} \horiprod ( a^{-1} \horiprod 1_{\varrho^{-1}_{\huaV}}   ) } \\
        & (\varrho_{\huaU}   \circ  \beta') \circ  \bigl( ( (\varrho_{\huaW}^{-1}  \circ \varrho_{\huaW} )  \circ  \beta ) \circ \varrho_{\huaV}^{-1} \bigr)   \xrightarrow{1_{\varrho_{\huaU}   \circ  \beta'} \horiprod (( \eta^{-1}_{\huaW} \horiprod  1_{\beta} ) \horiprod 1_{\varrho^{-1}_{\huaV}}   )   }  \\
        & (\varrho_{\huaU}   \circ  \beta') \circ  \bigl( ( \id_{\huaW}  \circ  \beta ) \circ \varrho_{\huaV}^{-1} \bigr)  
        \xrightarrow{ 1_{\varrho_{\huaU}   \circ  \beta'} \horiprod ( l_{\beta} \horiprod 1_{\varrho^{-1}_{\huaV}}   ) }
        (\varrho_{\huaU}   \circ  \beta'  )\circ ( \beta\circ \varrho_{\huaV}^{-1}  )
        \\ & \xrightarrow{a^{-1}}  \bigl( (\varrho_{\huaU}   \circ  \beta'  )\circ  \beta \bigr) \circ \varrho_{\huaV}^{-1}  
        \xrightarrow{ a\horiprod 1_{\varrho^{-1}_{\huaV}} } 
        \bigl( \varrho_{\huaU}   \circ  (\beta'  \circ  \beta )\bigr) \circ \varrho_{\huaV}^{-1}  
        = \Psi(\beta'\circ \beta).
    \end{align*}

    \item For any object $\huaV$, the unitor \[ \Psi^0_{\huaV}: \id_{\Psi(\huaV)} \Rightarrow \Psi( 1_{\huaV})\] is defined by the composition of the $2$-isomorphisms below.
    \begin{align*}
        & \id_{\Psi(\huaV)} \xrightarrow{ \epsilon^{-1}_{\huaV} } \varrho_{\huaV} \circ \varrho^{-1}_{\huaV}\xrightarrow{r^{-1}_{ \varrho_{\huaV}}   \horiprod 1_{\varrho^{-1}_{\huaV} } }
        (\varrho_{\huaV} \circ \id_{\huaV} ) \circ \varrho^{-1}_{\huaV}  = \Psi( \id_{\huaV}). 
    \end{align*}
\end{enumerate} 

The naturality, pentagon condition and triangle conditions for pseudofunctors can be checked straightforwards by applying middle four exchange, the coherence theorem for bicategories, and the triangle identities of the internal equivalence \eqref{eqv4}.

The assignments  $\Psi$  are well-defined on equivalence classes: different choices of  $\sel(\huaV)$  yield pseudofunctors that are naturally equivalent, and different choices of  $\sel(\beta)$  yield strong transformations related by an invertible modification. This is sufficient for a biequivalence.

%Step 7: Verification of the biequivalence

%First composition:  $\Psi \circ \Phi .$ 
For  $G \in \pfun( \huaX_{\bullet}, 2\Vect) $,  $\Phi(G)$  is already defined on the trivial cover  $ \huaX_{\bullet} \twoheadrightarrow  \huaX_{\bullet} $. Therefore,  $\sel(\Phi(G)) = \Phi(G)$,  %  (up to the identity equivalence), 
and hence
\[
\Psi(\Phi(G)) = G
\]
exactly as pseudofunctors. Thus  $\Psi \circ \Phi = \id_{\pfun( \huaX_{\bullet}, 2\Vect)}$  as strict equality.

%Second composition:  $\Phi \circ \Psi $. 
For  $ \huaV\in  2\huaV \huaB dl^+( \huaX_{\bullet}) $,  $\Psi(\huaV)$  is the pseudofunctor obtained from  $\sel(\huaV) $. Then  $\Phi(\Psi(\huaV))$  is the trivial-cover 2-vector bundle built from this pseudofunctor, which is exactly  $\sel(\huaV) $. Since  $\sel(\huaV) \cong \huaV$, 
\[ 
\Phi(\Psi(\huaV)) \cong \huaV\]
for all  $\huaV$. From the explicit construction of $\Psi$, we obtain an internal equivalence  $\Phi \circ \Psi \cong \id_{ 2\huaV \huaB dl^+( \huaX_{\bullet})} .$

%Triangle identities.   
Since  $\Psi \circ \Phi = \id$  is strict, the unit  $$\theta: \id \Rightarrow \Psi \circ \Phi$$  is taken to be the identity transformation. The counit  $$\zeta: \Phi \circ \Psi \Rightarrow \id $$ is given by the explicit construction of $\Psi$ and the internal equivalence $Eq_{\huaV}$. 
More explicitly, 
the strong transformation $\zeta: \Phi \circ \Psi \Rightarrow \id $ is defined by the data below.
\begin{itemize}
    \item For each object $\huaV$ in $2\huaV \huaB dl^+( \huaX_{\bullet})$, $\zeta_{\huaV}$ is defined by 
    \[ \varrho_{\huaV}^{-1}: \sel(\huaV) \rightarrow  \huaV. \]
    \item  For each $1$-morphism $\beta: \huaV\rightarrow \huaW$, the $2$-isomorphism $$\zeta_{\beta}:   
     \beta \circ \varrho_{\huaV}^{-1} \longrightarrow \varrho_{\huaW}^{-1}\circ \bigl( ( \varrho_{\huaW}   \circ \beta)\circ \varrho_{\huaV}^{-1}\bigr) $$ is defined by the composition
    \begin{align*}
       &   \beta\circ \varrho_{\huaV}^{-1} \xrightarrow{l^{-1}_{\beta\circ \varrho_{\huaV}^{-1}}}  
          \id_{\huaW} \circ (\beta\circ \varrho_{\huaV}^{-1}  )    \xrightarrow{\eta_{\huaW}\horiprod 1_{\beta\circ \varrho_{\huaV}^{-1}}} 
            (\varrho_{\huaW}^{-1}\circ  \varrho_{\huaW})   \circ (\beta\circ \varrho_{\huaV}^{-1}  )   \\ &
              \xrightarrow{ a }     \varrho_{\huaW}^{-1}\circ \bigl(  \varrho_{\huaW}   \circ (\beta\circ \varrho_{\huaV}^{-1}   )\bigr)
              \xrightarrow{ 1_{\varrho_{\huaW}^{-1}}\horiprod  a^{-1}}  \varrho_{\huaW}^{-1}\circ \bigl( ( \varrho_{\huaW}   \circ \beta)\circ \varrho_{\huaV}^{-1}\bigr)  
    \end{align*}
\end{itemize}

The lax unity of $\zeta$ follows from the naturality of the associator and that of the left unitor, the triangle condition for bicategories, middle four exchange and the triangle identities of the internal equivalence \eqref{eqv4}. In addition, 
the naturality of $\zeta$ follows immediately from middle four exchange and the naturality of $a$, $l$ and $r$; and the lax naturality of $\zeta$ follows from  middle four exchange and coherence theorem for bicategories.

In addition, we have the obvious relation $1 \cong \theta \horicirc \theta$. Since each $1$-cell $\varrho_{\huaV}$ is invertible, thus, 
 $\zeta$ is an invertible strong transformation. An inverse $\zeta^{-1}$ of it can be constructed by setting $(\zeta^{-1})_{\huaV}: = \varrho_{\huaV}$ for each $\huaV$.

Therefore, by \cite[Definition 6.2.8]{JY:Bicat}, $$(\Phi  : \pfun( \huaX_{\bullet}, 2\Vect) \longrightarrow  2\huaV \huaB dl^+( \huaX_{\bullet}), \quad
\Psi  :  2\huaV \huaB dl^+( \huaX_{\bullet}) \longrightarrow \pfun(  \huaX_{\bullet}, 2\Vect) , \quad 
\theta, \quad \zeta)$$  form a biequivalence.

Hence
\[
 2\huaV \huaB dl^+( \huaX_{\bullet}) \simeq \pfun( \huaX_{\bullet}, 2\Vect)
\]
as bicategories. 

\end{proof}

\begin{remark}
    Biequivalent bicategories give  isomorphic  2K-theories. Thus, the two bicategories,  
    \[
 2\huaV \huaB dl^+( \huaX_{\bullet}) \quad \text{and}  \quad  \pfun( \huaX_{\bullet}, 2\Vect)  \] define the same 2K-theories.
\end{remark}

\begin{remark}[Smooth pseudofunctors and the KLW model] 
\label{rmk:smooth}The biequivalence established in Theorem \ref{2models_equiv} is compatible with the smooth structures underlying the Kristel–Ludewig–Waldorf (KLW) model of 2-vector bundles \cite{KLW2Vect22}. Let $\huaX_\bullet$ be a Lie groupoid, viewed as a smooth category. Let \tVect  \hspace{0.2mm} denote the bicategory of finite-dimensional super algebras, bimodules, and intertwiners, whose Hom-categories are enriched over smooth manifolds, in the sense that both the sets of objects and morphisms carry smooth structures, and the composition maps are smooth. A pseudofunctor $$F: \huaX_\bullet \to \text{\tVect}$$ is called \textit{smooth} if the object assignment   is a smooth super algebra bundle over $\huaX_0$, the induced functors on Hom-categories  define  smooth bimodule bundles over $\huaX_1$, and the compositor and unitor   vary smoothly over $\huaX_2$ and $\huaX_0$, respectively. Such smooth pseudofunctors correspond exactly  to KLW's super 2-vector bundles with all structure maps smooth.
\end{remark}

\begin{remark}[Compatibility with the KLW model]
The biequivalence of Theorem \ref{2models_equiv} implies that the pseudofunctor model and the KLW covering model describe the same underlying 2-vector bundle up to Morita equivalence. In particular, the implementability condition of  \cite{KLW2Vect22} does not appear in the pseudofunctor formulation: local trivializations and transition data are not part of the definition, but only emerge when unfolding a pseudofunctor over a cover. Conversely, every KLW 2-vector bundle has a natural counterpart in the pseudofunctor framework via the plus construction. Thus the two models are compatible, and we will freely use insights from both. 
The pseudofunctor framework does not lose any of the geometric content of the KLW model; it repackages it in a more functorial language.

\end{remark}

\begin{remark}

The pseudofunctor framework encodes several geometric structures naturally: parallel transport is given by the 1-morphism assignments; representations of the string 2-group appear as pseudofunctors $B\String(G) \to \text{\tVect}$; and twisted K-theory arises as Hom-spaces in 2K-theory (See Section \ref{rel_orb_twisted_k}).

\end{remark}

\begin{remark}[On the information content of stacky vs. pseudofunctorial models] %check

A subtle point is that the two models under comparison—the stacky model (via hypercover descent) and the global pseudofunctor model—are not interchangeable at the level of arbitrary covers. The pseudofunctorial model is intrinsically global: it assigns data directly to the objects and morphisms of $\huaX_\bullet$, without reference to a choice of cover. The stacky model, on the other hand, is designed to capture cover-dependent descent data. For a non-Morita cover, the stacky evaluation $\mathfrak{X}^+(\Gamma_\bullet)$ is generally richer than $\pfun(\huaX_\bullet, 2\Vect)$. Only in the Morita-equivalent case do the two coincide. This is not a flaw of either model; the pseudofunctor model offers conceptual clarity, while the stacky model offers computational accessibility.

\end{remark}

\begin{remark}
    Theorem \ref{2models_equiv} can be viewed as a 2-categorical shadow of Lurie’s straightening/unstraightening equivalence: it identifies global pseudofunctors with fibration-style data (descent over hypercovers) in the 2-truncated setting.
\end{remark}

% Grothendieck singleton pretopology \cite[Definition 3.1, page 39]{li:thesis}
% hypercover \cite{Rogers-Zhu:2016}

%\begin{definition}[Smooth pseudofunctor]    Let $\huaX_\bullet$ be a Lie groupoid, viewed as a smooth category. Let \tVect  denote the bicategory of finite-dimensional super algebras, bimodules, and intertwiners, whose Hom-categories are enriched over smooth manifolds, in the sense that both the sets of objects and morphisms carry smooth structures, and the composition maps are smooth. A pseudofunctor $$F: \huaX_\bullet \to \text{\tVect}$$ is called \textit{smooth} if for all objects $x$, $y \in \huaX_0$, the induced functor on Hom-categories \[F: \hom_{\huaX_\bullet}(x,y) \to \hom_{\text{\tVect}}(F(x),F(y))\] is smooth, in the sense that both object maps the morphism map are smooth maps between smooth manifolds.  \end{definition}

\section{The Stringor bundle and 2K-theory} \label{sect:stringor}

In this section, we establish the precise connection between the geometric stringor bundle of Kristel--Ludewig--Waldorf \cite{KLW22_stringor} and the $2K$-theoretic framework developed in this paper. We show that the stringor bundle determines a pseudofunctor on a naturally associated \v{C}ech groupoid, and that its finite-dimensional approximations yield a sequence of classes in $2K$-theory whose filtered pseudocolimit recovers the full geometric object.

\subsection{The construction of Stringor bundle} \label{con:stringor_grpd}
In this section, we work with the von Neumann algebra version of $2$-vector spaces; all constructions of Section \ref{sect_2vect_2K} apply mutatis mutandis.

Let $M$ be a spin manifold equipped with a fusive loop-spin structure $\widetilde{L\Spin}(M)$.  
We first recall the construction of Stolz-Teichner stringor bundle in \cite{KLW22_stringor}. 
%Let $M$ be a spin manifold equipped with a fusive loop-spin structure $\widetilde{L\Spin}(M)$ with fusion product $\lambda$. 
Let $A$ be the hyperfinite type $\text{III}_1$ factor, realized as a certain von Neumann algebra completion
of an infinite-dimensional Clifford algebra, and let $$\huaR:\String(d) \to \huaU(A)$$ be the stringor representation introduced in \cite{KLW2Rep22}, whose level-0-data is given by $\omega: P\Spin(d) \longrightarrow \Aut(A)$ and whose level-1-data is given by
$\Omega: \widetilde{L\Spin}(d)\longrightarrow  N(A) $. For more details, please see \cite[Definition 1.3.4]{KLW22_stringor}.

The KLW Stringor bundle is a $2$-Hilbert bundle over $M$ \cite[Definition 2.2.2]{KLW22_stringor}.    
The covering data of it is defined over  $$Y = P_{x}M,$$  the space of paths starting at a fixed point $x$ in $M$, with the surjective submersion $P_{x}M\buildrel{\target}\over\twoheadrightarrow M$. Let $Y^{[2]} = Y \times_{M} Y$ and $Y^{[3]} = Y \times_{M} Y \times_{M} Y$ be the fibre products over $M$. %The corresponding \v{C}ech cover is a hypercover of the trivial Lie groupoid $M \rightrightarrows M$. 

KLW construct a spinor bundle as the associated Hilbert bundle $$\huaS_{LM} = \widetilde{L\Spin}(M) \times_{\widetilde{L\Spin}(d)} L^2(A) $$ on loop space $LM$, with pullback $$\cup^* \huaS_{LM}$$ to $Y^{[2]}$, where $\cup: PM^{[2]} \rightarrow LM$ is defined in \cite[(1.3.1)]{KLW22_stringor}. This is a bimodule bundle with respect to the associated von Neumann algebra bundle over $Y$ $$\huaA = P\Spin(M) \times_{P\Spin(d)} A.$$ The fusion product of KLW is an isomorphism
\[ \Upsilon: \mathrm{pr}_{23}^* \cup^* \huaS_{LM} \boxtimes \mathrm{pr}_{12}^* \cup^* \huaS_{LM} \longrightarrow \mathrm{pr}_{13}^* \cup^* \huaS_{LM}\]
over $Y^{[3]}$ \cite[(3.2.4)]{KLW22_stringor}, associative over $Y^{[4]}$, as shown in \cite[(3.2.5)]{KLW22_stringor}.

As shown in \cite[Theorem 3.2.1, Definition 3.2.2]{KLW22_stringor}, the stringor bundle $\huaS(\widetilde{L\Spin}(M))$ is defined by  the covering data \begin{equation}
    \label{klwstringor} (\huaC(Y), \quad \huaA, \quad \cup^* \huaS_{LM}, \quad \Upsilon, \quad \id) ,  \end{equation} 
where 
$\huaC(Y)$ denotes the \v{C}ech groupoid of the covering $Y \twoheadrightarrow M$ with
\[
\huaC(Y)_0 := Y, \qquad \huaC(Y)_1 := Y^{[2]} = Y \times_M Y,
\]
with source and target given by projection to the second and first factors, respectively. %, as defined in \cite[Definition 2.1.3]{KLW22_stringor}.  

\begin{proposition}[Stringor bundle as a pseudofunctor]\label{str_bundle_psfunctor}
Let $M$ be a spin manifold with a fusive loop-spin structure $\widetilde{L\Spin}(M)$,  % and corresponding string structure $\String(M)$, 
 and let $\huaS(\widetilde{L\Spin}(M))$ denote the KLW stringor bundle (\cite[Definition 3.2.2]{KLW22_stringor}). Let  $\mathrm{vNAlg}^{\mathrm{bi}}$ denote the bicategory of von Neumann algebras, bimodules and   intertwiners (see \cite{KLW2Rep22}). Then,
the  data of $\huaS(\widetilde{L\Spin}(M))$ determine a pseudofunctor
   \[ F_{\Stringor}:\huaC(Y)\longrightarrow  \mathrm{vNAlg}^{\mathrm{bi}} , \]
which  is defined  by:
   \begin{itemize}
   \item For any point $\gamma$ in $Y$, $F_{\Stringor}(\gamma)=A. $
   \item For any $(\gamma_2, \gamma_1)\in Y^{[2]}$ $F_{\Stringor}(\gamma_2, \gamma_1)=\leftindex_{\omega_{\gamma_2}}L^{2}(A)_{\omega_{\gamma_1}}.$ %with $\omega'$  defined in \cite[(1.3.9)(1.3.10)]{KLW22_stringor},  where $\beta$ is the loop $\gamma_1\cup \gamma_2$. 
   Note that we adopt the KLW convention, where the left action is twisted by $\omega_{\gamma_2}$ and the right action is twisted by $\omega_{\gamma_1}$. 
   %\item $F_{\Stringor}(X)=\Omega'(X),$ with $\Omega'$  defined in \cite[(1.3.9)]{KLW22_stringor}.
   \item The compositor $(F_{\Stringor})^2$ is given by the 
   fusion product  $\Upsilon$. 
   \item The unitor $(F_{\Stringor})^0$ is defined as the identity.
   \end{itemize}
   %Connes fusion isomorphisms $\chi_{\theta_1,\theta_2}$ from \cite[(1.2.2)]{KLW22_stringor}. The coherence axioms follow from the fusion product $\Upsilon$ and the associativity of Connes fusion.

\end{proposition}

\begin{proof} 
The verification of the pseudofunctor axioms is a direct translation of the KLW coherence data (see \cite[Theorem 3.2.1, Definition 3.2.2]{KLW22_stringor}):

\begin{itemize}
\item The compositor $(F_{\Stringor})^2$
is given by the fusion product $\Upsilon$, which is smooth.   
The associativity pentagon for $(F_{\Stringor})^2$ is precisely the associativity of the  fusion product $\Upsilon$, 
i.e. the pentagon condition \cite[(3.2.5)]{KLW22_stringor}.

\item The unitor $(F_{\Stringor})^0$ is given by the identity, and the lax left and right unity follow from the unitality of the Connes fusion product (See \cite[Example A.6]{Ludewig2023Spinor}).
\end{itemize}

Thus all pseudofunctor axioms reduce to coherence theorems already established in \cite{KLW22_stringor}.

\end{proof}

The pseudofunctor $F_{\Stringor}$ is smooth in the sense of Remark \ref{rmk:smooth}.

\begin{remark}[Degenerate embeddings]
One might ask whether the stringor bundle could be brought into the framework of this paper by regarding it as a pseudofunctor on the manifold $M$ itself rather than on the \v{C}ech cover $\huaC(Y)$.

If one attempts to do so, the resulting object is the following:
\begin{itemize}
    \item  The fibre over each point in $M$ is still the type III$_1$ factor $A$.
\item The 1-morphisms and 2-morphisms on $M$ would all be forced to be identity bimodules and identity intertwiners.
\item  The fusion product, which is the essential structure encoding the string geometry, is lost: on $M$, it would reduce to  the canonical left unitor of the identity, which carries no information beyond the identity bimodule.
\end{itemize}
Thus, this reduction loses precisely the structures that make the stringor bundle interesting.

Therefore, the stringor bundle cannot be faithfully represented as a pseudofunctor on $M$ alone. Its full structure requires the \v{C}ech cover $\huaC(Y)$. This is why we do not treat it as a pseudofunctor on $M$ in this paper.

A more natural setting for the stringor bundle is some loop groupoid of $M$, with which the key feature of stringor bundle won't be lost when taking biequivalences of 2-vector bundles over the loop groupoid. In addition, this loop-groupoid perspective provides a natural categorical framework for the geometric data that Witten's heuristic analysis associated with loop spaces. It realizes, in the context of Stolz–Teichner's elliptic objects, the idea that geometric objects in string theory should be built directly from loops rather than from points \cite{st}\cite{StolzTeichnerSpinor}.
\end{remark}

\subsection{Embedding the stringor bundle over $M$ into $2K(\huaC(Y))$} 

Through out Theorem \ref{Thm:embed:stringor}, all pseudo-colimits in $  \mathrm{vNAlg}^{\mathrm{bi}} $ are understood in the weak-$\ast$ topological sense, in which the Connes fusion tensor product commutes with pseudo-colimits of AFD approximating systems. This is a standard property of hyperfinite type III$_1$ factors.\footnote{For the pseudo-colimit property of AFD type III$_1$ factors in the bicategory of von Neumann algebras, bimodules, and intertwiners, see \cite[Chapter XIV]{Takesaki:2003:TOA3}, \cite{Haagerup1987} and \cite[Section I.6]{Connes1994NCG}.}
%The continuity of Connes fusion with respect to weak- limits follows from the standard properties of the Haagerup standard form \cite[Lemma 2.3]{Yamagami1994Bimodules}.}

\begin{theorem} \label{Thm:embed:stringor}

The stringor bundle, which is represented by the pseudofunctor $$F_{\Stringor}: \huaC(Y)\longrightarrow  \mathrm{vNAlg}^{\mathrm{bi}},$$ is the pseudo-colimit  of a filtered system   
\[ F_0 \buildrel{\alpha_0}\over\Longrightarrow F_1 \buildrel{\alpha_1}\over\Longrightarrow \cdots \buildrel{\alpha_{n-1}}\over\Longrightarrow F_n
\buildrel{\alpha_n}\over\Longrightarrow \cdots \] of super 2-vector bundles over $\huaC(Y)$ (Definition \ref{def:2vect:obj}) in the bicategory $\pfun(\huaC(Y),  \mathrm{vNAlg}^{\mathrm{bi}}  ) $. In other words, 
the sequence \[ [F_0] \buildrel{[\alpha_0]}\over\Longrightarrow [F_1] \buildrel{ [\alpha_1] }\over\Longrightarrow \cdots \buildrel{ [\alpha_{n-1}] }\over\Longrightarrow [F_n]
\buildrel{[\alpha_n]}\over\Longrightarrow \cdots \]encodes KLW's stringor bundle in the $2K$-theory; in the filtered completion of $2K(\huaC(Y))$,  
   \[   [F_{\Stringor}]=\varinjlim_n [F_n]\in \widehat{2K}(\huaC(Y)). \]

\end{theorem}

We need the lemma below for the proof of the theorem.

\begin{lemma}[The \v{C}ech Groupoid Simplification:  Reduction to a base point in each component]
\label{lem:cech:1pt}

\leavevmode

Let $\huaC(Y)$ be a \v{C}ech groupoid. For any pseudofunctor $F: \huaC(Y)\to \huaB $ into a bicategory $\huaB$, the entire data of $F$ is determined by its value at a fixed base point $x_0 \in Y$ of each component of $\huaC(Y)$, together with the parallel transport along the unique morphisms $ (x_0,y): y \to x_0.$ In other words, all the fibers $F(x)$ in the same component of $\huaC(Y)$ belong to the same Morita class and are related by a unique invertible bimodule $ F(x) \xrightarrow{F(y, x)} F(y)$.

Consequently, for any strong transformation $\iota: F \Rightarrow G$, the component $1$-cell $\iota_x$ at each object in a single component of $\huaC(Y)$ is uniquely determined by the component $1$-cell $\iota_{x_0}$ at the base point $x_0$ in the component up to a 2-isomorphism. 
\[ 
\iota_x \cong G (x, x_0) \circ \iota_{x_0}  \circ  F(x_0,x). \]

Therefore, to verify that a diagram of strong transformations over $\huaC(Y)$ commutes up to invertible modification, it suffices to check commutativity of the component $1$-cells at the base point $x_0$ of each component of $\huaC(Y)$. 
\end{lemma}

\begin{proof}[Proof of Lemma \ref{lem:cech:1pt}]

 Let $y_0$ be a fixed base point of a component $\huaY_0$ of $\huaC(Y)$. For any $y \in \huaY_0$, there is a unique 1-morphism $  (y_0, y): y \rightarrow y_0$ in $\huaC(Y)_1$. By pseudofunctoriality of $F$, the fiber $F(y)$ is equivalent to $F(y_0)$ via the invertible bimodule to $F(y_0)\xrightarrow{F(y, y_0)} F(y)$. Hence all object data are determined by the values at $y_0$ up to equivalence via the unique $1$-morphisms to the base point.
Moreover, for any $1$-morphism $x \xrightarrow{(y,x)} y$ in $\huaY_0$, $F(y, x)$ is $2$-isomorphic to the composition
\[ F(y, y_0) \circ  F(y_0, x)\] where $F(y_0, x)$ is an inverse of $F(x, y_0)$.

Now let $\iota: F \Rightarrow G$ be a strong transformation over $\huaC(Y)$, and suppose the component $1$-cell $\iota_{y_0}$ is given. For any $y \in \huaY_0$, by the square 
\[ \xymatrix{ F(y_0) \ar[d]_{F(y, y_0)} \ar[r]^{\iota_{y_0}}  &  G(y_0)   \ar[d]^{G(y, y_0)} \\
F(y) \ar@{}[ru]|{ \buildrel{\iota_{(y, y_0)}}\over\Leftarrow}  \ar[r]_{\iota_{y}} &G(y) } \]
the component $\iota_y$ is canonically $2$-isomorphic $ G(y, y_0)   \circ  \iota_{y_0}   \circ  F(y_0, y)$.

Each $\iota_{(y_0, y)}$ is determined by the value of $\iota_{(y, y_0)}$, the unitor of $F$ and $G$,  and the lax unity and lax naturality of strong transformation.
Moreover, for any $1$-morphism $(y, x)$ in $\huaY_0$, the component $2$-cell $\iota_{(y, x)}$ is determined  by the value of  $\iota_{(y, y_0)}$
and $\iota_{(y_0, x)}$, the compositor of $F$ and that of $G$, and the lax naturality of strong transformation.

If two strong transformations $\alpha$ and $ \beta$ agree only up to a fixed invertible modification at the base point of each component, i.e. there is a $2$-isomorphism $\Gamma_{y_0}: \alpha_{y_0}\Rightarrow \beta_{y_0}$ for each base point $y_0$, then, the same parallel-transport construction defines a modification \[ \Gamma: \alpha\Rightarrow \beta\] whose component at each $y$ is \[\Gamma_y := G(y, y_0)\circ \Gamma_{y_0} \circ F(y_0, y). \]
The modification axiom follows by pasting along the \v{C}ech cover and reducing to the basepoint, where it is trivially satisfied.
    
\end{proof}

\begin{proof}[Proof of Theorem \ref{Thm:embed:stringor}]
Because $A$ is an AFD type III$_1$ factor, there exists an increasing sequence of finite-dimensional super algebras
\[ 
A_1 \xrightarrow{\huaH_1} A_2 \xrightarrow{\huaH_2} \cdots \xrightarrow{\huaH_{n-1}} A_n\xrightarrow{\huaH_n}\]
where each $\huaH_n$ is an infinite-dimensional $A_{n+1}-A_n$-bimodule
such that $A$ is the weak-operator closure of the union $\bigcup\limits_n A_n $.
In the bicategory $\mathrm{vNAlg}^{\mathrm{bi}}$, this is equivalently expressed as
\[ A \cong \pcolim_n A_n. \] More explicitly, for each $n\in \mathbb{N}$, there is the triangle of bimodules which commutes up to 2-isomorphism,
\begin{equation} \label{colim:A:legs}
    \xymatrix{ A_n \ar[rr]^{\huaH_n}  \ar[rd]_{\iota_n } && A_{n+1} \ar[ld]^{\iota_{n+1}} \\
    & A \ar@{}[u]|>>>>>{\cong}&  }
\end{equation} and $$\{A, \quad l_n\}_{n\in \mathbb{N}}$$ satisfy the universal property of pseudo-colimit.

For each $n$, define a pseudofunctor $F_n: \huaC(Y)\to \mathrm{vNAlg}^{\mathrm{bi}}$ by:
\begin{itemize}
    \item 
For each object $y\in Y$, $F_n(y) := A_n$;
\item For each arrow $x\xrightarrow{(y, x)} y$ in $\huaC(Y)$, $F_n(y, x) := A_n$;
\item The compositor $F^2_n$ at each pair of composable arrows $\bigl((z, y), (y, x)   \bigr)$ is defined by the left unitor $l_{A_n}$;
\item The unitor $F^0_n := \id$.
\end{itemize}

Moreover, by Lemma \ref{lem:cech:1pt}, the  $A_{n+1}-A_n$-bimodule $\huaH_n$ induces a strong transformation
\[ \alpha_n: F_n \Longrightarrow F_{n+1}\] with the component $1$-cell $(\alpha_n)_{x_0}:= \huaH_n$ at a fixed object $x_0$ of each component of $\huaC(Y)$. In addition, each leg $\iota_n$ in \eqref{colim:A:legs} induces a strong transformation $\beta_n: F_n \Longrightarrow F_{\Stringor}$  with the component $1$-cell $(\beta_n)_{x_0}:= \iota_n$. By Lemma \ref{lem:cech:1pt} and the commutativity up to $2$-isomorphism in \eqref{colim:A:legs}, each triangle 
\begin{equation} \label{colim:stringor:legs}
    \xymatrix{ F_n \ar@{=>}[rr]^{\alpha_n}  \ar@{=>}[rd]_{\beta_n } && F_{n+1} \ar@{=>}[ld]^{\beta_{n+1}} \\
    & F_{\Stringor} \ar@{}[u]|>>>>>{\cong}  &  }
\end{equation} commutes up to an invertible modification. Thus, $\{ \alpha_n; \beta_n \}_{n\in \mathbb{N}}$ defines a cocone with vertex $F_{\Stringor}$. 
Next we show the cocone is indeed a pseudo-colimit cocone.

Let \[ \{\lambda_n: F_n \Rightarrow \huaT\}_{n\in\mathbb{N}}\] be any   cocone under the filtered system $\{ \alpha_n: F_n\Rightarrow F_{n+1}\}_{n\in \mathbb{N}}$. 

For each $x \in Y$, the object-level data of the cocone gives, for every $n$, a 1-morphism
\[ (\lambda_n)_x: A_n  \longrightarrow \huaT(x), \]
i.e. a $\huaT(x)-A_n$-bimodule. Since $A \cong \pcolim_n A_n $ in the bicategory $\mathrm{vNAlg}^{\mathrm{bi}}$, by Lemma \ref{lem:cech:1pt},  the universal property of the pseudo-colimit gives a unique 1-morphism 
\[ \Phi_x: A \longrightarrow \huaT(x) \]
up to equivalence
such that \[ 
\Phi_x \circ (\beta_n)_x \cong (\lambda_n)_x \]
for all $n$.

For each arrow $(y, x): x \to y$ in the \v{C}ech groupoid $\huaC(Y)$, we can define the $2$-isomorphisms 
\[ 
\Phi_{(y,x)}: \huaT(y,x) \circ \Phi_x \longrightarrow \Phi_y \circ F_{\Stringor}(y, x)\] as in Lemma  \ref{lem:cech:1pt}. 
For arbitrary $x$,  $y$, $\Phi_{(x,y)}$ can be defined as in the proof of Lemma \ref{lem:cech:1pt}. So $\Phi$ extends uniquely to a strong transformation.

If $\Phi'$ is another strong transformation satisfying $\Phi' \circ \iota_n \cong  \lambda_n$, then at each object $x$, $\Phi'_x$ and $\Phi_x$ both satisfy the same universal property of the pseudo-colimit $A \cong \pcolim_n A_n$, hence they are equivalent. By Lemma \ref{lem:cech:1pt}, the $1$-morphism components are then also uniquely determined up to equivalence. Therefore $\Phi \cong \Phi'$.

Thus, $F_{\Stringor}$ satisfies the universal property of a pseudo-colimit.

Therefore,
 $F_{\Stringor}$ is the pseudo-colimit of the filtered system
\[
F_0 \buildrel{\alpha_0}\over\Longrightarrow F_1  \buildrel{\alpha_1}\over\Longrightarrow  F_2 \buildrel{\alpha_2}\over\Longrightarrow \cdots .
\]

\end{proof}

\begin{remark}

The proof of Theorem  \ref{Thm:embed:stringor} relies on the fact that  
the fibre $A$ of the stringor bundle is an AFD (equivalently, hyperfinite) type III$_1$ factor.  And the proof crucially uses that the base of the pseudofunctor is over a \v{C}ech groupoid, rather than an arbitrary Lie groupoid, so that Lemma \ref{lem:cech:1pt} applies. 
%    \item Each 1-morphism is an invertible $A-A$-bimodule, given by a twisted standard bimodule $\leftindex_{\omega_x}L^2(A)_{\omega_y}$.
The argument does not generalize to settings where either of these conditions fails.

\end{remark}

\begin{remark}
   Theorem  \ref{Thm:embed:stringor} shows that the stringor bundle, despite its infinite-dimensional fibres, is not an essentially new object in the bicategory of 2-vector bundles: it is the filtered pseudo-colimit of its finite-dimensional approximants. In particular, the hyperfinite type III$_1$ structure of the fibre and the infinite-dimensional nature of the bimodules do not obstruct approximation; they are precisely the type of data that arises from an AFD system. This also gives a rigorous bicategorical formulation of the original heuristic picture of Stolz and Teichner  \cite{st}\cite{StolzTeichnerSpinor}, in which the Stringor bundle was expected to arise as a limit of finite-dimensional geometric data.
\end{remark}

\section{General Calculation} \label{sect:comp_lurie}

In this section, we specialize the 2-equivariant $2K$-theory framework of Section \ref{subsect:def:2eq2k} to the case where the base groupoid is the trivial groupoid  $\pt\git \pt$. In this setting, the classification of 2-equivariant 2-vector bundles reduces to a purely algebraic problem: the 2-representation theory of the acting 2-group $\huaG_{\bullet}$. Whereas the general framework of Section \ref{subsect:def:2eq2k} (and related works such as \cite{KLW2Vect22}) classifies 2-vector bundles over a Lie groupoid via \v{C}ech cohomology, the point-space case eliminates the geometric descent data and leaves only the group cohomological data of the $\huaG_{\bullet}$-action on the fibre. 
The main purpose of the calculations carried out in this section is to provide explicit algebraic counterparts to the 2-equivariant elliptic cohomology table sketched in \cite[Section 5]{Lurie_elliptic_survey}.

\subsection{$2K$-Theory of a Point for  Abelian Lie Groups} \label{subsect_comp_ba}

%Based on the general principle in Example \ref{2rep_concrete}, we compute some concrete examples.

In this subsection we specialize the equivariant framework of Section \ref{subsect:def:2eq2k} to the case where the base groupoid is a point and the acting 2-group is the delooping $BA$ of an abelian Lie group $A$. The resulting representation bicategory $2\Rep(BA)$ is purely algebraic, and its Grothendieck completion $2K_{BA}(\pt)$ can be computed explicitly. We classify the internal equivalence classes of pseudofunctors $B(BA) \rightarrow \text{\tVect}$ and identify the resulting objects and morphisms in the homotopy category. For $A = U(1)$ and $A = \Z/n$, the classification recovers the representation rings $\Z[t,t^{-1}]$ and $\Z[t]/\langle t^n-1 \rangle$, matching the algebraic structures predicted by Lurie's 2-equivariant elliptic cohomology program \cite[Section 5]{Lurie_elliptic_survey}. The computation also illustrates how the Grothendieck completion treats invertible and non-invertible 1-morphisms differently.

Here we let $A$ denote a finitely generated compact abelian Lie group. Then the delooping $BA$ is a strict 2-group with the multiplication functor $m_A$ defined by
\begin{align*}
    m_A: BA \times BA &\longrightarrow BA \\
    (\ast, \ast) &\mapsto \ast \\
    (a, a')&\mapsto a+a'
\end{align*} In other words, the multiplication functor is defined the same as the composition functor. The unit $\ast\git\ast \xrightarrow{u} BA$ is determined by the   only object $\ast$ in $BA$ and the identity morphism $\id_{\ast}$. 

We first give a classification of the internal equivalence classes in the bicategory $2\Rep(BA)$ and compute the 2-equivariant $2K$-theory 
$2K_{BA}(\ast\git \ast)$, which is the Grothendieck completion $K_0(\ho(2\Rep(BA)))$.

Let $F$ denote any pseudofunctor from $BBA $ to  \tVect. We describe its structure explicitly below.
\begin{itemize}
    \item On level $0$, the only object $\ast$ in $BBA$ is sent to a super algebra $\huaD$;
    \item on level $1$, the only $1$-morphism $\id_{\ast}$ is sent to an invertible $\huaD-\huaD$-bimodule $\huaM$;
    \item on level $2$, there is a group homomorphism $\alpha: A\rightarrow \Aut_{\huaD-\huaD}(\huaM)$ (i.e. the group of the invertible $\huaD-\huaD$-intertwiners from $\huaM$ to itself.)
    \item  the compositor and the unitor are given by the $2$-isomorphisms below.
     \[ F^2_{\id_{\ast}, \id_{\ast}}: \huaM\otimes_{\huaD}\huaM\rightarrow \huaM; \quad F^0_{\ast}: \huaD\rightarrow \huaM.\]
     When there is no confusion, we will use the symbols, $F^2$ and $F^0$, to denote the compositor and the unitor respectively, in this section.
     \item The compositor and the unitor satisfy the coherence conditions: lax associativity, lax left and right unity.  
\end{itemize}

In other words, $F$ is determined by the data \[ \bigl(  \huaD, \huaM, F^2, F^0, \alpha  \bigr).\] Note that  
\[ \xymatrix{\huaM\otimes_{\huaD} \huaM \ar[r]^>>>>>{F^2}_>>>>>{\cong} &\huaM \ar[r]^{(F^0)^{-1}}_{\cong} &\huaD }\] is a $2$-isomorphism. Thus, $\huaM$ is an inverse of itself.

\bigskip

In addition, we can define a pseudofunctor $G$ of the form  \[ \bigl(\huaD, \huaD, l_{\huaD}, \id_{\huaD}, \beta \bigr)\] where $\beta$ is any group homomorphism from $A$ to $\Aut_{\huaD-\huaD}(\huaD)$. It's straightforwards to check that the structure $\bigl(\huaD, \huaD, l_{\huaD}, \id_{\huaD}, \beta \bigr)$ satisfy the coherence conditions. %checked, correct.

\begin{lemma}\label{obj_2rep_rep}
    For any object $F$ in the bicategory $2\Rep(BA)$, assuming it is  corresponding to  $\bigl(  \huaD, \huaM, F^2, F^0, \alpha  \bigr)$, there is a unique object $G$ of the form $ \bigl(\huaD, \huaD, l_{\huaD}, \id_{\huaD}, \beta \bigr)$ internal equivalent to $F$.
\end{lemma}
\begin{proof}

We first define below a strong transformation \[ \sigma= (\sigma_{\ast}, \sigma_{\id_{\ast}}):  \quad F\Rightarrow G \] with 
\begin{itemize}
    \item  the component $1$-cell $\sigma_{\ast}:  = \huaM$; 
    \item the component $2$-cell $\sigma_{\id_{\ast}}$ 
     defined by the composition \[\huaD\otimes_{\huaD} \huaM \xrightarrow{l_{\huaM}} \huaM \xrightarrow{(F^2)^{-1}} \huaM\otimes_{\huaD} \huaM .\]
\end{itemize}
The lax unity follows from the lax right unity of $F$; the lax naturality follows from the lax associativity axiom of $F^2$, the lax left unity axiom of $ F^2$ and $F^0$, the naturality of the left unitor $l$, and the naturality of the associator $a$.

In addition, applying the lax left unity of pseudofunctors, we have $F^0\horiprod 1_{\huaM} = (F^2)^{-1}\vertprod l_{\huaM}$. The naturality of $\sigma$ 
\begin{equation} \label{obj_2rep_nat_str}
    \xymatrix{\huaD \otimes_{\huaD} \huaM \ar[d]_{\beta(g)\horiprod 1_{\huaM} }  \ar[rr]^{F^0\horiprod 1_{\huaM}}
    & & \huaM\otimes_{\huaD} \huaM \ar[d]^{1_{\huaM} \horiprod \alpha(g)} \\ 
    \huaD\otimes_{\huaD} \huaM  \ar@{}[rru]|{\circlearrowright} 
    \ar[rr]_{F^0\horiprod 1_{\huaM}}  &&\huaM\otimes_{\huaD} \huaM } ,
\end{equation}    for any $g\in A $,  
uniquely determines the group homomorphism $\beta$ corresponding to the given $\alpha$, thus, the pseudofunctor $G$. %follows from the functoriality of $l_{\huaM}$ and the naturality of $F^2$. 
%More explicitly, for  any $m_1$, $m_2$ in $\huaM$,  \eqref{obj_2rep_nat_str} means  \[ \beta(g) ((F^0)^{-1}(m_1) \cdot m_2) = (F^0)^{-1}(m_1)\cdot\alpha(g)( m_2)\]
%Since $F^0$ is an invertible intertwiner, there exists a unique $m_0\in \huaM$ such that $(F^0)^{-1}(m_0)$ is the unit in $\huaA$. Thus, 
In addition, we obtain directly that \begin{equation}\label{rel_al_be} \beta(g)(1)\cdot m  = \alpha(g)(m)  \end{equation} for any $g\in A$ and $m\in \huaM$.

Next we define below a strong transformation \[ \sigma'= (\sigma'_{\ast}, \sigma'_{\id_{\ast}}): \quad G\Rightarrow F\] with 
 \begin{itemize}
     \item  $\sigma'_{\ast}:  = \huaM$; 
     \item $\sigma'_{\id_{\ast}}$ defined by the composition
        \[\huaM\otimes_{\huaD} \huaM \xrightarrow{F^2} \huaM \xrightarrow{r_{\huaM}^{-1}} \huaM\otimes_{\huaD} \huaD. \]
 \end{itemize}
The lax unity and lax naturality of $\sigma'$ can be checked similarly. % to be checked.
And the naturality of $\sigma'$ gives the same relation between $\alpha$ and $\beta$ as \eqref{rel_al_be}.

The horizontal composition $(\sigma'\horicirc\sigma)$  is canonically $2$-isomorphic to the identity via the   modification $\Gamma: \sigma'\horicirc \sigma \Rightarrow \id_{F}$  defined by the composition
\[\Gamma_{\ast}:  \huaM\otimes_{\huaD} \huaM\xrightarrow{F^2} \huaM \xrightarrow{(F^0)^{-1}} \huaD.\]

The composition $ (\sigma \horicirc \sigma') $ is analogous.
In this way we obtain an internal equivalence between $F$ and $G$.

\end{proof}

\begin{lemma}\label{morita_equiv_2rep}

Let $\huaD$ and $\huaD'$ denote two objects in \tVect \hspace{0.1mm} belonging to the same Morita class. Fix
 a Morita equivalence 
\begin{equation}\label{equiv:dd} \bigl( \huaD\xrightarrow{\huaM} \huaD', \quad \huaD' \xrightarrow{\huaM'} \huaD, \quad 
\huaD \xRightarrow[\cong]{\eta} \huaM'\otimes_{\huaD'} \huaM, \quad 
\huaM\otimes_{\huaD}\huaM' \xRightarrow[\cong]{\eta'} \huaD'
\bigr) \end{equation} between  them. 
For any  $G = \bigl(\huaD, \huaD, l_{\huaD}, \id_{\huaD}, \alpha \bigr)  $, 
  there is a   pseudofunctor $H = \bigl(\huaD', \huaD', l_{\huaD'}, \id_{\huaD'},    \alpha' \bigr) $
  internal equivalent to $G$. 
The homomorphism $\alpha'$ is uniquely determined by $\alpha$ and the chosen Morita equivalence data: explicitly, 
 for any $d'\in \huaD'$ and $m\in \huaM$, $\alpha'(g)(d')\cdot m = (d'\cdot m )\cdot \alpha(g)(1)$. 

Conversely, different choices of Morita equivalence data may yield different $\alpha'$, but the resulting $H$'s are all internally equivalent.
\end{lemma}

\begin{proof}

From the Morita equivalence \eqref{equiv:dd}, we construct   an internal equivalence from $G$ to $H$.

First we define a strong transformation $\sigma: G\Rightarrow H$ with %respect to the direct sum of super algebras.
\begin{itemize}
    \item $\sigma_{\ast} = \huaM$;
    \item $\sigma_{\id_{\ast}} = r^{-1}_{\huaM}\circ l_{\huaM}$.
\end{itemize}

Similarly, we define a strong transformation $\sigma': H\Rightarrow G$ with
\begin{itemize}
    \item $\sigma'_{\ast} = \huaM'$;
    \item $\sigma'_{\id_{\ast}} = r^{-1}_{\huaM'}\circ l_{\huaM'}$.
\end{itemize}

The lax naturality and lax unity of $\sigma$ and $\sigma'$ follows immediately from the unity axiom, \cite[Proposition 2.2.4]{JY:Bicat} and the naturality of left and right unitors. In addition,   the naturality of $\sigma$ and $\sigma'$ 
\begin{equation}
    \xymatrix{  \huaD'\otimes_{\huaD'} \huaM \ar[rr]^{ r^{_1}_{\huaM}\circ l_{\huaM}} \ar[d]_{\alpha'(g)\horiprod 1_{\huaM}} 
    && \huaM\otimes_{\huaD} \huaD \ar[d]^{ 1_{\huaM}\horiprod \alpha(g)} \\
    \huaD'\otimes_{\huaD'} \huaM \ar@{}[rru]|{\circlearrowright}  \ar[rr]_{ r^{_1}_{\huaM}\circ l_{\huaM} } && \huaM\otimes_{\huaD} \huaD},
\end{equation} for any $g\in A$, 
uniquely determines the group homomorphism $\alpha'$ corresponding to the given $\alpha$, thus, the pseudofunctor $H$. Explicitly, for any $d'\in \huaD'$ and $m\in \huaM$, \begin{equation}\alpha'(g)(d')\cdot m = (d'\cdot m )\cdot \alpha(g)(1).\end{equation} Since $l_{\huaM}$ is a $2$-isomorphism, the element $\alpha'(g)(d')$ is uniquely determined. 
%Note that the naturality of $\sigma'$ gives the same formula of $\alpha'$ as above since $\huaM'\otimes_{\huaD'} \huaM$ is isomorphic to $\huaD$.

%\comment{The naturality of $\sigma'$  gives the same relation, only $\huaM$ replaced by $\huaM'$ and $\huaD$, $\huaD'$ interchanged.}

The horizontal composition $(\sigma'\horicirc\sigma)$  is canonically $2$-isomorphic to the identity via   the invertible modification $\Gamma: \id_{\huaD}\longrightarrow \sigma'\horicirc \sigma$ with $\Gamma_{\ast}:= \eta$; 
and the horizontal composition $(\sigma\horicirc\sigma')$  is canonically $2$-isomorphic to the identity via the invertible modification $\Gamma':  \sigma\horicirc \sigma' \longrightarrow \id_{\huaD'}$ with $\Gamma'_{\ast}:= \eta'$.

Thus, we   established an internal equivalence between $G$ and $H$
\[ \bigl( G\buildrel{\sigma}\over\Rightarrow H, \quad 
H\buildrel{\sigma'}\over\Rightarrow G, \quad 
\id_{\huaD}\xRightarrow[\cong]{\Gamma} \sigma'\horicirc \sigma, \quad 
\sigma\horicirc \sigma' \xRightarrow[\cong]{\Gamma'} \id_{\huaD'}\bigr).\]

\end{proof}

\begin{lemma} \label{intert_center}
    The group  $\Aut_{\huaD-\huaD}(\huaD)$ of the $\huaD-\huaD$-intertwiners from $\huaD$ to itself is isomorphic to $Z(\huaD)^{\ast}$, the group of invertible elements in the center of $\huaD$. 
    %\[Z(\huaD)^{\ast}=  \{u\in \huaD\mid u\text{ is invertible; } \hspace{0.1cm} u\cdot a = a\cdot u, \forall a\in \huaD. \} \] 
\end{lemma}
\begin{proof}

    Let $1$ denote the unit in $\huaD$. 
    For any $\phi \in \Aut_{\huaD-\huaD}(\huaD)$, for any $a$, $b$, $x$ in $\huaD$, $\phi(a\cdot x\cdot b) = a\cdot \phi(x) \cdot b$.
    \begin{align*}
        \phi(a) & = \phi(a\cdot 1\cdot 1) = a\cdot \phi(1); \\
        \phi(a) & = \phi(1\cdot 1\cdot a) = \phi(1)\cdot a.
    \end{align*}

    In addition, since $\phi$ is an invertible intertwiner, $\phi(1)$ is an invertible element.
    Thus, there is a one-to-one correspondence between $\Aut_{\huaD-\huaD}(\huaD)$ and $Z(\huaD)^{\ast}$.
    Moreover, for any $\phi$, $\psi$ in $\Aut_{\huaD-\huaD}(\huaD)$, $(\phi\circ\psi)(1) = \phi (1\cdot \psi(1) ) = \phi(1)\psi(1)$. Thus, the correspondence $\phi\mapsto \phi(1)$ defines a group isomorphism.

\end{proof}

\begin{remark} \label{Ba:obj:k}
    Combining Lemma \ref{obj_2rep_rep} and Lemma \ref{morita_equiv_2rep}, we know each internal equivalence class in $2\Rep(BA)$ can be represented by an object of the form   \[  \bigl(\huaD, \huaD, l_{\huaD}, \id_{\huaD}, \alpha \bigr)  \] where $\huaD$ is unique up to Morita equivalence.

    Moreover, if $\huaD=  k $, \[  \bigl(\huaD, \huaD, l_{\huaD}, \id_{\huaD}, \alpha \bigr) \text{    and    } \bigl(\huaD, \huaD, l_{\huaD}, \id_{\huaD}, \alpha' \bigr) \] represents different internal equivalence classes if $\alpha$ and $\alpha'$ are different group homomorphisms. 
    Thus, there is a $1-1$-correspondence between the internal equivalence classes in $2\Rep(BA)$ which sends the only object to  $k$ and 
   \begin{equation}   \hom_{\grp}(A, k^{\ast}) . \end{equation}
\end{remark}

\begin{example}\label{ba:1mor:d}

For any two pseudofunctors $BBA\rightarrow \text{\tVect}$ of the form  \[ F:= \bigl(\huaD, \huaD, l_{\huaD}, \id_{\huaD}, \alpha \bigr)  \quad \text{and} \quad  F':=  \bigl(\huaD', \huaD', l_{\huaD'}, \id_{\huaD'}, \alpha' \bigr) , \]  we give a classification  of the strong transformations $F\Rightarrow F'$. Any strong transformation $\xi: F\Rightarrow F'$ consists of: 
\begin{itemize}
    \item a $D'-D$-bimodule $\xi_{\ast}$;
    \item an invertible $D'-D$-intertwiner $\xi_{\id_{\ast}}: \huaD'\otimes_{\huaD'}\xi_{\ast} \longrightarrow \xi_{\ast}\otimes_{\huaD} \huaD$.
    
\end{itemize}
By lax unity axiom, $\xi_{\id_{\ast}}$ has to be $ r^{-1}_{\xi_{\ast}}\circ l_{\xi_{\ast}}$, which satisfies the lax naturality axiom.  In addition, by the naturality of $\xi$, we have the commutative diagram, for any $a\in A$, 
\[ \xymatrix{\huaD'\otimes_{\huaD'} \xi_{\ast} \ar[d]_{\alpha'(a)\horiprod 1_{\xi_{\ast}}}  \ar[rr]^{r^{-1}_{\xi_{\ast}} \circ l_{\xi_{\ast}}}  && \xi_{\ast}\otimes_{\huaD} \huaD \ar[d]^{1_{\xi_{\ast}} \horiprod \alpha(a)}  \\ \huaD'\otimes_{\huaD'} \xi_{\ast}   \ar[rr]^{r^{-1}_{\xi_{\ast}} \circ l_{\xi_{\ast}}}  && \xi_{\ast}\otimes_{\huaD} \huaD }\]
Thus, the group homomorphism $\alpha$ and $\alpha'$ satisfy the relation
\begin{equation} \label{m_eq_alpha}  \alpha'(a)(1)\cdot m = m\cdot \alpha(a)(1) \quad  \text{  for any  } a\in A, \text{  } m\in \xi_{\ast}. \end{equation} 

Thus, there is a one-to-one correspondence between the $1$-morphisms $F\Rightarrow F'$ in the bicategory  $2\Rep(BA)$  and the   $\huaD'-\huaD$-bimodules $\xi_{\ast}$ satisfying \eqref{m_eq_alpha}.  

\end{example}

In addition, we have seen in the proof of Lemma \ref{morita_equiv_2rep}, if $\huaM$ is invertible, then $D$ and $D'$ are Morita equivalent, and there is a one-to-one correspondence between $\alpha$ and $\alpha'$, i.e. \[ \hom_{\grp}(A, Z(\huaD)^{\ast}) \xrightarrow{\cong} \hom_{\grp}(A, Z(\huaD')^{\ast}) . \]
When we compute  the Grothendieck completion $K_0(\ho(2\Rep(BA)))$, only the $2$-isomorphism classes of  $1$-morphisms are considered. 
In the Grothendieck completion, invertible 1-morphisms between distinct objects are identified at the object level via internal equivalence; thus, apart from the identity morphisms, only non-invertible morphisms between different objects contribute new generators to the morphism groups.
%Thus,   only non-invertible morphisms between different objects in $\ho(2\Rep(BA))$ (other than the identities) are considered in the classification.  
Especially, when $\huaD = k$, such $2$-isomorphism classes corresponds to the dimensions of super $k$-vector spaces.

\begin{example} \label{ex_ba_mor_k}
    In this example we compute the morphisms in  $\ho(2\Rep(BA))$ between objects of the form $(k, \hspace{2mm} \alpha: A\rightarrow k^{\ast})$.

    Let $  (k, \alpha) \longrightarrow (k, \alpha')$ denote a $1$-morphism in $2\Rep(BA)$, which is given by a $k-k$-bimodule, i.e. a super $k$-vector space $V$.  
It should satisfy the condition \eqref{m_eq_alpha}, which implies $\alpha' = \alpha$. Thus, the only $1$-morphism between the two objects in $2\Rep(BA)$
is $0$ if $\alpha \neq \alpha'$.  If $\alpha = \alpha'$, any super $k$-vector space gives a $1$-morphism between them. The $2$-isomorphisms between the $k-k$-bimodules are $k$-linear isomorphisms. Thus, the monoid of the $2$-isomorphism classes of the $1$-morphisms $(k, \alpha)\rightarrow (k, \alpha)$ is isomorphic to $\Z_{\geq 0}$. 

    Thus, those objects in $\ho(2\Rep(BA))$ of the form $(k, \alpha)$ contribute a subcategory of the Grothendieck completion, which is 
    \[ \bigsqcup_{\alpha\in \hom_{\grp}(A, k^{\ast})} B\Z_{\geq 0}.\]
\end{example}

\begin{remark}
    In general, the computation of the whole homotopy category $\ho(2\Rep(BA))$ is nontrivial. We compute the cases for $\huaD = k$ explicitly in Remark \ref{Ba:obj:k} and Example %\ref{ex_ba_k_obj} and 
     \ref{ex_ba_mor_k}. The case $\huaD =k^{\oplus m} $ is deferred Proposition \ref{km_clf} in Section \ref{2k_wr_general}, after
    the computation of the Picard group of $k^{\oplus m}$. 
\end{remark}

%Next we consider the component in  the coproduct $\bigsqcup\limits_{\huaD} \hom_{\grp}(A, Z(\huaD)^{\ast})$ corresponding to the Morita class represented by the field $k$. In this case, $Z(k)^{\ast} = k^{\ast}$. In addition, the factor $\hom_{\grp}(A, k^{\ast})$ is a group induced by the multiplication in $k$. 
%Explicitly, for any group homomorphism $\alpha, \alpha': A\rightarrow k^{\ast}$, $(\alpha \cdot \alpha')(a): =\alpha(a)\cdot \alpha'(a)$. 

Recall that the set of the Morita classes of \tVect \hspace{0.1 mm} is still symmetric monoidal with respect to the direct sum. Especially, we have the relation  \begin{equation}
    \bigl(\huaD, \huaD, l_{\huaD}, \id_{\huaD}, \alpha \bigr) \oplus \bigl(\huaB, \huaB, l_{\huaB}, \id_{\huaB}, \beta \bigr)  
    = \bigl(\huaD \oplus \huaB, \huaD \oplus \huaB, l_{\huaD\oplus \huaB}, \id_{\huaD\oplus \huaB}, \alpha \oplus \beta \bigr) .
\end{equation}
%where for any $g\in A$, $(\alpha\oplus\beta)(g):= \alpha(g)\oplus \beta(g)$. 

\begin{example}\label{ex_ba_k_obj}
%Note $Z(k^{\oplus m})^{\ast} \cong (k^{\ast})^{\oplus m }$. 

For any abelian group $A$, $\hom_{\grp}(A,  (k^{\ast})^{\oplus m})\cong 
\oplus_{1}^m\hom_{\grp}(A,  k^{\ast})$. In this case, the full sub-bicategory $2\Rep_{[k]}(BA)$ of $2\Rep(BA)$ which maps the only object to $k^{\oplus m}$ for some $m\in \Z_{\geq 0}$,  is closed under the direct sum $\oplus_{2\Rep(BA)}$ (See Lemma \ref{2rep_symm_monoidal}). 

When $k$ is the field $\mathbb{C}$, we have the computation below.
\begin{enumerate}
    \item When $A$ is $U(1)$, $\hom_{\grp}(U(1), \mathbb{C}^{\ast})\cong \Z$. Thus, $$\pi_0(2\Rep_{[k]}(BU(1)))\cong \bigoplus_{n\in\Z}
    \Z_{\geq 0}\cdot t^n . $$
    We consider the filtered colimit of the character groups $\hom_{\grp}(U(1), (\mathbb{C}^*)^{\oplus n})$ (with $n\in \Z_{\geq 0}$) under the canonical inclusions. Its Grothendieck group $$K_0(\pi_0(2\Rep_{[k]}(BU(1))))$$ is, up to canonical isomorphism, the free abelian group on the character set of $U(1)$, which recovers $\Z[t, t^{-1}]$. %\footnote{This construction can be understood as a stable equivalence: the inclusions and projections between finite stages become mutually inverse in the colimit. } 

 %   And the Grothendieck group $$K_0(\pi_0(2\Rep_{[k]}(BU(1))))\cong \Z[t, t^{-1}].$$
    \item When $A$ is $\Z/n$,  $\hom_{\grp}(\Z/n, \mathbb{C}^{\ast})$ is isomorphic to the group $\mu_n$ of the $n$-th roots of the unit.  Thus, $\pi_0(2\Rep_{[k]}(B\Z/n))\cong \Z_{\geq 0}[t]/\langle t^n-1 \rangle$. And the Grothendieck group $$K_0(\pi_0(2\Rep_{[k]}(B\Z/n)))\cong \Z[t]/\langle t^n-1 \rangle.$$
\end{enumerate}
\end{example}

%Consider $$N_{A, k } = \hom_{\grp}(A, k^{\ast})\sqcup \{0\}$$ where, for any element $\alpha$ in the first component, $\alpha\cdot 0 := 0$.  $N_{A,k}$ is still a monoid, though not a group anymore.

\subsection{$2K$-Theory of a Point for discrete 2-groups} \label{sect:g:disc}

We now take the acting 2-group to be an ordinary Lie group $G$, viewed as a discrete 2-group with only identity 2-morphisms. In this setting, the representation bicategory $2\Rep(G)$ classifies pseudofunctors $BG \rightarrow \text{\tVect}$. We first classify the internal equivalence classes of such pseudofunctors with fibre an invertible super algebra, showing that they are parametrized by a parity homomorphism $\rho: G \to \Z/2$ and a class $F^2 \in H^2_{\text{sm}}(G;k^{\ast})$. We then classify the strong transformations between them: when the parity homomorphisms agree, these correspond to ordinary or projective super representations of $G$, depending on whether the compositor cohomology classes coincide; when the parities differ, only the zero transformation exists. As a consequence, the endomorphism ring of the trivial object in the Grothendieck completion recovers the ordinary equivariant K-theory $K_G(\ast)$, i.e. the representation ring $R(G)$.

%Let $G$ denote a Lie group. We view it  as a discrete bicategory. %The $G$-equivariant $2K$-theory $2K_{G}(\pt\git\pt)$
 %$K_0(2\Rep(G))$    of the internal equivalence classes in $2\Rep (G)$.

Since each pseudofunctor is equivalent to a normal pseudofunctor, so below we assume all the pseudofunctors are normalized.

\begin{example}
\label{gps_1mor}

In this example we give a classification the level-1 data of pseudofunctors $BG\longrightarrow \text{\tVect}$ under some assumption.

Recall two pseudofunctors  are internally equivalent if and only if they are related by an invertible strong transformation. And a strong transformation is invertible if and only if each of its component 1-cells is invertible. Thus, if two pseudofunctors $F$ and $F'$: $BG\rightarrow \text{\tVect}$ are internally equivalent,  the super algebras $F(\ast)$ and $F'(\ast)$ are Morita equivalent.

 Let $F: BG\longrightarrow \text{\tVect}$ denote a normalized pseudofunctor, which maps the only object $\ast$ to a super algebra $\huaA$ such that $\picard(\huaA)$ is an abelian group. With this condition, for any  $\huaA-\huaA$-bimodule $\huaM$ and any invertible $\huaA'-\huaA$-bimodule $\huaN$, 
    $\huaN\otimes_{\huaA}\huaM\otimes_{\huaA}\huaN^{-1}$ is $2$-isomorphic to $\huaM$, where $\huaN^{-1}$ is an inverse of $\huaN$.

    In this case, for two internally equivalent pseudofunctors $F$ and $F'$: $BG\rightarrow \text{\tVect}$  which map the only object in $BG$ to a super algebra $\huaA$,  
    %if there is an invertible strong transformation $\eta: F\Rightarrow F'$, 
    for each $g\in G$, the invertible $\huaA-\huaA$-bimodules $F(g)$ and $F'(g)$ are $2$-isomorphic. %A lax transformation is invertible if and only if each $\alpha_X$ is invertible and each $\alpha_f$ is invertible. 
    The level-1 data of each internally equivalence class $[F]$ is  determined by a group homomorphism $$\rho: G\longrightarrow \picard(\huaA).$$
    When we record each $F(g)$ only as its equivalence class $\rho(g)$ in $\picard(\huaA)$, we discard the choice of framing for $F(g)$ as a concrete bimodule. Within a single isomorphism class, distinct bimodules differ by a $2$-isomorphism.
This discarded information is not lost. It shifts into the compositor $F^2_{g,h}$, for any $g$, $h\in G$. Explicitly, given a pseudofunctor $F: BG\longrightarrow \text{\tVect}$, for any pseudofunctor $F': BG\longrightarrow \text{\tVect}$ such that, for any $g\in G$, $F'(g) =\psi_g(F(g))$ with each $\psi_g$ an arbitrary $2$-isomorphism, define  $F^{'2}_{g, h}$ by the composition
\begin{equation}
    F'(g)\otimes_{\huaA} F'(h) \xrightarrow{\psi_g^{-1}\horiprod \psi_h^{-1}} F(g) \otimes_{\huaA} F(h) \xrightarrow{F^2_{g, h}} F(gh) \xrightarrow{\psi_{gh}} F'(gh).
\end{equation} The compositors $\{F^{\prime 2}_{g, h}\}_{g, h}$ satisfy the pentagon condition, 
%\[ a_{F'(h), F'(g), F'(f)}= (\id_{F'(h)} \horiprod (F'^2_{g, f})^{-1})     \vertprod (F'^2_{h, gf})^{-1}\vertprod F'^2_{hg, f} \vertprod (F'^2_{h, g}\horiprod \id_{F'(f)}),\] 
which can be proved directly from the pentagon condition of $F^2$ and the naturality of the associator. 

%Thus, in the final classification, what is quotiented out by H^2(G, Z(A)^*) corresponds precisely to the automorphism freedom that was discarded. Nothing is missing.

%In the next example, we  classify the structure 2-morphisms of a given pseudofunctor $BG\rightarrow \text{\tVect}$.

  %  Note,  for any $(g, h)\in G\times G$, $F^2_{g, h}$ is an invertible $k-k$-intertwiner $F(g)\otimes_k F(h)\rightarrow F(gh) $,    % \Aut_{k-k}(F(gh))\cong k^{\ast}$.  
%     which corresponds to an element in $Z(\huaA)^{\ast}$.    Since $F^2$ takes value in $Z(\huaA)^{\ast}$,  

  \end{example}

    \begin{example} \label{ps_equiv_k}

    In this example we give a classification of those pseudofunctors $ BG\longrightarrow \text{\tVect}$ which sends the only object to an invertible super algebra. 

 As indicated in \cite[(1.27)]{Freed_vienna}, for any invertible super algebra $\huaA$,  $\picard(\huaA) \cong \Z/2$. By Lemma \ref{intert_center}, for any invertible $\huaA-\huaA$-bimodule $\huaM$,  the group $\Aut_{\huaA-\huaA}(\huaM)$ is isomorphic to $k^{\ast}$.  
So a normalized pseudofunctor $F: BG\longrightarrow \text{\tVect}$ with $F(\ast)=k $ is classified by $\rho: G \to \Z/2$ and $F^2 \in Z^2_{\text{sm}}(G, k^{\ast})$. 
%    Thus,   the internal equivalence classes of those  pseudofunctors $BG\rightarrow \text{\tVect}$ with $F(\ast) = \huaA $ are classified by  a group homomorphism $\rho: G \to \Z/2$ and $F^2 \in Z^2_{\text{sm}}(G, k^{\ast})$. 
%\[\coprod_{\huaA} \bigl(  \hom_{\grp}(G, \Z/2) \times {H}^2_{\text{sm}}(G, k^{\ast}) \bigr)\] where $\huaA$ goes over the Morita classes of invertible algebras.
The compositor $F^2$ satisfies the 2-cocycle condition by the pentagon condition.

    Two pseudofunctors $F$ and $F'$ are internally equivalent if and only if there exists a family $\{\eta_g \in k^{\ast} \}_{g\in G} $  such that
\[ F'^2_{g,h} = \eta_{gh}^{-1} \cdot F^2_{g,h} \cdot \eta_g  \cdot \eta_h. \] 
This is precisely the condition that $F$ and $F'$ are equivalent up to a coboundary. Hence the equivalence classes of compositors are parametrized by ${H}^2_{\text{sm}}(G, k^{\ast})$.

Thus, the internal equivalence classes of normalised pseudofunctors $F: BG \to \text{\tVect}$  sending the only object to an invertible super algebra   are in one-to-one correspondence with 
$$\hom_{\grp}(G, \Z/2) \times {H}^2_{\text{sm}}(G, k^{\ast}).$$
\end{example}

\begin{example}\label{ps_equiv_k_mor}

For any representatives  \[ F:= (\rho: G\rightarrow \Z/2, \quad F^2) \quad \text{and}  \quad  F':= (\rho': G\rightarrow \Z/2, \quad F'^2) \]
representing two normalized pseudofunctors which map the only object to an invertible super algebra $\huaA$, we give a classification of the morphisms $\alpha: F\Rightarrow F'$ in $\ho(2\Rep(G))$. 

The strong transformation $\alpha$ consists of:
\begin{itemize}
    \item a $\huaA- \huaA$-bimodule $\alpha_{\ast}: F(\ast)\rightarrow F'(\ast)$,  which is of the form $\huaA\otimes_k V$ with $V$ a super $k$-vector space;
    \item for each $g\in G$, an invertible $\huaA-\huaA$-intertwiner $\alpha_g: F'(g)\otimes_{\huaA} \alpha_{\ast}\longrightarrow \alpha_{\ast}\otimes_{\huaA} F(g)$.
\end{itemize}

Since $F(g)\cong \huaA\otimes_k L_g$  and $F'(g)\cong \huaA\otimes_k L'_g$ with $L_g$ and $L'_g$ isomorphic to $ k$ or $\Pi k$, applying $\huaA\otimes_k\huaA \cong \huaA$, the intertwiner $\alpha_g$ reduces to \[ \overline{\alpha}_g: L'_g\otimes_k V\xrightarrow{\cong} V\otimes_k L_g.\] Here $\overline{\alpha}_g$ is an even linear isomorphism of super vector spaces. Hence $\overline{\alpha}_g\in \GL(V)$.

If $\rho\neq \rho'$, the only strong transformation $F\Rightarrow F'$ is that given by $\alpha_{\ast} =\{0\}$ and $\alpha_g = \id_{ \{0\}}$, for any $g\in G$.  Note, the zero transformation with $\alpha_{\ast} = \{0\} $ and each $\alpha_g = \id_{\{0\}}$ always exists whether $\rho = \rho'$ or not. We discuss the cases when $\alpha_{\ast} \neq 0$ below.

When $\rho = \rho'$, $L'_g  = L_g$ as super vector spaces for any $g\in G$. So $\overline{\alpha}_g\in \GL(V)$.  The lax naturality of $\alpha$  yields: 
\begin{equation}
    \alpha_{gf}\circ F'^2_{g, f}= F^2_{g, f}\circ \alpha_g \circ \alpha_f, \quad \text{for any } g, f\in G.
\end{equation} Since $F'^2_{g, f}$  and $ F^2_{g, f} $ are both in $k^{\ast}$ and they 
both act as scalar multiplication, thus, we have 
\begin{equation}
    \alpha_{gf}= (F^2_{g, f}(F'^2_{g, f})^{-1}) \circ \alpha_g  \circ \alpha_f.
\end{equation}

Set \(z_{g,h} := F^2_{g,h} (F'^2_{g,h})^{-1} \in k^{\ast}\); then \(z \in Z^2_{\text{sm}}(G,k^{\ast})\).

\begin{itemize}
    \item  If \(z\) is a coboundary (i.e. \([F^2] = [F'^2]\) in \(H^2_{\text{sm}}(G,k^{\ast})\)),  
    we may adjust \(\alpha_g\) by a 1-cochain to obtain \(\alpha_{gh} = \alpha_g \alpha_h\).  
    Thus \(\alpha\) is an ordinary \emph{super representation} of \(G\) on the super vector space \(V\).  
    Its isomorphism classes correspond to elements of the super representation ring \(s\Rep_k(G)\).

    \item  If \(z\) is not a coboundary, then \(\alpha\) is a \emph{projective super representation} of \(G\) with Schur multiplier \([z] \in H^2_{\text{sm}}(G,k^{\ast})\).  
\end{itemize}

\end{example}

\begin{example}

  %  Especially, when $k=\mathbb{C}$,  the computation of the internal equivalence classes of $G$-representations for some $G$ are given in  Figure \ref{G_U1}. 

% align with the program outlined in \cite{Lurie_elliptic_survey}, where the $2$-representation of $G$ recover the classical representation theory of $G$. 

The results for $k=\mathbb{C}$ are summarized in Figure \ref{G_U1}.
In Figure \ref{G_U1} we use the direct sum $ B (s\Rep_k(G) )\oplus B (s\Rep_k(G)) $ of the delooping to denote the category whose objects are the disjoint union of the object in each delooping, and whose morphisms between distinct summands are the zero morphisms.

\begin{figure}
\begin{center}

\begin{tabular}{|c | c |  c | c |c |  } 
 \hline 
$G$   & $\hom_{\grp}(G, \Z/2) $  &  ${H}^2_{\text{sm}}(G, \mathbb{C}^\ast)$     & Grothendieck group of the objects  \\
\hline \hline 
$\Z/n$ with $n$ odd & $0$ & $0$     &    $ 0 $    \\ \hline 
$\Z/n$ with $n$ even & $\Z/2$ & $0$     &    $ \Z/2 $  \\ \hline 
$U(1)$  & $0$ & $0$  & $0$  \\ \hline 
\end{tabular}     

\bigskip

\bigskip

\begin{tabular}{|c  | c| c| } 
 \hline 
$G$     & the full subcategory of $\ho(2\Rep(G))$  &  Its Grothendieck completion\\
  & consisting of those $F$  with $F(\ast) = k$ & \\ 
\hline \hline 
$\Z/n$ with $n$ odd   &  $ B (s\Rep_k(\Z/n) )$  & $B(K_{\Z/n}(\pt))$  \\ \hline 
$\Z/n$ with $n$ even     & $B (s\Rep_k(\Z/n))\oplus B (s\Rep_k(\Z/n) )  $  &  $B(K_{\Z/n}(\pt))\oplus B(K_{\Z/n}(\pt))$\\ \hline 
$U(1)$   & $B(s\Rep_k(U(1))) $ & $B(K_{U(1)}(\pt))$\\ \hline 
\end{tabular}

\caption{  Examples of the classification for $U(1)$ and $\Z/n\Z$ }  
\label{G_U1}
\end{center}  \end{figure}

\end{example}

Generally, applying the argument in Example \ref{ps_equiv_k} and \ref{ps_equiv_k_mor}, we obtain the conclusion below.
\begin{proposition}

For any Lie group $G$, consider the object $F= (\rho, F^2)$ in $\ho(2\Rep(G))$ which sends the only object $\ast$ in $BG$ to $k$, and sends each $g\in G$ to the $k-k$-bimodule $k$ with $F^2 = l_k\in H^2_{\text{sm}}(G, k^\ast)$. Then the Grothendieck completion of the full subcategory of $\ho(2\Rep(G))$ consisting of this single object $F$ is isomorphic to the delooping $BK_{G}(pt)$ of $G$-equivariant K-theory of a single point. 
    
\end{proposition}

\begin{remark}[Smoothness conditions]

We work with smooth pseudofunctors: all assignments $g \mapsto F(g)$ and compositors $F_{g,h}^2$ are required to be smooth. This is consistent with the standard smooth structure on the target bicategory \tVect (cf. \cite{KLW2Rep22} \cite{KLW2Vect22}). Under this convention, the compositor defines a class in the smooth group cohomology $\mathrm{H}^n_{\text{sm}}(G; k^{\ast})$.

\end{remark} %checked

\subsection{$2K$-theory of a point: the case of $k^{\oplus m}$ and the wreath product} \label{2k_wr_general}

We now generalize the classification of Section \ref{sect:g:disc} from fibre $k$ to fibre $k^{\oplus m}$, the direct sum of $m$ copies of the ground field. This requires a detailed understanding of the Picard group of invertible $k^{\oplus m}-k^{\oplus m}$-bimodules, which we show is isomorphic to the wreath product $(\Z/2) \wr \Sigma_m$ (Proposition \ref{picard_km_wr}). Using this, we classify the internal equivalence classes of pseudofunctors $F: BG \to \text{\tVect}$ with $F(\pt) = k^{\oplus m}$ for $\huaD = k$. The classification is given by orbits of $\Sigma_m$ on $\hom_{\grp}(G, (k^\ast)^{\oplus m})$ (Proposition \ref{km_clf}). This completes the computation of the object-level data in the representation bicategory $2\Rep(G)$ for fibre $k^{\oplus m}$. The non-abelian case is briefly discussed in Remark \ref{rmk_final}, where we note that a full classification would require non-abelian group cohomology.

In this subsection,  we give a  generalization of Section \ref{sect:g:disc}.  Let $m$ be a fixed positive integer.  We classify the internal equivalence classes of pseudofunctors $F: BG \rightarrow \text{\tVect}$ with $F(\ast) = k^{\oplus m}$ and the strong transformations between them.

Let $\huaD$ be a super algebra with $\picard(\huaD)$ abelian. Let $F: BG\longrightarrow \text{\tVect}$ denote a normalized pseudofunctor, which maps the only object $\ast$ to the super algebra $\huaD^{\oplus m }$,  maps each $g\in G $ to an invertible $\huaD^{\oplus m}-\huaD^{\oplus m}$-bimodule $F(g)$.
To simplify the symbols, we use $\huaA$ to denote $\huaD^{\oplus m}$ in this subsection.

\bigskip

We first give a classification of invertible $\huaD^{\oplus m}-\huaD^{\oplus m}$-bimodules.

Let $e_i\in \huaA$ denote the element with the $i$-th component being the unit $1$ and the other components zero. Note $\bigoplus\limits_{i=1}^m e_i$  is the unit of $\huaA$; each $e_i\in Z(\huaA)$; $e_ie_j=0$ if $i\neq j$. 

For any $\huaA-\huaA$-bimodule $\huaM$, let \[ \huaM_{ij} := e_i\cdot \huaM\cdot e_j.\] Thus, $\huaM = \bigoplus\limits_{i, j = 1}^m \huaM_{ij}$. Under this symbol, $\huaA = \bigoplus\limits_{i=1}^m \huaA_{ii}$  with each $\huaA_{ii}=\huaD$. If $\huaM$ is invertible, then each row and each column contains exactly one non-zero summand, each isomorphic to $\huaD$.  Hence,  \begin{equation}\label{km:bimodule:form}\huaM = \bigoplus\limits_{i=1}^m \eta_{\tau(i) \sigma(i)},\end{equation} for some elements $\tau$, $\sigma$ in the symmetric group $\Sigma_m$, where $\eta_{\tau(i)\sigma(i)} \cong \huaD$ as a super $k$-vector space, and  for any $(a_1, \cdots a_m)\in \huaA$ and any $(k_1, \cdots k_m)\in \huaM$ with $k_i \in \eta_{\tau(i)\sigma(i)} $,
\begin{align*}
    (a_1, \cdots a_m)\cdot (k_1, \cdots k_m) & = (a_{\tau(1)}k_1, \cdots a_{\tau(m)}k_m); \\ 
    (k_1, \cdots k_m) \cdot  (a_1, \cdots a_m) & =  (k_1 a_{\sigma(1)}, \cdots k_m a_{\sigma(m)}).
\end{align*}

Each component $\eta_{\tau(i)\sigma(i)}$ is an invertible $\huaD-\huaD$-bimodule, classified by an element $\lambda_i$ in $\picard(\huaD)$, as shown in Example \ref{ps_equiv_k}. Thus, each invertible $\huaA-\huaA$-bimodule $\bigoplus\limits_{i=1}^m \eta_{\tau(i) \sigma(i)}$ is classified by an element $$(\lambda_1, \cdots \lambda_m, \tau, \sigma)$$
in $\picard(\huaD)^{\times m} \times \Sigma_m \times\Sigma_m$.   The correspondence between the collection of invertible $\huaA-\huaA$-bimodules and the set \[ \picard(\huaD)^{\times m}\times \Sigma_m\times\Sigma_m\] is one-to-one.

Let $\alpha\huaM_{\gamma}\beta$ denote the $\huaA-\huaA$-bimodule corresponding to the element \begin{equation}(\lambda_{\gamma(1)}, \cdots \lambda_{\gamma(m)}, \tau\alpha, \sigma\beta). \label{amb_ele_correct}\end{equation}

In fact, $$\underline{\alpha}: \huaM \rightarrow \alpha\huaM_{\alpha}\alpha, \quad (k_1, \cdots k_m)\mapsto (k_{\alpha(1)}, \cdots 
k_{\alpha(m)})$$ defines an invertible $\huaA-\huaA$-intertwiner. It simply permutes the $m$ summands by $\alpha$, without changing the individual bimodule structure on each summand.

Next, we explore a little more the structure of the set  $\picard(\huaD)^{\times m} \times \Sigma_m \times\Sigma_m$. Let  \[ \huaM = \bigoplus\limits_{i=1}^m \eta_{\tau(i) \sigma(i)} \quad \mbox{and} \quad \huaM' = \bigoplus\limits_{i=1}^m \eta'_{\tau'(i) \sigma'(i)} \] denote two invertible $\huaA-\huaA$-bimodules classified by    $$(\lambda_1, \cdots \lambda_m, \tau, \sigma), \quad (\lambda'_1, \cdots \lambda'_m, \tau', \sigma')$$ respectively 
in $\picard(\huaD)^{\times m} \times \Sigma_m \times\Sigma_m$.
The tensor product \begin{align*}
    \huaM\otimes_{\huaA}\huaM' & = \bigl( \bigoplus\limits_{i=1}^m \eta_{\tau(i) \sigma(i)} \bigr) \otimes_{\huaA} \bigl( \bigoplus\limits_{i=1}^m \eta'_{\tau'(i) \sigma'(i)} \bigr) \\
    & =   \bigoplus\limits_{i=1}^m \eta_{\tau(i) \sigma(i)}\otimes_{\huaD}  \eta'_{\tau'\bigl(\tau'^{-1}\sigma(i) \bigr)  \sigma'\bigl(\tau'^{-1}\sigma(i) \bigr)}.
    \end{align*} which is classified by
    \begin{equation} \label{multi_kn_bimodule}
    (\lambda_1, \cdots \lambda_m, \tau, \sigma) \cdot (\lambda'_1, \cdots \lambda'_m, \tau', \sigma')
    : = (\lambda_1\lambda'_{\tau'^{-1}\sigma(1)}, \cdots \lambda_m\lambda'_{\tau'^{-1}
    \sigma(m)}, \tau, \sigma' \tau'^{-1} \sigma ). 
\end{equation}

It is straightforwards to check that the multiplication defined in \eqref{multi_kn_bimodule} is strictly associative.
In addition, the element $$\ukm:= (1, \cdots 1, \id, \id ) , $$ which classifies the weak unit object $\bigoplus\limits_{i=1}^m \huaD $,    satisfies the strict right unitality.
\begin{align*}
   & \ukm\cdot (\lambda_1, \cdots \lambda_m, \tau, \sigma) 
    = (\lambda_{\tau^{-1}(1)}, \cdots \lambda_{\tau^{-1}(m)}, \id, \sigma\tau^{-1})   \xrightarrow[\cong]{\underline{\tau} } 
    (\lambda_1, \cdots \lambda_m, \tau, \sigma). \\
    & (\lambda_1, \cdots \lambda_m, \tau, \sigma)\cdot \ukm 
    =  (\lambda_1, \cdots \lambda_m, \tau, \sigma). 
\end{align*}  

Thus, the sub-bicategory of \tVect \hspace{0.1mm}with only one object $\huaA$, invertible $\huaA-\huaA$-bimodules and the intertwiners between them forms a fair $2$-category in the sense of \cite{kock2006weak}, since the composition is strictly associative and strictly right unital, but the left unit law holds only up to $2$-isomorphisms.

It's also straightforwards to check that the map $\underline{\alpha}$ from $\huaA-\huaA-\text{invBimod}$ to $\huaA-\huaA-\text{invBimod}$ preserves the tensor product and sends the identity $(1,\cdots 1, \id, \id )$ to $(1,\cdots 1, \alpha, \alpha )$.
Note that $\ukm\cdot \ukm=\ukm$. The set  of invertible  $\huaA-\huaA$-bimodules, identified with $\picard(\huaD)^{\times m}\times \Sigma_m\times\Sigma_m$, has the structure of an $H$-space.

Moreover,  any invertible $\huaA-\huaA$-bimodule $\huaM$ identified with an element $ (\lambda_1, \cdots \lambda_m, \tau, \sigma) $ 
is $2$-isomorphic,  via $\underline{\tau^{-1}}$,  to $\tau^{-1}\huaM_{\tau^{-1}}\tau^{-1}$, which corresponds to $(\lambda_{\tau^{-1}(1)}, \cdots \lambda_{\tau^{-1}(m)}, \id, \sigma\tau^{-1} )$.  Thus, 
each invertible $\huaA-\huaA$-bimodule is $2$-isomorphic to one represented by an element $(\lambda_1, \cdots \lambda_m, \tau, \sigma) $ with $\tau = \id$. 

For any $(\lambda_1, \cdots \lambda_m, \id, \sigma) $ and $(\lambda'_1, \cdots \lambda'_m, \id, \sigma') $, \begin{equation}
\label{multi_idtau} 
(\lambda_1, \cdots \lambda_m, \id, \sigma) \cdot (\lambda'_1, \cdots \lambda'_m, \id, \sigma')  =  (\lambda_1\lambda'_{\sigma(1)}, \cdots \lambda_m\lambda'_{\sigma(m)}, \id, \sigma'\sigma). \end{equation}
And $(1, \cdots 1, \id, \id) \cdot (\lambda_1, \cdots \lambda_m, \id, \sigma)  = (\lambda_1, \cdots \lambda_m, \id, \sigma) $. By \eqref{multi_idtau} and the strict associativity and the strict right unitality, the set $$ \{ (\lambda_1, \cdots \lambda_m, \id, \sigma) \mid \mbox{   each   }\lambda_i\in \picard(\huaD); \hspace{0.2mm} \sigma\in \Sigma_m.\} $$
is the wreath product \[ \picard(\huaD)\wr\Sigma_m,\] with the product \eqref{multi_idtau} corresponding to the left action convention of the wreath product.

we can define a surjective $H$-homomorphism $$\Pi$$ from the set of all the invertible $\huaA-\huaA$-bimodules to the wreath product $\picard(\huaD)\wr\Sigma_m$ by sending $(\lambda_1, \cdots \lambda_m, \tau, \sigma) $ to $(\lambda_{\tau^{-1}(1)}, \cdots \lambda_{\tau^{-1}(m)}, \id, \sigma\tau^{-1} )$, which   preserves the multiplication and  the identity.
Thus, $\Pi$ is a strict $H$-homomorphism. And $\ker\Pi = \{(1, \cdots, 1, \alpha, \alpha) \mid \alpha\in \Sigma_m\} $.

Thus, we have proved the conclusion.

\begin{proposition}\label{picard_km_wr}
The Picard group of invertible $\huaD^{\oplus m}- \huaD^{\oplus m}$-bimodules is isomorphic to the wreath product 
    \[ \picard(\huaD) \wr \Sigma_m.\]
\end{proposition}

\begin{proposition}\label{km_clf}
For $\huaD = k$, the internal equivalence classes of pseudofunctors $F: BBA\longrightarrow \text{\tVect}$ with $F(\ast) = k^{\oplus m}$ are in bijection with $\hom_{\grp}(A,  (k^{\ast})^{ \oplus m }) /\Sigma_m$.  \end{proposition} 

\begin{proof}
 By Lemma \ref{obj_2rep_rep}, every class has a normal form $(\huaD, \alpha)$. By Proposition \ref{picard_km_wr}, $\picard(k^{\oplus m}) \cong (\Z/2) \wr \Sigma_m$, and, as implied in \eqref{m_eq_alpha},  the naturality condition of a strong transformation $\sigma: (\huaD, \alpha)\Rightarrow 
 (\huaD', \alpha')$ forces, componentwise,
\[  \alpha'_i(a) \cdot s_i  = s_i \cdot \alpha_{\tau (i)}(a), \quad  \text{ for some }\tau\in \Sigma_m\] where $s=(s_1, \cdots s_m)$ is any element in $\sigma_{\ast}$, which is of the form \eqref{km:bimodule:form}. We can assume there is a nonzero component $s_i$.
Since $s_i$ is non-zero in a $1$-dimensional super $k$-vector space and $k$ is commutative, this implies
\[ \alpha'_i(a) = \alpha_{\tau(i)}(a). \] 
Hence $\alpha' = \tau^{-1} \cdot \alpha$ for some $\tau \in \Sigma_m$. Conversely, every such $\tau$ arises from an invertible bimodule. Therefore the equivalence classes are precisely the orbits of $\Sigma_m$ on $\hom_{\grp}(A,(k^\ast)^{\oplus m})$.     
\end{proof}

\begin{remark}[On the classification of pseudofunctors for non-abelian Picard groups] \label{rmk_final}
For a non-abelian Picard group, the classification of pseudofunctors $F: BG \longrightarrow \text{\tVect}$ with $F(\ast) = \huaD^{\oplus m}$ is more subtle than the product \[ \hom_{\grp}(G, \picard(\huaD)  \wr\Sigma_m ) \times  H^2_{\text{sm}}( G; Z(\huaD)^{\ast})^{\times m} \] suggests. In general,   a complete classification would require an analysis of   non-abelian group cohomology. 
 We do not pursue this general classification here. The computations in this section, in particular the explicit description of $\picard(k^{\oplus m})$, are independent of this issue and will serve as the foundation for future work on the full classification and its applications to power operations in $2K$-theory.
\end{remark}

\section{Related Work and Future Directions} \label{future_direct}

We conclude with two future directions suggested by the framework developed in this paper. These are not needed for the main results, but indicate natural extensions of the theory.

\subsection{Higher transgression}

In this section, we sketch a connection between the 2-categorical framework developed in this paper and the classical notion of transgression in group cohomology. Given a 2-group $\Gamma^\alpha$ arising from a 4-class $\alpha$ on a compact Lie group $G$, its higher inertia groupoid $\pfun(B\Z, B\Gamma^{\alpha})$ encodes the $S^1$-transgression of $\alpha$. We show that the horizontal composition of $1$-morphisms in this higher inertia groupoid recovers the transgression map explicitly. This provides a categorical perspective on transgression and suggests further applications to the $2K$-theory of loop spaces.

Here we start by the basic example below.

\begin{example} \label{skeletal_2grp_recall}
        Let $G$ be a compact Lie group and $\alpha$ an element in $H^4(BG; \Z)$ % \cong H^3_{SM}(G; U(1))$  where $H_{SM}^*(-)$ is the smooth Segal-Mitchison cohomology.  It is shown in \cite[Theorem 99]{schommer:string-finite-dim} that the
%isomorphism classes of the central extensions of the 2-groups \[ 0 \rightarrow \pt \git \mathbb{T} \xrightarrow{f} \Gamma^{\alpha} \xrightarrow{g} G\git G\rightarrow 0 \] 
%are in natural bijection with $H^3_{SM}(G; U(1))$. 
   % The central extension $$\Gamma^{\alpha}$$ corresponding to $\alpha$ is 
   We describe the skeletal 2-group determined by $(G, U(1), \alpha)$  below.
    For more details, see \cite[Section 8.3]{baez:2gp}.
    
    The objects of $\Gamma^{\alpha}$ are $G$ and the morphisms of $\Gamma^{\alpha}$ are the product $G\times U(1)$, with the source and target maps the projection to $G$. The automorphisms of each object are identified with the fiber $U(1)$ as a group. The monoidal structure of $\Gamma^{\alpha}$ is given by the group multiplication in $G$ on objects and the group multiplication of $G\times U(1)$. Although the monoidal structure is indeed strictly associative, we can still equip it with a non-trivial associator defined by $\alpha$, that is, for any $g_1$, $g_2$ and $g_3$ in $G$,
    \[ (g_1g_2)g_3\xrightarrow{ \alpha(g_1, g_2, g_3) \in U(1)} g_1(g_2g_3).\]
    That $\alpha$ is a $3$-cocycle ensures that the pentagon identity is satisfied. The left and right unitors, $l$ and $r$, are uniquely determined by $\alpha$ so that the triangle identities hold. More explicitly, for any elements $g_1$, $g_2$ in $G$, the triangle identity reads
    \begin{equation} \label{alr_rel}
        r_{g_2} + \alpha(g_1, e, g_2) = l_{g_1},
    \end{equation} for any $g_1$, $g_2$ in $G$.

\end{example}

\begin{example}
    Let $G$ be a finite group and $\alpha$ an element in $H^4(BG; \Z) \cong H^3_{SM}(G; U(1))$  where $H_{SM}^*(-)$ is the smooth Segal-Mitchison cohomology. Recall from Example \ref{skeletal_2grp_recall} the skeletal $2$-group $\Gamma^{\alpha}$ determined by $(G, U(1), \alpha)$.  We consider the higher inertia groupoid \[ \pfun(B\Z, B\Gamma^{\alpha})\] 
 whose objects are pseudofunctors from $B\Z$ to $B\Gamma^{\alpha}$,   whose $1$-morphisms are strong transformations and whose $2$-morphisms are modifications. 
    
  An object of $\pfun(B\Z, B\Gamma^{\alpha})$ consists of an element $g\in G$, together with elements 
   $F^{2}_{m, n} $, $F^0  \in U(1)$ satisfying
    \begin{align*}
        F^2_{m,n} + F^2_{m+n, k} & = F^2_{m, n+k} + F^2_{n, k} + \alpha(g^m, g^n, g^k); \\
        F^0+ F^2_{0, m} & = l_{g^m}; \quad 
        F^0 + F^2_{m, 0} =r_{g^m},
    \end{align*} for any $m$, $n$, $k$ in $\Z$.

    A $1$-morphism in $\pfun(B\Z, B\Gamma^{\alpha})$ from $( g, \{F^{2}_{m, n}\}_{m, n}, F^0)$ to 
    $( g', \{F'^{2}_{m, n}\}_{m, n}, F'^0)$ consists of $h\in G$   and elements  $\beta_{ n}\in U(1)$ such that  \[g' = h^{-1} g h, \] and the following cocycle condition holds:
    \begin{equation} \label{lax_nat_inertia_eq}
        \beta_{n+m } + F'^2_{m, n} = F^2_{m, n}+ \alpha(g^{\prime m}, g^{\prime n}, h) + \beta_{ n} - \alpha(g^{\prime m}, h, g^n) + \beta_{ m}  + \alpha (h, g'^m, g'^n). 
    \end{equation}
%In particular, for an automorphism (i.e. $( g, \{F^{2}_{m, n}\}_{m, n}, F^0)=( g', \{F'^{2}_{m, n}\}_{m, n}, F'^0)$), this reduces to  \begin{equation}    \beta_{ n} +  \beta_{ m} =  \beta_{ n+m} -( \alpha(g^m, g^n, h)+\alpha(h, g^m, g^n)-\alpha(g^m, h, g^n)) . \end{equation} 

Now consider the horizontal composition of two $1$-morphisms:
\[(g, F^2, F^0) \xrightarrow{(h, \beta_{ n})} (g', F^{\prime 2}, F^{\prime 0})
\xrightarrow{(h', \beta'_{ n})} (g'', F^{\prime\prime 2}, F^{\prime\prime 0}). \]

A direct computation from the definition of horizontal composition in the bicategory of pseudofunctors \cite[Definition 4.2.15]{JY:Bicat} yields
\begin{equation}\label{tran_hori_comp} (\beta' \horicirc  \beta)_{n}
= (\beta'_{n} + \beta_{ n}) - \tau(\alpha)([h'|h]g_n), \end{equation}
where $\tau: Z^3(G;U(1)) \to Z^2(I(BG);U(1))$ is the $S^1$-transgression map  \cite[Theorem 3]{willerton2008}, and $[h'|h]g_n \in I(BG)$ denotes the loop associated to the triple $(h',h,g^n)$. This proves Theorem \ref{transgression_inertia}.

 \label{higher_inertia_string}
\end{example}

\begin{theorem}
    \label{transgression_inertia}
    Let $G$ be a finite group and $\alpha$ an element in $H^3_{SM}(G; U(1))$. Let $\Gamma^{\alpha}$ denote the skeletal 2-group determined by $(G, U(1), \alpha)$. Then the image of $\alpha$ under the $S^1$-transgression \[ \tau: Z^3(G; U(1))\longrightarrow Z^2( I(BG); U(1))\] is determined by the horizontal composition of the 1-morphisms, which are the strong transformations, in the higher inertia groupoid \[ \pfun(B\Z, B\Gamma^{\alpha}).\]
This point is expressed explicitly in \eqref{tran_hori_comp}.
\end{theorem}

\begin{remark}

The discussion in Example \ref{higher_inertia_string}  and Theorem \ref{transgression_inertia} holds verbatim if $U(1)$ is replaced by any abelian Lie group $H$. In that case, one considers the skeletal 2-group $\Gamma_H^\beta$ determined by a class $\beta \in H^3_{SM}(G;H)$, and the transgression map becomes
\[  \tau_H: Z^3(G; H)\longrightarrow Z^2( I(BG); H) \]
The horizontal composition formula \eqref{tran_hori_comp} then takes values in $H$ instead of $U(1)$.

\end{remark}

\begin{example}[Categorical Transgression over a point, sketch]
A natural next step is to apply the transgression formula of Theorem \ref{transgression_inertia} within the $2K$-theoretic framework. Given a $\Gamma^\alpha$-equivariant 2-vector bundle   $\Phi: B\Gamma^\alpha \rightarrow \text{\tVect}$ over a point, one may precompose with a pseudofunctor $F: B\Z \to B\Gamma^\alpha$ to obtain a 2-vector bundle on the higher inertia groupoid. The bicategorical equalizer of $\Phi \circ F(1)$ with the constant functor at the fibre is expected to produce a 2-vector bundle over $\pfun(B\Z, B\Gamma^{\alpha})$, whose product structure is governed by the transgression $\tau(\alpha)$. The details, including the precise construction of the equalizer and the verification of the coherence conditions, will be addressed in future work.
     
\end{example}

\textbf{Future directions: twisted Tate K-theory} The transgression picture outlined above is closely related to the twisted Tate K-theory of loop spaces.  A forthcoming project will establish a precise relation between the $\String(G)$-equivariant $2K$-theory of this loop groupoid and the twisted $G$-equivariant Tate K-theory of the free loop space, as anticipated in the introduction. This would provide a direct link between the higher categorical framework developed here and the geometry of loop spaces and elliptic cohomology.

\subsection{2-orbifold 2-vector bundles} \label{sect_orb_vb_2K}

In this section, we outline a framework for 2-orbifold 2-vector bundles, extending the theory of Section \ref{sect_2vect_2K} to the setting of weak groupoid objects internal to a bicategory. This provides a higher-categorical analogue of orbifold K-theory, where the usual orbifold groupoid is replaced by a weak groupoid object in the bicategory $\bgpd$ of Lie groupoids, bibundles, and equivariant bibundle maps. The resulting 2-orbifold $2K$-theory simultaneously generalizes the ordinary $2K$-theory of Lie groupoids (Section \ref{sect_2vect_2K}) and the $2$-equivariant $2K$-theory of Lie 2-group actions (Section \ref{sect_2eq_2vb_2k}).

The model we take for 2-orbifolds is a Segalic pseudofunctor \cite[Definition 8.1.2]{Paoli2019}, which we sketch below. 

\begin{definition}
    A weak groupoid object in a bicategory $\huaC$ (with bi-pullbacks) is a pseudofunctor
    \[ \huaX: \Delta^{op}\longrightarrow \huaC\] such that \begin{enumerate}
    \item $\huaX_0$ is discrete;
        \item for each $n\geq 2$, the Segal map $\huaX_n\longrightarrow \huaX_1\times_{\huaX_0}\cdots\times_{\huaX_0} \huaX_1$ is an equivalence in $\huaC$;  (Segal condition);
        \item the weak globularity conditions of \cite[Definition 8.1.2]{Paoli2019} hold in $\huaC$: the structure maps induced by the degeneracy maps of $\Delta$ are equivalences.

        \item the groupoid condition: there exists an inversion 1-morphism $\mathrm{inv}: \huaX_1 \to \huaX_1$ with 2-isomorphisms
  \begin{align*}
        \mu \circ (\mathrm{inv} \times\id_{\huaX_1}) & \cong u \circ \source, \\    
      \mu \circ (\id_{\huaX_1} \times \mathrm{inv}) & \cong u \circ \target.   \end{align*}  
    \end{enumerate}
\end{definition}

\begin{construction}

Let $\huaX:\Delta^{op}\to\huaC$ be a weak groupoid object. Define a bicategory $\mathbb{X}$ as follows. Here $\huaX_0$ and $\huaX_1$ are Lie groupoids.
\begin{itemize}
    \item  objects:  $\huaX_0$;
\item the category of 1-morphisms:  $\huaX_1$;
\item  composition of 1-morphisms: induced by the Segal map $\huaX_1\times_{\huaX_0}\huaX_1 \xrightarrow{\simeq} \huaX_2$;
\item  identity 1-morphisms: induced by the degeneracy $s_0:\huaX_0\to\huaX_1$;
\item  horizontal composition of 2-morphisms: induced by the face maps of $\huaX_2$;
\item  the associators and unitors are given by the pseudofunctoriality of $\huaX$.

\end{itemize}

We will   refer to the bicategory   $\mathbb{X}$ as \textit{the bicategory associated to the weak groupoid object $\huaX$}.

\end{construction}

\bigskip

\begin{definition}[2-orbifold 2-vector bundles] \label{def:2orb2vectbl}
Let $\huaX$ be a weak groupoid object in $\bgpd$, and let $\mathbb{X}$ be the associated bicategory. 
A  2-orbifold 2-vector bundle over $\huaX$ is a pseudofunctor
\[
\huaV: \mathbb{X} \longrightarrow \text{\tVect}.
\]

A  1-morphism between two 2-orbifold 2-vector bundles $\huaV, \huaW$ is a strong transformation
\[
\alpha: \huaV \Rightarrow \huaW.
\]

A  2-morphism between two 1-morphisms $\alpha, \beta$ is a modification
\[
\Gamma: \alpha \Rightarrow \beta.
\]

The resulting bicategory is denoted by
\[
(\textsf{s}2\huaV \huaB dl_k)_{\mathrm{orb}}(\huaX).
\]
\end{definition}

\begin{definition}[2-orbifold 2K-theory] \label{def:2orb_2k}
The  2-orbifold 2K-theory of $\huaX$ is the Grothendieck completion (Definition \ref{def:gr_comp}) of a skeleton of the homotopy category of $(\textsf{s}2\huaV \huaB dl_k)_{\mathrm{orb}}(\huaX)$:
\[
2K_{\mathrm{orb}}(\huaX) := K_0\left(\ho((\textsf{s}2\huaV \huaB dl_k)_{\mathrm{orb}}(\huaX))\right).
\]
\end{definition}

\begin{example}[Lie groupoids]
Let $\huaX_\bullet$ be a Lie groupoid. It gives a weak groupoid object $\huaX$ in $\bgpd$ by taking:
\[
\huaX_0 := \huaX_0\git  \huaX_0, \qquad
\huaX_1 := \huaX_1\git  \huaX_1,
\] viewed as discrete Lie groupoids,
with source and target induced by $\source,\target:\huaX_1\to\huaX_0$, composition induced by the groupoid composition, and all structure 2-morphisms identities.

The associated bicategory $\mathbb{X}$ has:
\begin{itemize}
    \item objects:  $\huaX_0$;
    \item 1-morphisms:  $\huaX_1$;
    \item 2-morphisms: identity 2-morphisms.
\end{itemize}
Thus $\mathbb{X}$ is the ordinary Lie groupoid $\huaX_\bullet$ regarded as a bicategory with only identity 2-morphisms. 

Consequently, a 2-orbifold vector bundle over $\huaX$ is precisely a pseudofunctor $\huaV:\huaX_\bullet\to\text{\tVect}$, i.e., a 2-vector bundle in the sense of Definition \ref{def:2vect:obj}.
\end{example}

\begin{example}[Equivariant case]
Let $\huaG_\bullet$ be a coherent Lie 2-group acting on a Lie groupoid $\huaY_\bullet$. The action groupoid 
\[
\huaG_\bullet \ltimes \huaY_\bullet
\]
is a weak groupoid object in $\bgpd  $, with
\[
\huaX_0 := \huaY_\bullet, \qquad \huaX_1 := \huaG_\bullet \times \huaY_\bullet,
\]
source the projection to $\huaY_\bullet$, target the action map, and composition induced by the 2-group multiplication.

The associated bicategory $\mathbb{X}$ has:
\begin{itemize}
    \item objects:   $\huaY_\bullet$;
    \item 1-morphisms:   $\huaG_\bullet   \times \huaY_\bullet$;
\end{itemize}

A 2-orbifold 2-vector bundle over $\huaX$ is precisely a pseudofunctor
\[
\huaV: \mathbb{X} \longrightarrow \text{\tVect},
\]
which is exactly a $\huaG_\bullet$-equivariant 2-vector bundle over $\huaY_\bullet$ in the sense of Definition \ref{def:2eq2vectgrpd}.
\end{example}

\begin{remark}
The two examples above show that the notion of 2-orbifold 2-vector bundle unifies:
\begin{enumerate}
    \item \textbf{Ordinary 2-vector bundles}: when $\huaX$ comes from a Lie groupoid;
    \item \textbf{2-Equivariant 2-vector bundles}: when $\huaX$ is a Lie groupoid acted by a coherent Lie 2-group.
\end{enumerate}
Thus the 2-orbifold $2K$-theory $2K_{\mathrm{orb}}(\huaX)$ simultaneously generalizes both $2K(\huaY_\bullet)$ and $2K_{\huaG_\bullet}(\huaY_\bullet)$.
\end{remark}

\textbf{Future directions} The framework outlined above opens several avenues for further investigation. A natural next step is to develop the analytic underpinnings of 2-orbifold 2-vector bundles, including notions of sections, metrics, and Dirac operators on 2-orbifolds, in analogy with the classical orbifold setting. On the topological side, the 2-orbifold $2K$-theory defined here is expected to admit a proper equivariant spectrum lifting, along the lines of Appendix \ref{subsect:2k_spectrum}, which would provide a richer invariant than the Grothendieck completion of the homotopy category. Finally, the relation between 2-orbifold 2K-theory and the stringor bundle of Section \ref{sect:stringor}, both of which involve higher geometric structures on loop spaces, deserves a systematic study; it is plausible that the 2-orbifold framework provides a natural home for the stringor bundle when regarded over a suitable loop groupoid.
\newpage

\appendix

\section{A K-theory spectrum} \label{subsect:2k_spectrum}

In this section we formulate a K-theory spectrum.
The idea is  based on \cite{Osorno2012SpectraAGT}.

 We first define the bicategory of $\huaG_{\bullet}$-representations, %in the spirit of \cite{KLW2Rep22} and \cite{Huan:2Rep_2Vect}. 
 $2\Rep \huaG_{\bullet}$, which is the bicategory of pseudofunctors from $B\huaG_{\bullet}$ to \tVect, strong transformations  and modifications. 

\begin{lemma} \label{2rep_symm_monoidal}
    The bicategory $2\Rep \huaG_{\bullet}$ is strict symmetric monoidal in the sense of \cite[Definition 1.3]{Osorno2012SpectraAGT}.
\end{lemma}

\begin{proof}

We can define a strict functor $\oplus_{2\Rep \huaG_{\bullet}}: 2\Rep \huaG_{\bullet}\times 2\Rep \huaG_{\bullet} \longrightarrow 2\Rep \huaG_{\bullet}$, which is defined from the strict functor (see Lemma \ref{2vect_strict}) \[ \oplus: \text{\tVect} \times \text{\tVect} \longrightarrow \text{\tVect}.\]
For each pair of pseudofunctors $\mathbb{F}:= (F, F^2, F^0)$ and  $\mathbb{G}:=(G, G^2, G^0)$ from $B\huaG_{\bullet}$ to \tVect, 
$\mathbb{F}\oplus_{2\Rep \huaG_{\bullet}}\mathbb{G}$ is defined to be the composition 
\begin{equation}
    B\huaG_{\bullet} \xrightarrow{\Delta} B\huaG_{\bullet}\times B\huaG_{\bullet} \xrightarrow{\mathbb{F}\times \mathbb{G}} \text{\tVect} \times \text{\tVect} \xrightarrow{\oplus} \text{\tVect}.
\end{equation} In the composition $\Delta$ is the diagonal functor, which is strict; and the pseudofunctor \[\mathbb{F}\times \mathbb{G} : = \Bigl(
F\times G, \quad F^2\times G^2, \quad F^0\times G^0 
\Bigr).  \]

A $1$-morphism in $2\Rep \huaG_{\bullet}\times 2\Rep \huaG_{\bullet} $ is given by a pair of strong transformations. Let $\mathbb{F}, \mathbb{G}, \mathbb{F}', \mathbb{G}'$ denote four pseudofunctors from $B\huaG_{\bullet}$ to \tVect, and let $\alpha: \mathbb{F}\Rightarrow \mathbb{G}$ and $\alpha': \mathbb{F}'\Rightarrow \mathbb{G}'$ denote two strong transformations. Then $\alpha\oplus_{2\Rep \huaG_{\bullet}}  \alpha': \mathbb{F}\oplus_{2\Rep \huaG_{\bullet}} \mathbb{G} \Rightarrow \mathbb{F}'\oplus_{2\Rep \huaG_{\bullet}} \mathbb{G}'$ is defined to be the horizontal product
\begin{equation}
    \id_{\Delta}\horicirc (\alpha \times \alpha') \horicirc \id_{\oplus}. 
\end{equation} We make it explicitly  what  this horizontal product means below. 
For the single object $\ast$, $(\alpha\oplus_{2\Rep \huaG_{\bullet}}  \alpha')_{\ast}$ is defined to be the direct sum \[ F(\ast)\oplus F'(\ast)
\xrightarrow{\alpha_{\ast}\oplus \alpha'_{\ast}} G(\ast)\oplus G'(\ast)\] 
For each object $b_0$ in $\huaG_{\bullet}$, $$(\alpha\oplus_{2\Rep \huaG_{\bullet}} \alpha')_{b_0}:= 
\alpha_{b_0}\oplus \alpha'_{b_0}: \quad \bigl( G(b_0) \oplus G'(b_0) \bigr) \circ (\alpha_{\ast} \oplus \alpha'_{\ast}) \longrightarrow 
(\alpha_{\ast}\oplus \alpha'_{\ast}) \circ \bigl( F(b_0) \oplus F'(b_0) \bigr)$$
Note that \begin{align*}\bigl( G(b_0) \oplus G'(b_0) \bigr) \circ (\alpha_{\ast} \oplus \alpha'_{\ast}) &= (G(b_0)\circ \alpha_{\ast})\oplus (G'(b_0) \circ \alpha'_{\ast}); \\ (\alpha_{\ast}\oplus \alpha'_{\ast}) \circ \bigl( F(b_0) \oplus F'(b_0) \bigr)
&= (\alpha_{\ast}\circ F(b_0) ) \oplus  ( \alpha'_{\ast} \circ F'(b_0)).\end{align*} 

The lax unity and lax naturality of $\alpha\oplus_{2\Rep \huaG_{\bullet}}  \alpha'$  can be obtained immediately from those of $\alpha$ and $\alpha'$. 

Since the cartesian product $\times$ is strictly associative and unital, and horizontal product preserves vertical product, %(\cite[Explanation 2.1.6, (2.1.9)]{JY:Bicat}), 
$\oplus_{2\Rep \huaG_{\bullet}}$ is strictly associative and unital, thus, a strict functor. 

\bigskip

The zero object in $2\Rep \huaG_{\bullet}$ is the the trivial pseudofunctor, which sends $\ast$ to the super algebra $0$, sends each object  of $\huaG_{\bullet}$    to %the $0-0$-bimodule $0$, i.e. 
$\id_0$, and sends each $1$-morphism in $\huaG_{\bullet}$ to the identity intertwiner between the bimodule $0$.

Let $\tau:2\Rep \huaG_{\bullet} \times 2\Rep \huaG_{\bullet} \longrightarrow 2\Rep \huaG_{\bullet}\times 2\Rep \huaG_{\bullet} $ denote the switch pseudofunctor. We can define a strong transformation $\varpi: \oplus_{2\Rep \huaG_{\bullet}} \longrightarrow \oplus_{2\Rep \huaG_{\bullet}} \circ \tau$ in the way below.
For each pair of objects $\mathbb{F}$ and $\mathbb{G}$ in $ 2\Rep \huaG_{\bullet}$, $$\varpi_{(\mathbb{F}, \mathbb{G})}:= \tau_{(\mathbb{F}, \mathbb{G})}: F(\ast)\oplus G(\ast) \longrightarrow G(\ast)\oplus F(\ast).$$ And for each pair of 1-morphisms $\alpha: \mathbb{F}\Rightarrow \mathbb{G}$ and $\alpha': \mathbb{F}'\Rightarrow \mathbb{G}'$ in $2\Rep \huaG_{\bullet}$, $$\varpi_{(\alpha, \alpha')}: (\alpha'\oplus_{2\Rep \huaG_{\bullet}} \alpha) \circ \tau_{(\mathbb{F}, \mathbb{F}')} \longrightarrow \tau_{(\mathbb{G}', \mathbb{G})} \circ (\alpha\oplus_{2\Rep \huaG_{\bullet}} \alpha')$$ is defined to be the identity. It's straightforward to check that $\varpi$ satisfy the strict diagrams in \cite[Definition 1.3]{Osorno2012SpectraAGT}, i.e. \eqref{strsymmvarpi1} and \eqref{strsymmvarpi2} below.

\begin{equation}\label{strsymmvarpi1}
    \xymatrix{ 2\Rep \huaG_{\bullet}\oplus_{2\Rep \huaG_{\bullet}} 2\Rep \huaG_{\bullet} \oplus_{2\Rep \huaG_{\bullet}} 2\Rep \huaG_{\bullet} \ar[rr]^{\varpi\oplus_{2\Rep \huaG_{\bullet}} \id }\ar[rrd]_{\varpi_{1, 2\oplus 3}}
     && 2\Rep \huaG_{\bullet}\oplus_{2\Rep \huaG_{\bullet}} 2\Rep \huaG_{\bullet}\oplus_{2\Rep \huaG_{\bullet}} 2\Rep \huaG_{\bullet}  \ar[d]^{ \id\oplus_{2\Rep \huaG_{\bullet}} \varpi} \\
     && 2\Rep \huaG_{\bullet} \oplus_{2\Rep \huaG_{\bullet}}  2\Rep \huaG_{\bullet} \oplus_{2\Rep \huaG_{\bullet}} 2\Rep \huaG_{\bullet}   }
\end{equation}

    \begin{equation} \label{strsymmvarpi2}
       \xymatrix{    2\Rep \huaG_{\bullet}\oplus_{2\Rep \huaG_{\bullet}}   2\Rep \huaG_{\bullet}  \ar@{=}[rr] \ar[rd]_{\varpi} && 2\Rep \huaG_{\bullet}\oplus_{2\Rep \huaG_{\bullet}}  2\Rep \huaG_{\bullet} \\ & 2\Rep \huaG_{\bullet}\oplus_{2\Rep \huaG_{\bullet}} 2\Rep \huaG_{\bullet} \ar[ru]_{ \varpi}&}
    \end{equation}

    This also shows that the direct sum is well-defined on internal equivalence classes: 
if $F \simeq F'$ and $G \simeq G'$, then $F \oplus G \simeq F' \oplus G'.$ 
\end{proof}

%Let $M(2\Rep \huaG_{\bullet})$ denote the sub-bicategory of $2\Rep \huaG_{\bullet}$ with  the same objects as $2\Rep \huaG_{\bullet}$, same $1$-morphisms as $2\Rep \huaG_{\bullet}$ and invertible $2$-morphisms. 
 By \cite[Theorem 2.1]{Osorno2012SpectraAGT},  we have an infinite loop space upon group completion. 
Similar to the construction in Section \ref{2K_spectra}, applying the construction of classifying space of bicategory, we can  define the  $2K$-theory of $2\Rep(\huaG_{\bullet})$. %w.r.t the 2-vector space model in \cite{KLW2Vect22}. 
\begin{definition}[K-theory space]
    \[\huaK(2\Rep(\huaG_{\bullet})) : = \Omega B (|N^D_{\bullet} (M(2\Rep (\huaG_{\bullet}) ))|)  . \]
        
\end{definition}

\begin{definition}[K-theory spectrum]
%Moreover, by \cite[Theorem 2.1]{Osorno2012SpectraAGT}, there exists a spectrum $\mathbb{A} ( M(2\Rep \huaG_{\bullet}) )$. 
The  K-theory spectrum of $2\Rep (\huaG_{\bullet})$ is defined to be \[ \mathbb{K}(2\Rep (\huaG_{\bullet})) : = \mathbb{A}( M(2\Rep (\huaG_{\bullet})) ). \] By \cite[Theorem 2.1]{Osorno2012SpectraAGT}, its zeroth space of this spectrum is $\huaK(2\Rep(\huaG_{\bullet}))$.  \end{definition}

\begin{remark}
This spectrum is the equivariant analogue of the non-equivariant spectrum constructed in Section \ref{2K_spectra}. Indeed, when $\huaG_\bullet$ is the trivial 2-group, $2\Rep(\huaG_\bullet)$ reduces to \tVect, and $\mathbb{K}(2\Rep(\huaG_\bullet))$ agrees with the spectrum $\mathbb{K}(\text{\tVect})$ from Section \ref{2K_spectra}. Moreover, the zeroth homotopy group of this spectrum recovers the Grothendieck group of objects of $\ho(2\Rep(\huaG_\bullet))$, i.e. the object set of the $\huaG_{\bullet}$-equivariant $2K$-theory $2K_{\huaG_\bullet}(\ast)$ of a point  from Example \ref{2rep_concrete} .
\end{remark}

\newpage

\subsection*{Acknowledgment}
This research is based upon the work  supported by the Young Scientists Fund of the National Natural Science Foundation of China (Grant No. 11901591), the General Program of the National Natural Science Foundation of China (Grant No. 12371068)  for the project 
``Quasi-elliptic cohomology and its application in topology and mathematical physics", the National Science Foundation under Grant Number DMS 1641020,  and the research funding from Huazhong University of Science and Technology.

%The research was inspired during the Mathematical Research Community (MRC) workshop at Rhode Island in June 2019. The author feels indebted to Dileep Menon, Matthew Bruce Young for valuable supports and advice, especially the organizers Nora Ganter and Daniel Berwick-Evans.
%The author thanks Hisham Sati and Urs Schreiber for suggesting the author compute twisted quasi-elliptic cohomology of spheres, and thanks Center for Quantum and Topological Systems at New York University Abu Dhabi for hospitality and support. In addition, 
The author thanks the Max Planck Institute for Mathematics and   Beijing International Center for Mathematical Research for hospitality and support. Part of this work was done during the author's visit at MPIM and  BICMR. 

The author would like to thank Paolo Tomasini for helpful discussion on elliptic cohomology and higher geometry at MPIM early in 2024. 
The author also thank Konrad Waldorf for helpful discussions. In addition, the author would like to  thank Niles Johnson and Donald Yau for their book \textit{2-Dimensional Categories}, which is exceptionally thorough and crystal clear.

\newpage
\bibliographystyle{alpha}
\bibliography{2equiv2K}
\end{document}